\documentclass[12pt,reqno]{amsart}

\usepackage{setspace}

\usepackage{enumitem}
\usepackage{tikz-cd}
\usepackage{tikz}
\usetikzlibrary{positioning} 

\usepackage[linktocpage=true]{hyperref}
\makeatletter
\@namedef{subjclassname@2020}{%
  \textup{2020} Mathematics Subject Classification}
\makeatother
\newtheorem*{theorem*}{Theorem A}
\newtheorem*{theorem**}{Theorem B}

\usepackage[T1]{fontenc}
\usepackage{ae,aecompl}
\usepackage{graphics,cancel,graphicx}
\usepackage{epsfig,psfrag,manfnt}
\usepackage{mathrsfs}
\usepackage{graphicx,mathabx}
\usepackage{bbm,subfigure,rotating,bm,float,mathdots,wasysym,scalerel}

\newlength{\wdth}

\usepackage{stmaryrd}
\usepackage[all,cmtip]{xy}
\usepackage{amssymb,amsmath,amsfonts,amsthm,color}

\usepackage{scalerel,stackengine}
\stackMath
\newcommand\reallywidehat[1]{%
\savestack{\tmpbox}{\stretchto{%
  \scaleto{%
    \scalerel*[\widthof{\ensuremath{#1}}]{\kern-.6pt\bigwedge\kern-.6pt}%
    {\rule[-\textheight/2]{1ex}{\textheight}}%WIDTH-LIMITED BIG WEDGE
  }{\textheight}% 
}{0.5ex}}%
\stackon[1pt]{#1}{\tmpbox}%
}
\newcommand{\calA}{\mathcal{A}}

\newcommand{\calO}{\mathcal{O}}
\newcommand{\calP}{\mathcal{P}}

\newcommand{\mA}{\mathbb{A}}
\newcommand{\mB}{\mathbb{B}}
\newcommand{\mC}{\mathbb{C}}
\newcommand{\mD}{\mathbb{D}}

\newcommand{\mF}{\mathbb{F}}

\newcommand{\mH}{\mathbb{H}}

\newcommand{\mN}{\mathbb{N}}
\newcommand{\mO}{\mathbb{O}}
\newcommand{\mP}{\mathbb{P}}

\newcommand{\mR}{\mathbb{R}}

\newcommand{\mT}{\mathbb{T}}

\newcommand{\mZ}{\mathbb{Z}}

\newcommand{\bba}{\bm{a}}

\newcommand{\bbh}{\bm{h}}

\newcommand{\bbp}{\bm{p}}

\newcommand{\bbw}{\bm{w}}

\newcommand{\bbz}{\bm{z}}

\newcommand{\com}{\complement}
\newcommand{\balpha}{\bm{\alpha}}

\newcommand{\bcdot}{{\scaleobj{1.2}{\bm{.}}}}

\newtheorem{theorem}{Theorem}[section]
\newtheorem{lemma}[theorem]{Lemma}
\newtheorem{corollary}[theorem]{Corollary}
\newtheorem{proposition}[theorem]{Proposition}

\theoremstyle{definition}

\newtheorem{remark}[theorem]{Remark}
\newtheorem*{Th}{Theorem}

\theoremstyle{definition}

\theoremstyle{definition}

\theoremstyle{definition}

\begin{document}



\title[]{On algebras of Dirichlet series invariant under permutations  of coefficients}

\author[]{Alexander Brudnyi}
\address{Department of Mathematics and Statistics\\
University of Calgary\\
Calgary, Alberta, Canada T2N 1N4}
\email{abrudnyi@ucalgary.ca}

\author[]{Amol Sasane}
\address{Department of Mathematics \\London School of Economics\\
     Houghton Street\\ London WC2A 2AE\\ United Kingdom}
\email{A.J.Sasane@lse.ac.uk}
 
\maketitle
  
\begin{abstract} 
Let $\mathscr O_u$ be the algebraof holomorphic functions on 
$\mC_+:=\{s\in \mC: \textrm{Re }s>0\}$ that are limits of Dirichlet series 
$D=\sum_{n=1}^\infty a_n n^{-s}$, $s\in \mC_+$, that converge uniformly 
on proper half-planes of $\mC_+$. We study algebraic-topological properties 
of natural topological subalgebras of $\mathscr O_u$, the Banach algebras 
$\mathscr W, \mathscr A, \mathscr H^\infty$ and the Fr\'echet algebra 
$\mathscr O_b$.Here $\mathscr W$ consists of functions in 
$\mathscr O_u$ represented by absolutely convergent  Dirichlet series 
on the closure of $\mC_+$, $\mathscr A$ is the uniform closure of 
$\mathscr W$, $\mathscr H^\infty$ is the algebra of all bounded 
functions in $\mathscr O_u$, and $\mathscr O_b$ is set of all 
$f(s)=\sum_{n=1}^\infty a_n n^{-s}$ in $\mathscr O_u$ for which 
$f_r\in \mathscr H^\infty$, $r\in (0,1)$, where 
$f_r(s):=\sum_{n=1}^\infty a_n r^{\Omega(n)} n^{-s}$ and 
$\Omega(n)$ is the number of prime factors of $n$, 
counted with multiplicity. Let $S_\mN$ be the group of permutations 
of $\mN$. Each $\sigma\in S_\mN$ determines a completely 
multiplicative permutation $\hat\sigma\in  S_{\mN}$ 
(i.e., such that $\hat\sigma(mn)=\hat\sigma(n)\hat\sigma(m)$ 
for all $m,n\in \mN$) via the fundamental theorem of arithmetic. 
Then for a Dirichlet series $D=\sum_{n=1}^\infty a_n n^{-s}$, 
and $\sigma \in S_{\mN}$,  the formula  
$S_\sigma(D)=\sum_{n=1}^\infty a_{\hat{\sigma}^{-1}(n)} n^{-s}$, 
determines an action of the group $S_\mN$ on the set of all Dirichlet series. 
It is shown that each of the algebras above is invariant with 
respect to this group action. Given any subgroup $G$ of $S_\mN$, 
the set of $G$-invariant subalgebras of these algebras are also studied, 
and their maximal ideal spaces are described. These descriptions are 
used to characterise for these algebras the groups of units and their 
subgroups of invertible elements possessing logarithms, 
calculate the Bass stable rank,  show that some of these algebras 
are projective free, and describe when the special linear group can 
be generated by elementary matrices, with bounds on the number 
of factors when a factorisation exists.
\end{abstract}

\newpage
 
{\small \tableofcontents }
 
%==============
\section{Introduction}
\label{Section_1}
%==============

In recent years, there has been an increased interest in the study 
of Dirichlet series, i.e., series of the form $\sum_{n=1}^\infty a_nn^{-s}$, 
from the point of view of functional and harmonic analysis. 
Historically, the study of Dirichlet series  arose mainly in 
analytic number theory since the 19$^{\text{th}}$ century,  
and quickly became an active and fruitful area of analysis 
(see, e.g., \cite{HarRie}). The turning point in its development 
was the  work of Harald Bohr, who realised in 1913 that Dirichlet 
series and formal power series in infinitely many variables are closely related.
This discovery subsequently made it possible to use the tools of functional 
and multivariate complex analysis to answer questions about Dirichlet 
series via power series, see, e.g., the monograph \cite{DGMS} 
and references therein.

In this article, we focus on algebraic-topological properties 
of some algebras of holomorphic functions that arise in the 
study of series obtained by rearranging coefficients of Dirichlet series 
via permuting primes. 

Specifically, let $S_\mN$ be the group of permutations of $\mN$. 
Each $\sigma\in S_\mN$ determines a completely multiplicative 
permutation $\hat\sigma\in  S_{\mN}$ (that is, such that 
$\hat\sigma(mn)=\hat\sigma(n)\;\!\hat\sigma(m)$ for all $m,n\in \mN$) 
via  the fundamental theorem of arithmetic. Then for a Dirichlet series 
$D=\sum_{n=1}^\infty a_n n^{-s}$, and $\sigma \in S_{\mN}$, 
the formula $S_\sigma(D)=\sum_{n=1}^\infty a_{\hat{\sigma}^{-1}(n)} n^{-s}$, 
determines an action of the group $S_\mN$ on the set of all Dirichlet series.  

Let $\mathscr O_u$  be the algebra of holomorphic functions 
on the open right half-plane $\mC_+$ which are limits of Dirichlet 
series converging uniformly on proper half-planes. We study 
four natural topological subalgebras 
$\mathscr W, \mathscr A, \mathscr H^\infty, \mathscr O_b$ 
of  the algebra $\mathscr O_u$ which are invariant with respect to 
the above group action. Here the algebra $\mathscr W$ consists 
of functions in $\mathscr O_u$ represented by absolutely convergent  
Dirichlet series on the closure of $\mC_+$, the algebra $\mathscr A$ 
is the uniform closure of $\mathscr W$, $\mathscr H^\infty$ is the 
algebra of all bounded functions in $\mathscr O_u$, and 
$\mathscr O_b$ is the algebra of all $f(s)=\sum_{n=1}^\infty a_n n^{-s}$ 
in $\mathscr O_u$ for which $f_r\in \mathscr H^\infty$, $r\in (0,1)$,  
where $f_r(s):=\sum_{n=1}^\infty a_n r^{\Omega(n)} n^{-s}$ and 
$\Omega(n)$ is the number of prime factors of $n$, counted with multiplicity.

Analogous to algebras of symmetric polynomials in invariant theory, 
given a subgroup $G\subset S_\mN$, we introduce and study 
$G$-invariant subalgebras 
$\mA_G :=\!\{f\in\mA\, :\, S_\sigma (f)=f\ {\rm for\ all}\, \sigma\in G\}$ 
of algebras $\mA\in \{\mathscr O_b,\mathscr{W}, \mathscr{A}, \mathscr{H}^\infty\}$. 
We prove that the interesting $G$-invariant functions come only 
from non one-point finite orbits of the action of $G$ 
(see Theorems~\ref{te1.6} and \ref{te3.3}). We also show the existence 
of a continuous linear projection $\pi_G:\mA\to\mA_G$ such that 
$\pi_G(fg)=f\pi_G(g)$ for all $f\in\mA_G$ and $g\in\mA$ 
(see Theorem~\ref{teo1.7}).  

Then, in Theorems \ref{teo1.11}, \ref{teo1.12} and \ref{teo_contractibility}, 
we give a complete topological description of the maximal ideal spaces 
of algebras $\mA_G$ (as quotient spaces of finite- or infinite-dimensional 
open or closed polydiscs under the action of a profinite group). 
This is further used to characterise for these algebras the groups of units 
and their subgroups of invertible elements possessing logarithms, 
calculate the Bass stable rank (the concept introduced in \cite{Bas}, 
which plays an important role, e.g., in some stabilisation problems of $K$-theory), 
show that  some of these algebras are projective free (a notion that arose 
in connection with the Serre's problem from 1955 of determining 
whether finitely generated projective modules over  polynomial rings are free), 
and describe when the special linear group can be generated by 
elementary matrices, with new lower bounds on the number of factors 
when a factorisation exists (the topic goes back to 
the Gromov's Vaserstein problem, settled in \cite{IvaKut}). 
It is worth noting that the obtained results are new even for the case 
when $G$ is the trivial group (i.e., for the original topological algebras).

The paper is organised as follows. Sections~\ref{Section_formulation_results} 
and \ref{Sec2} contain formulations of our main results. 
In particular, in Section~\ref{Section_formulation_results}, 
basic definitions and results are formulated, including results on 
the topological structure of the maximal ideal spaces of algebras 
of invariant Dirichlet series. In turn, the aforementioned 
algebraic properties of these algebras are presented in Section~\ref{Sec2}. 
The subsequent sections mostly contain the proofs of these results. 
In Section~\ref{Section_3}, we collect and prove auxiliary results on 
holomorphic functions on open balls of complex Banach spaces $c_0$ 
and $\ell^\infty$, which will be needed to prove our main results. 
A detailed table of contents is provided above for the convenience of the reader. 

%==============
\section{Basic Definitions and Results}
\label{Section_formulation_results}
%==============

%---------------------------------------------
\subsection{Formal Dirichlet series} 
%---------------------------------------------

A {\em formal Dirichlet series} is a  series of the form 
\begin{equation}
\label{eq1.1}
\textstyle
D=\sum\limits_{n=1}^\infty a_n n^{-s},
\end{equation}
where $\{a_n\}_{n\in \mN}$ is a sequence of complex numbers, 
and $s$ is a complex variable. The set $\mathscr{D}$ of all 
formal Dirichlet series is a unital algebra with addition 
and multiplication given by the formulae: 
\begin{equation}
\label{eq1.2}
\begin{array}{l}
\sum\limits_{n=1}^\infty a_n n^{-s}+\sum\limits_{n=1}^\infty b_n n^{-s}
=\sum\limits_{n=1}^\infty (a_n+b_n)n^{-s},\medskip\\
\Big(\sum\limits_{n=1}^\infty a_n n^{-s}\Big) 
\cdot \Big(\sum\limits_{n=1}^\infty b_n n^{-s}\Big) 
= \sum\limits_{n=1}^\infty \Big(\sum\limits_{d\mid n} a_d b_{\frac{n}{d}}\Big)  n^{-s},
\end{array}
\end{equation}
and the unit $\mathbf{1}:=\sum_{n=1}^\infty \delta_{n1}n^{-s}$, 
where $\delta_{n1}=0$ for $n\neq 1$ and $ \delta_{11}=1$. 
It is a local algebra with the maximal ideal consisting of all 
series \eqref{eq1.1} with $a_1=0$. 

Recall that for the series  $D=\sum_{n=1}^\infty a_n n^{-s}\in \mathscr{D}$, 
the {\em abscissa of convergence} is defined by 
$$
\sigma_c(D)=
\inf\{\sigma \in \mR: D\textrm{ converges in}\ \{s\in\mC: {\rm Re}\, s>\sigma\}\}
\in [-\infty,\infty].
$$
If $\sigma_c(D)<\infty$, then the limit function $f$ of $D$ 
is holomorphic on the open half-plane $\{s\in\mC: {\rm Re}\,s>\sigma_c(D)\}$, 
see, e.g., \cite[Thm.~1.1]{DGMS}. Similarly, the 
abscissae of {\em absolute} and {\em uniform} convergence of $D$ 
are given by
\begin{eqnarray*}
\sigma_a(D)
\!\!\!&=&\!\!\!
\inf\{\sigma \in \mR: D\textrm{ converges absolutely in } 
\{s\in\mC: {\rm Re}\, s>\sigma\}\},\\
\sigma_u(D)
\!\!\!&=&\!\!\!
\inf\{\sigma \in \mR: D\textrm{ converges uniformly} 
\hspace{0.72mm}\textrm{ in } \{s\in\mC: {\rm Re}\, s>\sigma\}\}.
\end{eqnarray*}
Obviously, we have $-\infty\le\sigma_c(D)\le \sigma_u(D)\le \sigma_a(D)\le\infty$. 
It is also known that $\sigma_a(D)\le \sigma_c(D)+1$, 
see, e.g., \cite[Prop.~1.3]{DGMS}. 
 
If $\sigma_c(D)<\infty$, and  $f$ is the limit function of $D$, 
then the coefficients of $D$ can be  reconstructed by the formula
\begin{equation}
\label{coefficients}
a_n=\lim\limits_{R\rightarrow \infty} {\scaleobj{0.9}{\frac{1}{2Ri}}} 
{\scaleobj{0.9}{\int_{\kappa -iR}^{\kappa +iR} }}
f(s) n^s \text{d}s,
\end{equation}
for all $\kappa>\sigma_a(D)$ and $ n\in \mN$, 
where the path integral  is along the segment joining 
$\kappa-iR$ and $\kappa+iR$, see, e.g., \cite[Prop.~1.9]{DGMS}.   
Thus, if two Dirichlet series converge pointwise to the same 
function on a half-plane of $\mC$, then they are identical.
 
\phantom{$p_1$}Let the set of primes $\mP\subset\mN$ 
be labelled in ascending order as $p_1<p_2<p_3<\cdots$. 
By the fundamental theorem of arithmetic, every $n\!\in\! \mN$ 
may be written uniquely in the form 
$n\!=\!\prod_{p_i\in \mP} p_i^{\nu_{p_i}(n)}$, 
where $\nu_{p_i}(n)\!\in\! \mN\!\cup\!\{0\}$ denotes the largest integer 
$m$ such that $p_i^m$ divides $n$. 
 
The  {\em Bohr transform} $\mathfrak B$ associates with each series 
of $\mathscr D$ a formal power series over $\mC$ in 
countably many indeterminates $x_1,x_2,\dots$. 
It is defined as follows. If 
$$
\textstyle
D=\sum\limits_{n=1}^\infty a_n n^{-s}\in \mathscr D,
$$ 
then
\begin{equation}
\label{eq_pg_3}
\textstyle
\mathfrak{B} (D)(x_1,x_2,\dots):= 
\sum\limits_{n=1}^\infty a_n\Big(\prod\limits_{i\in \mN} x_i^{\nu_{p_i}(n)}\Big).
\end{equation}
If $\bbz=(z_1,z_2,\dots )\in\mC^\infty$ is such that the series 
$$
\textstyle
\sum\limits_{n=1}^\infty a_n\Big(\prod\limits_{i\in \mN} z_i^{\nu_{p_i}(n)}\Big)
$$
converges, then the limit of the series is called the 
{\em value of $\mathfrak B(D)$ at $\bbz$} 
(denoted by $\mathfrak B(D)(\bbz)$). 

By $\mC\llbracket x_1,x_2,\dots \rrbracket$ we denote 
the algebra of all formal power series as above. 
Then $\mathfrak B:\mathscr D\to \mC\llbracket x_1,x_2,\dots \rrbracket$ 
is an isomorphism of algebras, see, e.g., \cite[Sect.~14]{CasEve}.
 
%------------------------------- 
\subsection{The Fr\'echet algebras $\mathscr{O}_u$ and $\mathscr{O}_b$, 
and the Banach algebras $\mathscr{H}^\infty$, $\mathscr{W}$ and $\mathscr{A}$}
\label{subsect1.2}
%-------------------------------

By $\mathscr O_u$ we denote the set of holomorphic functions on 
$\mC_+:=\{s\in \mC: \textrm{Re }s>0\}$ that are limits of Dirichlet series 
that converge uniformly on proper half-planes of $\mC_+$. 
(Note that for any function in $\mathscr O_u$ there is only one 
Dirichlet series that converges to it, see \eqref{coefficients}.) 
Endowed with the metrisable locally convex topology 
defined by the system of seminorms
$$
P_\varepsilon(f):=\sup_{\textrm{Re}\,s>\varepsilon} |f(s)|,
\quad f\in \mathscr O_u,\quad \varepsilon>0,
$$
the space $\mathscr O_u$ becomes a Fr\'echet algebra 
with respect to pointwise operations, 
see, e.g., \cite[Thm.~2.6]{Bon} and references therein.
 
For $n\!\in\! \mN$, let  $\Omega(n)\!=\!\sum_{i\;\!\in\;\! \mN}\nu_{p_i}(n)$ 
be the number of prime factors of $n$, counted with multiplicity. 
For each $r\!\in\!\mC^*\!:=\!\mC\setminus\{0\}$, 
we consider a map $\cdot_{r}\!:\!\mathscr D\!\to\!\mathscr D$ 
sending $D\!=\!\sum_{n=1}^\infty a_n n^{-s}$ to 
$D_{ r}\!:=\!\sum_{n=1}^\infty  r^{\Omega(n)}a_n n^{-s}$. 
It is easy to see that the correspondence $r\mapsto \cdot_r$ 
determines a monomorphism of the multiplicative group $\mC^*$ 
in the group ${\rm Aut}(\mathscr D)$ of automorphisms 
of the algebra $\mathscr D$, cf. \eqref{eq1.2}.

Then it is easy to check that 
\begin{equation}
\label{eq_2_10}
\mathfrak B\circ \cdot_r=\cdot_{\texttt{r}} \circ \mathfrak B, 
\end{equation}
where $\cdot_{\texttt{r}}$ acts on formal power series in 
$\mC\llbracket x_1,x_2,\dots \rrbracket$ via the substitution 
$x_n\mapsto rx_n$, $n\in \mN$.  

\begin{proposition}\label{prop1.1}
If $D\in\mathscr D$ is such that $\sigma_u(D)\le 0,$ 
then $\sigma_u(D_r)\le 0$ for all $|r|\le 1$.
\end{proposition}

Thus, if a Dirichlet series $D$ converges uniformly on 
proper half-planes of $\mC_+$ to a function $g\in\mathscr O_u$, 
then for each $|r|\le 1$ the Dirichlet series $D_r$ converges uniformly 
on proper half-planes of $\mC_+$ to a function in $\mathscr O_u$ 
denoted by $g_r$.

We denote by $\mathscr O_b\subset\mathscr O_u$ the subset of functions 
$g$ such that the corresponding functions $g_r\in \mathscr O_u$, 
are bounded for all $r\in (0,1)$. Then since each $\cdot_r\in {\rm Aut}(\mathscr D)$, 
$\mathscr O_b$ is a Fr\'echet algebra with the metrisable locally convex topology 
defined by the system of seminorms
\begin{equation}
\label{equation_1.5}
P_r(g)=\sup_{s\;\!\in\;\!\mC_+}|g_r(s)|,\quad 0<r<1.
\end{equation}

\begin{proposition}\label{prop1.2}
For every nonconstant $g\in\mathscr O_b,$ the correspondence 
$t\mapsto \log P_{e^t}(g)$ determines an increasing convex function on $(-\infty, 0)$. 
Also$,$ $\lim\limits_{r\to 0}P_r(g)=\lim\limits_{\textrm{\em Re}\;\!s\to\infty}|g(s)|$ and 
$\lim\limits_{r\to 1}P_r(g)=\sup\limits_{s\;\!\in\;\!\mC_+}|g(s)|$.
\end{proposition}

\noindent In particular, $P_{r_1}(g)<P_{r_2}(g)<\sup_{s\;\!\in\;\!\mC_+}|g(s)|$ 
for all $0<r_1<r_2<1$.

Let $\mathscr{H}^\infty$ be the Banach algebra of bounded holomorphic functions 
from $\mathscr O_u$ equipped with the supremum norm. 
Due to Proposition \ref{prop1.2}, $\mathscr{H}^\infty$ is a subalgebra of $\mathscr O_b$. 
We denote by $\mathscr{W}$ the set of limit functions of Dirichlet series 
$$
\textstyle
D=\sum\limits_{n=1}^\infty a_nn^{-s} 
$$
with $\bba=(a_1,a_2,\dots)\in\ell^1$. Each such $D$ is absolutely convergent 
for all $s$ in the closure $ \overline{\mC}_+$ of $\mC_+$. 
Therefore, $\mathscr W\subset\mathscr H^\infty$. 
If $f\in\mathscr W$ is the limit function of $D$ as above, 
then we define 
\begin{equation}
\label{eq2.12}
\textstyle
\|f\|_1:=\|\bba\|_1:=\sum\limits_{n=1}^\infty |a_n|.
\end{equation}
With pointwise operations and the norm $\lVert\cdot\rVert_1,$  
$\mathscr{W}$ is a unital commutative complex Banach algebra, 
see, e.g., \cite{HW}. 

The closure of $\mathscr{W}$ in $\mathscr H^\infty$ 
is a uniform Banach algebra denoted by $\mathscr{A}$. 
It  can also be described as the subalgebra of  $\mathscr H^\infty$ 
of functions that are uniformly continuous on $\mC_+$, see, e.g., \cite{ABG}.
  
 %--------------------------
 \subsection{Action of the permutation group of $\mN$ on $\mathscr O_b$}
\label{subsection_3.1}
%----------------------------

Let $S_\mN$ be the group of permutations of $\mN$. 
Each $\sigma\in S_\mN$ determines a completely multiplicative permutation 
$\hat\sigma\in  S_{\mN}$ (i.e., such that $\hat\sigma(mn)=\hat\sigma(n)\hat\sigma(m)$ 
for all $m,n\in \mN$) via  the fundamental theorem of arithmetic 
by setting\footnote{Recall that $\mP=\{p_1,p_2,\dots \}$ where $p_1<p_2<\cdots$.}
\begin{equation}
\label{eq1.8b}
\textstyle
\hat\sigma(n):=\prod\limits_{i \in   \mN} \bigl(p_{\sigma(i)}\bigr)^{\nu_{p_i}(n)}.
\end{equation}

There is a faithful representation of $S_\mN$ on $\mathscr D$, described as follows.

\begin{proposition}
\label{prop3.1a}
The correspondence $\sigma\mapsto S_\sigma,$ 
$$
\textstyle
S_{\sigma}(D):=\sum\limits_{n=1}^\infty a_{\hat\sigma^{-1}(n)} n^{-s}
\;\text{ for }\; D=\sum\limits_{n=1}^\infty a_n n^{-s}, 
$$ 
determines a monomorphism $S$ of $S_\mN$ in the group ${\rm Aut}(\mathscr D)$.
\end{proposition}
 
Let $\mA$ be one of the algebras 
$\mathscr{O}_b$, $\mathscr{W}, \mathscr{A}, \mathscr{H}^\infty$, 
and $\mathscr D_\mA\subset\mathscr D$ be the subalgebra of Dirichlet series 
converging to functions in $\mA$.

\begin{proposition}
\label{prop1.4}
The algebra $\mathscr D_{\mA}$ is invariant with respect to the action $S,$ 
i.e.$,$ $S_\sigma(\mathscr D_{\mA})\subset \mathscr D_{\mA}$ 
for all $\sigma \in S_\mN$.
\end{proposition}  

In what follows, for a map $T: \mathscr D_\mA\to \mathscr D_\mA$, 
we retain the same symbol $T$ for the induced map of the algebra $\mA$ 
of limit functions of $\mathscr D_\mA$.

\begin{proposition}
\label{prop1.5}
The correspondence $\sigma\mapsto S_\sigma$ 
determines a monomorphism of $S_\mN$ 
in the group ${\rm Aut_I}(\mA)$ of isometric algebra isomorphisms of $\mA$.
\end{proposition}

In the case of $\mA=\mathscr O_b$, the latter means that 
for every seminorm $P_r$, $r\in (0,1)$, see \eqref{equation_1.5}, 
we have
$$
P_r(S_\sigma(g))=P_r(g)\ 
\textrm{  for all }g\in\mathscr O_b,\ \sigma\in S_\mN.
$$

Given a subgroup $G$ of $S_\mN$, a series $D\in \mathscr{D}$ 
is said to be {\em $G$-invariant} if $D=S_\sigma (D)$ 
for all $\sigma\in G$. Equivalently, $D=\sum_{n=1}^\infty a_n n^{-s}$ 
is $G$-invariant if and only if 
\begin{equation}
\label{eq3.18}
a_{\hat\sigma(n)}=a_n\ \textrm{ for all }\, n\in\mN,\ \sigma\in G.
\end{equation}
In turn, a function $f\in\mathscr O_b$ is called {\em $G$-invariant} 
if the Dirichlet series converging to $f$ is $G$-invariant 
(i.e., $S_\sigma(f)=f$ for all $\sigma\in G$).

For $\mA\in\{\mathscr O_b,\mathscr{W}, \mathscr{A}, \mathscr{H}^\infty\}$, 
the collection of all $G$-invariant elements in $\mA$ 
will be denoted by $\mA_G$. According to Propositions \ref{prop1.4}, \ref{prop1.5}, 
$\mA_G$ is a closed subalgebra of $\mA$. 

In what follows, for a subset  $Y$ of a set $X$, 
we denote the complement of $Y$ by $Y^{\com}:=X\setminus Y$.

Let $\mO\subset\mN$ and $\mO_\infty:= \mO^{\com}$ 
\label{pagenumber_mO_mOinfty}
be the unions of all finite and infinite orbits, respectively, 
of the action of $G$ on $\mN$. Let $\mN_{\mO}\subset\mN$ 
be the multiplicative unital semigroup generated by all $p_i\in\mP$ with $i\in\mO$. 

\begin{theorem}
\label{te1.6} 
The limit function of a Dirichlet series 
$\sum\limits_{n=1}^\infty a_n n^{-s}\in \mathscr D_{\mathscr O_b}$ 
belongs to $\mathscr O_{b,G}$ if and only if 
$$
\begin{array}{ccl}
\textrm{\em (i)} &\ a_n=0 &\textrm{for all } \ n\in (\mN_\mO)^{\com}, \textrm{ and }\\
\textrm{\em (ii)} & a_{\hat{\sigma}(n)}=a_n& \textrm{for all }\ n\in \mN_\mO,\ \sigma \in G.
\end{array}
$$
\end{theorem}

\begin{remark}
Suppose the action of $G\subset\mN$ on $\mN$ has no one-point orbits. 
Then Theorem~\ref{te1.6} implies that a necessary condition for 
$\mathscr O_{b,G}$ to contain nonconstant functions is that $G$ 
has a proper finite index subgroup. (E.g., this is not true 
if $G$ is isomorphic to the Higman group, see \cite{Hig}, 
i.e., the group with generators $g_1,g_2,g_3,g_4$ that satisfy the relations 
$g_1^{-1}g_2 g_1\!=\!g_2^2$, $g_2^{-1}g_3 g_2\!=\!g_3^2$, 
$g_3^{-1}g_4 g_3\!=\!g_4^2$, $g_4^{-1}g_1 g_4\!=\!g_1^2$.)
\end{remark}

Given $n\in  \mN_\mO$,  let $\textsf{O}(n)=\{\hat\sigma(n)\in\mN\, :\, \sigma\in G\}$ 
be the (finite) orbit of $n$ under the action of $G$, 
and let $|\textsf{O}(n)|$ be its cardinality. 
Consider a linear map $\pi_G:\mathscr D\to\mathscr D$ defined as follows:

If $D=\sum\limits_{n=1}^\infty a_n n^{-s}$, 
then $\pi_G(D):=\sum\limits_{n=1}^\infty b_n n^{-s}$, where
\begin{equation}
\label{eq1.5}
b_n:=\left\{\begin{array}{ccl}
0&{\rm if}&n\in (\mN_\mO)^{\com},\\[0.09cm] 
{\scaleobj{1.2}{\frac{1}{|{\textsf O}(n)|}}} \sum\limits_{k\;\!\in\;\! {\textsf O}(n)}a_k 
&{\rm if}&n\in\mN_\mO.
\end{array}\right.
\end{equation}

\begin{theorem}
\label{teo1.7}
Let $\mA\in\{\mathscr O_b,\mathscr{W}, \mathscr{A}, \mathscr{H}^\infty\}$. 
Then  $\pi_G(\mathscr D_\mA)\subset \mathscr D_\mA,$ 
and the induced map $\pi_G:\mA\to\mA$ 
is a continuous linear projection of norm one onto $\mA_G,$ 
satisfying the property 
$$
\pi_G(fg)=f\pi_G(g) \ \textrm{ for all } \,f\in\mA_G,\ g\in\mA.
$$
\end{theorem}

In the above theorem, for the algebra $\mathscr O_b$, 
the {\em norm one} condition on  the continuous linear projection 
$\pi_G:\mathscr O_{b}\to\mathscr O_{b,G}$  means that 
for all seminorms $P_r$, $r\in (0,1)$, and all $f\in \mathscr O_b$, 
we have $P_r(\pi_G(f))\le P_r(f)$, with equality for $f$ of  constant value one.

%---------------------------------------
\subsection{Action of a profinite completion of $G_\mO$ on $\ell^\infty(\mO)$}
\label{sec1.4a}
%-------------------------------------
 
Let $\ker_\infty$ denote the kernel of the action of $G$ on $\mO$,
$$
\ker_\infty=\{\sigma \in G: \sigma(n)=n\textrm{ for all } n\in \mO\}.
$$
It is a normal subgroup of $G$, and  the quotient group 
$G_\mO:=G/\ker_\infty$ is naturally identified with 
a subgroup of the group $S_\mO$ of permutations of elements of $\mO$.

We label  the set of (finite)  orbits of the action of $G_\mO$ on $\mO$ 
as $O_1, O_2, O_3,\dots$  (for instance, such that  if $i<j$, 
then $\min O_i <\min O_j$).
 
The restriction of the action of $G_\mO$ to $O_i$ determines 
a homomorphism $\rho_i: G_\mO\to S_{O_i}$ 
(where $S_{O_i}$ denotes the permutation group of elements of $O_i$). 
Let
$$
\textstyle Q_\mO:=\prod\limits_{i}S_{O_i}.
$$ 
Endowed with the product topology, $Q_\mO$ becomes a {\em profinite group} 
(i.e., a compact, Hausdorff and totally disconnected group). 
We define the homomorphism $\rho: G_\mO\to Q_\mO$ by the formula
\begin{equation}
\label{eq1.6}
\rho(g):=(\rho_1(g),\rho_2(g),\dots)\in Q_\mO,\quad g\in G_\mO.
\end{equation}
It is clear that $\rho$ is a monomorphism. 
We denote by ${\scaleobj{0.9}{\widehat G_\mO}}$ 
the closure of $\rho(G)$ in $Q_\mO$. 
Then ${\scaleobj{0.9}{\widehat G_\mO}}$ is a profinite group as well, 
called a {\em profinite completion} of $G_\mO$.
 
Next, there is a monomorphism $\psi:Q_\mO\to S_\mO$ given by the formula
$$
(\psi(\textsl{q}))(n):=q_i(n),\quad 
\textsl{q}=(q_1,q_2,\dots)\in Q_\mO,\quad n\in O_i,\quad i=1,2,\dots.
$$ 
In particular, $\psi\circ\rho={\rm id}_{G_{\mO}}$, 
and so $G_\mO$ is a subgroup of $\psi({\scaleobj{0.9}{\widehat G_\mO}})$ 
(dense in the topology induced from $Q_\mO$).
 
The action $\psi$ of $Q_\mO$ on $\mO$ determines an (adjoint) action 
$\psi^*$ of $Q_\mO$  on the Banach space $\ell^\infty(\mO)$ 
(of bounded complex functions on $\mO$ equipped with the supremum norm),  
$$
(\psi^*(\textsl{q}))(f):= f\circ \psi(\textsl{q}),
\quad \textsl{q}\in Q_\mO,\quad  f\in \ell^\infty(\mO).
$$ 
Obviously, every $\psi^*(\textsl{q})$ is a linear isometry of $\ell^\infty(\mO)$. 
 
Further, the closed unit ball of  $\ell^\infty(\mO)$ is the set 
$$
{\scaleobj{0.9}{\overline{\mD}^{\mO}}}
:=\{f\in \ell^\infty(\mO)\, :\, f(\mO)\subset \overline{\mD}\},
$$  
where $\overline{\mD}$ is the closure of the open unit disc $\mD\subset\mC$. 
We endow ${\scaleobj{0.9}{\overline{\mD}^{\mO}}}$ with the product topology, 
denoted by $\tau_p$, so that it becomes a compact Hausdorff space.
 
\begin{proposition}
\label{prop1.8}
The correspondence 
$\textsl{q}\mapsto \psi^*(\textsl{q})|_{{\scaleobj{0.9}{\overline{\mD}^{\mO}}}}$ 
determines  a continuous group action of $Q_\mO$ on 
$({\scaleobj{0.9}{\overline{\mD}^{\mO}}},\tau_p)$.
\end{proposition}
 
Let ${\scaleobj{0.9}{\overline{\mD}^{\mO}}}/{\scaleobj{0.9}{\widehat{G}_\mO}}$ 
be the quotient space by the action of $\psi^*|_{{\scaleobj{0.9}{\widehat{G}_\mO}}}$ 
on ${\scaleobj{0.9}{\overline{\mD}^{\mO}}}$ 
(i.e., the set of all orbits of ${\scaleobj{0.9}{\overline{\mD}^{\mO}}}$ 
under the action $\psi^*|_{{\scaleobj{0.9}{\widehat{G}_\mO}}}$ of 
${\scaleobj{0.9}{\widehat{G}_\mO}}$) endowed with the quotient topology. 
Hence, if $\pi: {{\scaleobj{0.9}{\overline{\mD}}^{\mO}}}
\to {\scaleobj{0.9}{\overline{\mD}^{\mO}}}/{\scaleobj{0.9}{\widehat{G}_\mO}}$ 
is the quotient map, then 
$U\subset {\scaleobj{0.9}{\overline{\mD}^{\mO}}}/{\scaleobj{0.9}{\widehat{G}_\mO}}$ 
is open if and only if $\pi^{-1}(U)$ is an open subset of 
$({\scaleobj{0.9}{\overline{\mD}^{\mO}}},\tau_p)$. 
Then, since ${\scaleobj{0.9}{\overline{\mD}^{\mO}}}$ and 
${\scaleobj{0.9}{\widehat{G}_\mO}}$ are compact Hausdorff spaces, 
${\scaleobj{0.9}{\overline{\mD}^{\mO}}}/{\scaleobj{0.9}{\widehat{G}_\mO}}$ 
is a compact Hausdorff space as well, see, e.g., \cite[Prop.~12.24]{Lee}. 

%-----------------------------------
\subsection{Some algebras of holomorphic functions on $\mD^{\mO}$}
\label{subsec_1.5}
%----------------------------------

In what follows, we also consider ${\scaleobj{0.9}{\mD^\mO}}$ 
equipped with the topology $\tau_{hk}$ of a {\em hemicompact $k$-space} 
with respect to an exhaustion of ${\scaleobj{0.9}{\mD^\mO}}$ 
by proper compact balls $r {\scaleobj{0.9}{\overline{\mD}^{\mO}}}$, $r\in (0,1)$, 
see,  e.g., \cite[Ch.~3,\,p.68]{Gol}  for the general definition. 
Thus, every compact subset $K\subset ({\scaleobj{0.9}{\mD^\mO}},\tau_{hk})$ 
is contained in one of $r {\scaleobj{0.9}{\overline{\mD}^{\mO}}}$, and 
a subset $F\subset ({\scaleobj{0.9}{\mD^\mO}},\tau_{hk})$ is closed 
if and only if each $F\cap  r{\scaleobj{0.9}{\overline{\mD}^{\mO}}}$, $r\in (0,1)$, 
is compact (in the product topology). It is clear that the inclusion map 
$({\scaleobj{0.9}{\mD^\mO}},\tau_{hk})\hookrightarrow 
({\scaleobj{0.9}{\overline{\mD}^\mO}},\tau_{p})$ is continuous, 
and an embedding on each ball $r{\scaleobj{0.9}{\overline{\mD}^{\mO}}}$, $r\in (0,1)$. 
Also, a function on $({\scaleobj{0.9}{\mD^\mO}},\tau_{hk})$ is continuous 
if and only if it is continuous in the product topology on each 
$r{\scaleobj{0.9}{\overline{\mD}^{\mO}}}$, $r\!\in \!(0,1)$, see, e.g., \cite[Prop.~2.7]{Ste}. 
Similarly to Proposition \ref{prop1.8}, we get
 
\begin{proposition}
\label{prop1.9}
The correspondence $\textsl{q}\mapsto \psi^*(\textsl{q})|_{{\scaleobj{0.9}{\mD^{\mO}}}}$ 
determines a continuous group action of $Q_\mO$ on 
$({\scaleobj{0.9}{\mD^{\mO}}},\tau_{hk})$.
\end{proposition}
 
Recall that  a function $f : {\scaleobj{0.9}{\mD^\mO}} \to \mC$ is said to be 
{\em G\^ateaux holomorphic} if for every $\bbz \in {\scaleobj{0.9}{\mD^\mO}}$ 
and $\bbh \in \ell^\infty(\mO)$, the function $ w\mapsto f(\bbz+w\bbh)$ 
is holomorphic in a neighbourhood of $0\in \mC$.  
The algebra $\mathcal O_w({\scaleobj{0.9}{\mD^\mO}})$ 
consists of all continuous G\^{a}teaux holomorphic functions 
on $({\scaleobj{0.9}{\mD^\mO}},\tau_{hk})$. Endowed with the system of seminorms
$$
P_r(f)=\|f(r\,\cdot)\|_\infty:=\sup\limits_{\bbz\;\!\in\;\! {\mD}^\mO} |f(r\bbz)|,
\quad 0<r<1,
$$ 
$\mathcal O_w({\scaleobj{0.9}{\mD^\mO}})$ becomes a Fr\'echet algebra. 
 
Next, $H_w^\infty({\scaleobj{0.9}{\mD^\mO}})$ and $A({\scaleobj{0.9}{\mD^\mO}})$ 
denote the subalgebras of $\mathcal O_w({\scaleobj{0.9}{\mD^\mO}})$ of bounded functions, 
and functions admitting continuous extensions to ${\scaleobj{0.9}{\overline{\mD}^{\mO}}}$, 
respectively, equipped with the supremum norm. 
 
For $i\in \mO$, define $z_i\in \ell^1(\mO)$ by $z_i(j)=\delta_{ij}$, 
$j\in\mO$ (the Kronecker delta). 
Under the identification of $\ell^1(\mO)$ with a subspace of the dual space 
$(\ell^\infty(\mO))^*$, the restriction $z_i|_{{\scaleobj{0.9}{\mD^\mO}}}$ 
belongs to $A({\scaleobj{0.9}{\mD^\mO}})$. 
The subalgebra $W({\scaleobj{0.9}{\mD^{\mO}}})$ of 
$A({\scaleobj{0.9}{\mD^{\mO}}})$ consists of pointwise limits on 
${\scaleobj{0.9}{\mD^{\mO}}}$ of complex power series 
in the variables $z_i$, $i\in \mO$, that are absolutely convergent 
on ${\scaleobj{0.9}{\overline{\mD}^{\mO}}}$. Define the norm 
$\lVert f\rVert_1$ of $f\in W({\scaleobj{0.9}{\mD^{\mO}}})$ 
as the sum of moduli of coefficients of the power series expansion of $f$. 
Then $(W({\scaleobj{0.9}{\mD^\mO}}),\lVert\cdot\rVert_1)$ is a 
complex unital Banach algebra. 

Let $\mathcal A\in\{\mathcal O_w({\scaleobj{0.9}{\mD^\mO}}), 
H_w^\infty({\scaleobj{0.9}{\mD^\mO}}), A({\scaleobj{0.9}{\mD^\mO}}), 
W({\scaleobj{0.9}{\mD^\mO}})\}$. It is easy to see that the action 
$\psi^*$ of $Q_\mO$ on ${\scaleobj{0.9}{\mD^\mO}}$ 
determines an adjoint continuous action, denoted by $\Psi$ of $Q_\mO$, 
on the algebra $\mathcal A$ given by the formula
\begin{equation}
\label{eq1.7}
\Psi_{\textsl{q}}(f):=f\circ\psi^*(\textsl{q}),
\quad \textsl{q}\in Q_\mO,\quad f\in\mathcal A.
\end{equation}
  
Let $c_{\mO}:\mC_+\to ({\scaleobj{0.9}{\mD^\mO}},\tau_{hk})$ 
be the continuous map defined by 
\begin{equation}
\label{eq1.8}
(c_\mO(s))(i):=p_i^{-s}, \quad s\in\mC_+,\quad  i \in \mO.
\end{equation}

Let $\mathscr D[\mO]\subset\mathscr D$ be the subalgebra of 
Dirichlet series $D=\sum_{n=1}^\infty a_n n^{-s}$ 
with $a_n=0$ for all $n\in (\mN_\mO)^{\com}$.  
For $\mA\in\{\mathscr O_b,\mathscr{W}, \mathscr{A}, \mathscr{H}^\infty\}$, 
we denote by $\mA[\mO]\subset\mA$ the closed subalgebra of limit functions 
of Dirichlet series in $\mathscr D[\mO]\cap\mathscr D_\mA$. 
We consider $\mA[\mO]$ with the induced Fr\'echet/Banach algebra structure.

Let $(\mathcal A,\mA)$ be one of the pairs 
$(\mathcal O_w({\scaleobj{0.9}{\mD^\mO}}),\mathscr O_b[\mO])$, 
$(H_w^\infty({\scaleobj{0.9}{\mD^\mO}}), \mathscr H^\infty [\mO])$, 
$(A({\scaleobj{0.9}{\mD^\mO}}), \mathscr A[\mO])$, 
$(W({\scaleobj{0.9}{\mD^\mO}}),\mathscr W[\mO])$.

\begin{theorem}
\label{te1.9}
The pullback by $c_\mO$ determines an isometric isomorphism of algebras 
$c_\mO^*:\mathcal A\to\mA$ such that
$$
S_{\psi(\textsl{q})}\circ c_\mO^*=c_\mO^*\circ\Psi_{\textsl{q}},
\quad \textsl{q}\in Q_\mO.
$$
\end{theorem}
  
See Section~\ref{subsection_3.1} for the definition of the action $S$.
 
For $\mathcal A\in
\{\mathcal O_w({\scaleobj{0.9}{\mD^\mO}}), H_w^\infty({\scaleobj{0.9}{\mD^\mO}}), 
A({\scaleobj{0.9}{\mD^\mO}}), W({\scaleobj{0.9}{\mD^\mO}})\}$,  
we denote by $\mathcal A_{{\scaleobj{0.9}{\widehat G_\mO}}}$ 
the subalgebra of functions  invariant with respect to the action 
$\Psi|_{{\scaleobj{0.9}{\widehat G_\mO}}}$, i.e., $f\in \mathcal A$ such that 
$$
\Psi_{\textsl{q}}(f)=f 
\ \textrm{ for all }\, \textsl{q}\in {\scaleobj{0.9}{\widehat G_\mO}}.
$$
Since the action $\Psi: {\scaleobj{0.9}{\widehat G_\mO}}\times\mathcal A\to\mathcal A$ 
is continuous, $\mathcal A_{{\scaleobj{0.9}{\widehat G_\mO}}}$ 
is a closed subalgebra of $\mathcal A$. Now, as $\rho(G)$ is a dense subgroup 
of ${\scaleobj{0.9}{\widehat G_\mO}}$,  Proposition~\ref{prop1.8} 
and Theorem~\ref{te1.9} imply:

\begin{corollary}
\label{cor1.10}
The homomorphism $c_\mO^*$ maps the algebra 
$\mathcal A_{{\scaleobj{0.9}{\widehat G_\mO}}}$ 
isomorphically onto the algebra $\mA_G$.
\end{corollary}
 
\begin{remark}
\label{rem1.11}
Let $\pi_G$ be the projection  of Theorem~\ref{teo1.7}. 
Theorem~\ref{te1.9} and Corollary~\ref{cor1.10} imply that
$$
P_{{\scaleobj{0.9}{\widehat G_\mO}}}:=(c_\mO^*)^{-1}\circ (\pi_G|_{\mA})\circ c_\mO^*,
$$
is a bounded linear projection of norm $1$ from $\mathcal A$ 
onto $\mathcal A_{{\scaleobj{0.9}{\widehat{G}_\mO}}}$, satisfying the property 
$$
P_{{\scaleobj{0.9}{\widehat G_\mO}}}(fg)
=fP_{{\scaleobj{0.9}{\widehat G_\mO}}}(g)
\text{ for all }f\in  \mathcal A_{{\scaleobj{0.9}{\widehat{G}_\mO}}} 
\text{ and }g\in\mathcal A.
$$
Let $m_{{\scaleobj{0.9}{\widehat{G}_\mO}}}$ be the normalised Haar measure 
on ${\scaleobj{0.9}{\widehat{G}_\mO}}$. 
Then $P_{{\scaleobj{0.9}{\widehat G_\mO}}}$ 
can also be described by the following:

\begin{proposition}
\label{prop1.11}
We have for all $f\in\mathcal A$,  
\begin{equation}
\label{proj}
(P_{{\scaleobj{0.9}{\widehat G_\mO}}}(f))(\bbz) 
=\int_{{\scaleobj{0.9}{\widehat G_\mO}}}(\Psi_{\textsl{q}}(f))(\bbz)\;\! 
\text{\em d}m_{{\scaleobj{0.9}{\widehat{G}_\mO}}}\!(\textsl{q}),\;\;\  
\bbz\in{\scaleobj{0.9}{\mD^\mO}}.
\end{equation}
\end{proposition}

For a compact subset $K$ of ${\scaleobj{0.9}{\mD^\mO}}$ 
which is $\psi^*|_{{\scaleobj{0.9}{\widehat G_\mO}}}$-invariant 
(that is, $(\psi^*(\textsl{q}))(K)\subset K$ for all 
$\textsl{q}\in {\scaleobj{0.9}{\widehat G_\mO}}$), 
let $C(K)_{{\scaleobj{0.9}{\widehat G_\mO}}}$ denote the subalgebra of $C(K)$ 
(the Banach algebra of complex-valued continuous functions on $K$) 
of functions invariant with respect to the action $\Psi|_{{\scaleobj{0.9}{\widehat G_\mO}}}$ 
(i.e., such that $\Psi_{\textsl{q}}(f)=f$ for all $\textsl{q}\in {\scaleobj{0.9}{\widehat G_\mO}}$). 
Then it is worth noting that due to Proposition~\ref{prop1.8}, 
formula \eqref{proj} gives  a linear continuous projection of norm $1$ 
also from $C(K)$ onto $C(K)_{{\scaleobj{0.9}{\widehat G_\mO}}}$.
\end{remark}
 
%-------------------------------------------
\subsection{Maximal ideal spaces of $\mathscr{W}_G$ and $\mathscr{A}_G$}
%-------------------------------------------
 
Recall that for a commutative unital complex Banach algebra $A$,  
the maximal ideal space $M(A)\subset A^\ast$ is the set of 
nonzero homomorphisms $A \!\rightarrow\! \mC$, 
endowed with the Gelfand topology, the weak$^*$ topology of  $A^\ast$. 
It is a compact Hausdorff space contained in the unit sphere of $A^\ast$. 
The {\em Gelfand transform}, defined by $\hat{a}(\varphi) := \varphi(a)$ 
for $a \in A$ and $\varphi \in M(A)$, is a morphism from $A$ into $C(M(A))$, 
the Banach algebra of complex-valued continuous functions on $M(A)$. 
Moreover, $\|\hat{a}\|_\infty:=\sup_{\varphi \;\!\in\;\! M(A)}|\varphi(a)|\leq \|a\|$ 
for all $a\in A$.

It is clear that there are subalgebras 
$W({\scaleobj{0.9}{\overline{\mD}^{\mO}}}/{\scaleobj{0.9}{\widehat G_\mO}})$ and 
$A({\scaleobj{0.9}{\overline{\mD}^{\mO}}}/{\scaleobj{0.9}{\widehat G_\mO}})$ 
of $C({\scaleobj{0.9}{\overline{\mD}^{\mO}}}/{\scaleobj{0.9}{\widehat G_\mO}})$,  
whose pullbacks by the quotient map $\pi:{\scaleobj{0.9}{\overline{\mD}^{\mO}}}
\to {\scaleobj{0.9}{\overline{\mD}^{\mO}}}/{\scaleobj{0.9}{\widehat G_\mO}}$ 
coincide with continuous extensions of 
$W({\scaleobj{0.9}{\mD^\mO}})_{{\scaleobj{0.9}{\widehat G_\mO}}}$ and 
$A({\scaleobj{0.9}{\mD^{\mO}}})_{{\scaleobj{0.9}{\widehat G_\mO}}}$ to
${\scaleobj{0.9}{\overline{\mD}^{\mO}}}$. 
These are Banach algebras with the norms induced from 
$W({\scaleobj{0.9}{\mD^\mO}})_{{\scaleobj{0.9}{\widehat G_\mO}}}$ and 
$A({\scaleobj{0.9}{\mD^{\mO}}})_{{\scaleobj{0.9}{\widehat G_\mO}}}$. 
Moreover, since the algebra $\mathscr A$ is the uniform closure of $\mathscr W$, 
Theorem~\ref{teo1.7} and Corollary~\ref{cor1.10} imply that 
$A({\scaleobj{0.9}{\overline{\mD}^{\mO}}}/{\scaleobj{0.9}{\widehat G_\mO}})$ 
is the uniform closure of 
$W({\scaleobj{0.9}{\overline{\mD}^{\mO}}}/{\scaleobj{0.9}{\widehat G_\mO}})$. 

\begin{theorem}
\label{teo1.11}$\;$
\begin{itemize}[leftmargin=0.72cm]
\item[{\em (1)}]
The maximal ideal spaces of $\mathscr W_G,$ $\mathscr A_G,$ 
$W({\scaleobj{0.9}{\overline{\mD}^{\mO}}}/{\scaleobj{0.9}{\widehat G_\mO}})$ 
and $A({\scaleobj{0.9}{\overline{\mD}^{\mO}}}/{\scaleobj{0.9}{\widehat G_\mO}})$ 
are homeomorphic. 

\smallskip 
 
\item[{\em (2)}]
$M(A({\scaleobj{0.9}{\overline{\mD}^{\mO}}}/{\scaleobj{0.9}{\widehat G_\mO}}))
={\scaleobj{0.9}{\overline{\mD}^{\mO}}}/{\scaleobj{0.9}{\widehat G_\mO}}$.
\end{itemize}
\noindent Here we identify points of 
${\scaleobj{0.9}{\overline{\mD}^{\mO}}}/{\scaleobj{0.9}{\widehat G_\mO}}$ 
with evaluation homomorphisms of the algebra 
$A({\scaleobj{0.9}{\overline{\mD}^{\mO}}}/{\scaleobj{0.9}{\widehat G_\mO}})$.
\end{theorem}

%------------------------------------------------------------ 
\subsection{Spectrum of  $\mathscr{O}_{b,G}$}
\label{subsec1.7}
%------------------------------------------------------------ 
 
Recall that for a commutative unital Fr\'echet algebra $A$, 
the {\em spectrum} $M(A)$ is the set of all 
nonzero continuous algebra homomorphisms $\varphi:A\rightarrow\mC$, 
endowed with the Gelfand topology, i.e., the weak$^*$ topology of $A^\ast$. 
For $a\in A$, we call the function $\hat{a}:M(A)\to \mC$, 
given by $\hat{a}(\varphi)=\varphi(a)$ for all $\varphi\in M(A)$, 
the {\em Gelfand transform of $a$}. 

For each $r\in (0,1)$, the compact ball $r{\scaleobj{0.9}{\overline\mD^{\mO}}}$ 
is invariant under the action 
$\psi|_{{\scaleobj{0.9}{\widehat G_{\mO}}}}$ of ${\scaleobj{0.9}{\widehat G_\mO}}$, 
and hence its image $\pi (r{\scaleobj{0.9}{\overline\mD^{\mO}}})$ is a compact subset of 
${\scaleobj{0.9}{\overline{\mD}^{\mO}}}/{\scaleobj{0.9}{\widehat G_\mO}}$ 
contained in ${\scaleobj{0.9}{\mD^\mO}}/{\scaleobj{0.9}{\widehat G_\mO}}$. 
We equip ${\scaleobj{0.9}{\mD^\mO}}/{\scaleobj{0.9}{\widehat G_\mO}}$ 
with the topology $\tau_{hk}$ of a hemicompact $k$-space 
with respect to the exhaustion of 
${\scaleobj{0.9}{\mD^\mO}}/{\scaleobj{0.9}{\widehat G_\mO}}$ 
by compact sets $\pi (r{\scaleobj{0.9}{\overline\mD^{\mO}}})$, $r\in (0,1)$. 
Thus, every compact subset $K\subset 
({\scaleobj{0.9}{\mD^\mO}}/{\scaleobj{0.9}{\widehat G_\mO}},\tau_{hk})$ 
is contained in one of $\pi (r{\scaleobj{0.9}{\overline\mD^{\mO}}})$, and a subset 
$F\subset ({\scaleobj{0.9}{\mD^\mO}}/{\scaleobj{0.9}{\widehat G_\mO}} ,\tau_{hk})$ 
is closed if and only if each $F\cap  \pi (r{\scaleobj{0.9}{\overline\mD^{\mO}}})$, $r\in (0,1)$, 
is compact in the topology of 
${\scaleobj{0.9}{\overline{\mD}^\mO}}/{\scaleobj{0.9}{\widehat G_\mO}}$.

Next, for each $r\in (0,1)$, the preimage of $\pi(r{\scaleobj{0.9}{\overline\mD^{\mO}}})$ 
under $\pi$ is $r{\scaleobj{0.9}{\overline\mD^{\mO}}}$. 
This shows that $\pi :({\scaleobj{0.9}{\mD^\mO}},\tau_{hk})\to 
({\scaleobj{0.9}{\mD^\mO}}/{\scaleobj{0.9}{\widehat G_\mO}},\tau_{hk})$ 
is a continuous proper map. Hence, there is a subalgebra 
$\mathcal O_w({\scaleobj{0.9}{\mD^{\mO}}}/{\scaleobj{0.9}{\widehat G_\mO}})$ 
of $C(({\scaleobj{0.9}{\mD^{\mO}}}/{\scaleobj{0.9}{\widehat G_\mO}},\tau_{hk}))$, 
such that its pullback  by $\pi$ coincides with the algebra 
$\mathcal O_w({\scaleobj{0.9}{\mD^\mO}})_{{\scaleobj{0.9}{\widehat G_\mO}}}$. 

We equip $\mathcal O_w({\scaleobj{0.9}{\mD^{\mO}}}/{\scaleobj{0.9}{\widehat G_\mO}})$ 
with the Fr\'echet topology transferred from 
$\mathcal O_w({\scaleobj{0.9}{\mD^\mO}})_{{\scaleobj{0.9}{\widehat G_\mO}}}$, 
so that $\pi^*: \calO_w({\scaleobj{0.9}{\mD^\mO}}/{\scaleobj{0.9}{\widehat{G}_\mO}}) 
\to \mathcal O_w({\scaleobj{0.9}{\mD^\mO}})_{{\scaleobj{0.9}{\widehat G_\mO}}}$ 
is an isometric algebra isomorphism.
  
\begin{theorem}
\label{teo1.12}$\;$
\begin{itemize}[leftmargin=0.6cm]
\item[{\em (1)}]
The spectra of the algebras $\mathscr O_{b,G}\!$ and 
$\mathcal O_w({\scaleobj{0.9}{\mD^{\mO}\!}}/{\scaleobj{0.9}{\widehat G_\mO}})$ 
are homeomorphic. 

\smallskip

\item[{\em (2)}]
$M(\mathcal O_w({\scaleobj{0.9}{\mD^{\mO}}}/{\scaleobj{0.9}{\widehat G_\mO}}))
=({\scaleobj{0.9}{\mD^{\mO}}}/{\scaleobj{0.9}{\widehat G_\mO}},\tau_{hk})$.
\end{itemize}
\end{theorem}

\subsection{Remarks}
\begin{enumerate}[leftmargin=*]
\item 
If the cardinality $|\mO|$ of $\mO$ is finite, then 
$({\scaleobj{0.9}{\mD^{\mO}\!}}/{\scaleobj{0.9}{\widehat G_\mO}},\tau_{hk})$ 
can be identified with a bounded Stein domain $\Omega\!\subset\!\mC^{|\mO|}$, 
and algebras 
$\mathcal O_w({\scaleobj{0.9}{\mD^{\mO}\!}}/{\scaleobj{0.9}{\widehat G_\mO}})$ and 
$A({\scaleobj{0.9}{\overline{\mD}^{\mO}}}/{\scaleobj{0.9}{\widehat G_\mO}})$ 
with algebras $\mathcal O(\Omega)$ and $A(\Omega)$, respectively, 
of holomorphic functions, and those that are continuous up to the boundary, 
see Section~\ref{subsekt2.3} below for the details. 

In general, $({\scaleobj{0.9}{\mD^\mO}}/{\scaleobj{0.9}{\widehat G_\mO}},\tau_{hk})$ 
is the inverse limit of a sequence of bounded Stein domains 
$\Omega_i\subset \mC^{|\mO_i|}$,  
where $\mO_i={\scaleobj{0.9}{\bigcup_{\;\!k=1}^{\;\!i}}} O_k$, $i\in \mN$, and 
the corresponding dual direct limits of  sequences of algebras 
$\mathcal O(\Omega_i)$ and $A(\Omega_i)$ are dense in 
$\mathcal O_w({\scaleobj{0.9}{\mD^\mO}}/{\scaleobj{0.9}{\widehat G_\mO}})$ and 
$A({\scaleobj{0.9}{\overline{\mD}^\mO}}/{\scaleobj{0.9}{\widehat G_\mO}})$ respectively, 
see Corollary~\ref{cor12.2} and Remark~\ref{rem13.3} below.

\item 
We also remark that the covering dimension of topological spaces 
${\scaleobj{0.9}{\overline{\mD}^{\mO}}}/{\scaleobj{0.9}{\widehat G_\mO}}$ and 
$({\scaleobj{0.9}{\mD^{\mO}}}/{\scaleobj{0.9}{\widehat G_\mO}},\tau_{hk})$ is $2|\mO|$. 
(Recall that the {\em covering dimension} of a topological space $X$ 
is the smallest integer $d$ such that every open covering of $X$ 
has an open refinement of order at most $d + 1$. 
If no such integer exists, then $X$ is said to have {\em infinite covering dimension}.) 

\item  
The structure of the maximal ideal space of the algebra $\mathscr H^\infty_G$ 
is unknown for all $|\mO|\ge 2$.  We only mention that 
$({\scaleobj{0.9}{\mD^{\mO}}}/{\scaleobj{0.9}{\widehat G_\mO}},\tau_{hk})$ 
is a subspace of $M(\mathscr H^\infty_G)$ and, in particular, 
the covering dimension of $M(\mathscr H^\infty_G)$ is $\ge 2|\mO|$. 
However, the answer to the question whether 
$({\scaleobj{0.9}{\mD^{\mO}}}/{\scaleobj{0.9}{\widehat G_\mO}},\tau_{hk})$ 
is dense in $M(\mathscr H^\infty_G)$ is unknown 
and is part of the corona problem for 
algebras of bounded holomorphic functions on polydiscs, 
one of the main open problems of  multivariate complex analysis. 
For $|\mO|=1$, the answer is affirmative due to the celebrated Carleson corona theorem 
for $H^\infty(\mD)$ (see \cite{Car}), 
and in this case the maximal ideal space is also well studied.  
\end{enumerate}

%===========================================
\section{Properties of Algebras of Invariant Dirichlet Series}
\label{Sec2} 
%===========================================

In this section,  we present some properties of algebras 
$\mathscr W_G$, $\mathscr A_G$, $\mathscr H_G^\infty$ and $\mathscr O_{b,G}$.
  
Recall that a topological space $X$ is {\em contractible} 
if the identity map $\textrm{id}_X:X \to X$ is null-homotopic, 
i.e., there exist an element $x_*\in X$ and a continuous map 
$\textrm{H}:[0,1]\times X \to X$ such that $\textrm{H}(0, \cdot)=\textrm{id}_X$ 
and $\textrm{H}(1,x)= x_*$ for all $x\in X$.
 
As an application of Theorems~\ref{teo1.11} and \ref{teo1.12}, we obtain
 
\begin{theorem} 
\label{teo_contractibility}
Let $\mA\in \{\mathscr W_G, \mathscr A_G, \mathscr O_{b,G}\}$. 
Then $M( \mA)$ is  contractible. 
\end{theorem}

%----------------------------------
\subsection{Group of units}
%----------------------------------

Theorem~\ref{teo_contractibility} can be applied, 
e.g., to describe groups of invertible elements of algebras of invariant Dirichlet series. 
In what follows, for a unital commutative complex Banach or Fr\'{e}chet algebra $A$, 
we denote by $A^{-1}$ the multiplicative group of all invertible elements of $A$, 
and by $e^A$ the subgroup of $A^{-1}$ of elements of the form $e^a$, $a\in A$. 
 
\begin{theorem}
\label{te2.2}
If $A\in \{\mathscr W_G, \mathscr A_G,\mathscr H_G^\infty\},$
then
$$
A^{-1}:=\big\{f\in A\, :\, \inf\limits_{{\scaleobj{0.9}{\mC_+}}} |f|>0\big\},
$$
and 
$$
(\mathscr O_{b,G})^{-1}:=\{f\in\mathscr O_{b,G}\, :\,  f_r\in (\mathscr H_G^\infty)^{-1},\, 
\forall r\in (0,1)\}.
$$ 
 
\noindent 
Moreover$,$ if $A\in \{\mathscr W_G, \mathscr A_G, \mathscr O_{b,G}\},$ 
then 
$$
A^{-1}=e^{A}.
$$
\end{theorem}
 
\noindent (Recall that if a function $f\in \mathscr O_b$ has 
Dirichlet series $D=\sum_{n=1}^\infty a_n n^{-s}$, 
then for $r\in (0,1)$, the function $f_r \in \mathscr H^\infty$, 
and has the Dirichlet series $D_{ r}:=\sum_{n=1}^\infty  r^{\Omega(n)}a_n n^{-s}$.)
 
\begin{remark}
\label{rem2.3} 
Theorem \ref{te2.2}  also gives the following: 
every $f=(\mathscr H_G^\infty)^{-1}$ has the form $f=e^{g}$ 
for some $g\in \mathscr O_{b,G}$ with a bounded real part. 
\end{remark}

%----------------------------------
\subsection{Stable rank} 
%----------------------------------

Let $A$ be an associative ring with unit. For $n\in \mN$, 
let $U_n(A)$ denote the set of {\em unimodular elements} of $A^n$, i.e., 
$$
U_n(A):= \{(a_1,\dots, a_n)\in A^n: Aa_1+\cdots+Aa_n=A\}.
$$
An element $(a_1,\dots, a_n)\in U_n(A)$ is {\em reducible} 
if there exists an element $(c_1,\dots, c_{n-1})\in A^{n-1}$ 
such that $(a_1+c_1a_n,\dots, a_{n-1}+c_{n-1}a_n)\in U_{n-1}(A)$. 
The {\em stable rank of $A$}, denoted by $\text{sr}\;\!A$, 
is the least $n\in \mN$ such that every element of $U_{n+1}(A)$ is reducible. 
The concept of the stable rank, introduced by Bass \cite{Bas}, 
plays an important role in some stabilisation problems of algebraic $K$-theory. 
For an $x\in \mR$,  $\lfloor x\rfloor$ denotes the greatest integer $\leq x$; 
any positive constant $c(x)$ depending on $x\in \mR$ 
will be extended as $c(\infty)=\infty$.

\begin{theorem}
\label{thm2.4} 
The following is true:
$$
\textstyle
\text{\em sr}\;\! \mathscr{W}_G
=\text{\em sr}\;\!\mathscr{A}_G
=\text{\em sr}\;\! \mathscr O_{b,G}
= \left\lfloor {\scaleobj{1.2}{\frac{|\mO|}{2}}} \right\rfloor +1
\le \text{\em sr}\;\!\mathscr H_G^\infty.
$$
\end{theorem}

%----------------------------------------
\subsection{Projective freeness} 
%----------------------------------------

For an associative ring $A$ with unit $1$, $\text{M}_n(A)$ denotes 
the $n \times n$ matrix ring over $A$, and 
$\text{GL}_n(A) \subset \text{M}_n(A)$ denotes the group of invertible matrices. 

A commutative unital ring $A$ is {\em projective free} if 
every finitely generated projective $A$-module is free. 
If $A$-modules $M,N$ are isomorphic, then we write $M\cong N$. 
If  $M$ is a finitely generated $A$-module, then  
(i) $M$ is {\em free} if $M \cong A^k$ for some integer $k\geq 0$, and 
(ii) $M$ is {\em projective} if there exists an $A$-module $N$ 
and an integer $n\geq 0$ such that $M\oplus N \cong A^n$. 
In terms of matrices (see, e.g., \cite[Prop.~2.6]{Coh}), 
the ring $A$ is projective free if and only if 
every idempotent matrix $P$ is conjugate (by an invertible matrix $S$) 
to a diagonal matrix with elements $1$ and $0$ on the diagonal, 
that is, for every $n\in \mN$ and every $P\in \text{M}_n(A)$ satisfying $P^2=P$, 
there exists an $S\in \text{GL}_n(A)$ such that  for some integer $k\geq 0$, 
$$
S^{-1} P S
=
[\begin{smallmatrix}
I_k & 0\\
0 & 0
\end{smallmatrix}].
$$ 
In 1976, it was shown independently by Quillen and Suslin that 
if $\mF$ is a field, then the polynomial ring $\mF[x_1, \dots , x_n]$ 
is projective free, settling Serre's conjecture from 1955 (see \cite{Lam}). 
In the context of a commutative  unital complex Banach algebra $A$, 
\cite[Thm.~4.1]{BruSas23} (see also \cite[Cor.~1.4]{BruSas}) 
says that the contractibility of the maximal ideal space $M(A)$ 
is sufficient for $A$ to be projective free. We have 

\begin{theorem}
\label{teo2.5a}
${\mathscr{W}}_G, {\mathscr{A}}_G, \mathscr O_{b,G}$ are projective free. 
\end{theorem}

\begin{remark}
The algebra $\mathscr H^\infty_G$ is projective free if $|\mO|=1$ 
(see, e.g., \cite[Cor.~3.30]{Qua} or \cite[Thm.~1.5]{BruSas}). 
However, nothing is known about the projective freeness of this algebra for $|\mO|\ge 2$. 
Recall that one of the obstructions for projective freeness is 
the nontriviality of the second \v{C}ech cohomology group 
$H^2(M(\mathscr H^\infty_G),\mZ)$, see, e.g., \cite[Thm.~4.1,\,Cor.~2.2]{BruSas23}.
\end{remark}

%----------------------------------
\subsection{Structure of the special linear group} 
\label{subsecttion2.4}
%----------------------------------

For an associative ring $A$ with unit $1$, 
\begin{itemize}[leftmargin=*]
\item[${\scaleobj{0.9}{\bullet}}$] 
$\text{SL}_n(A) \subset \text{GL}_n(A)$ denotes 
the subgroup of matrices with determinant $1$, and 
\item[${\scaleobj{0.9}{\bullet}}$] 
$\text{E}_n(A) \subset \text{SL}_n(A)$ denotes 
the subgroup of $\text{SL}_n(A)$ generated by all elementary matrices, 
i.e., matrices in $\text{SL}_n(A)$ which differ from the identity matrix $I_n$ 
by at most one non-diagonal entry,
\item[${\scaleobj{0.9}{\bullet}}$] 
$t_n(A)$ is the smallest integer $t$ such that every $\alpha\in \text{E}_n(A)$ 
is a product of at most $t$ {\em unitriangular matrices} in $\text{GL}_n(A )$ 
(that is, either upper triangular matrices with $1$ along the main diagonal 
or lower triangular matrices with $1$ along the main  diagonal)\footnote{By 
the Gaussian elimination method, each such matrix is in $\text{E}_n(A)$.}. 
We set $ t_n(A ) = \infty$ if such a $t$ does not exist.
\end{itemize}

It is well-known that if $A\!=\!\mF$, where $\mF$ is a field, 
then $\text{E}_n(\mF) \!=\! \text{SL}_n(\mF)$ and $t_n(\mF) = 4$ for all $n \ge 2$, 
and the same is true for all $A$ with $\text{sr}\;\!A=1$, see \cite[Lm.~9]{DenVas}. 
In general this is not always true. For instance, if $A$ is the polynomial ring 
$\mF[x_1, \dots , x_d]$, then $\text{E}_n(A) = \text{SL}_n(A)$ for $d = 1$ and all $n\in \mN$, 
as $A$ is a Euclidean ring, while for $d\ge 2$, 
$\text{E}_2(A) \subsetneq \text{SL}_2(A)$, see \cite[Prop.~(7.3)]{Coh0}, 
but $\text{E}_n(A) = \text{SL}_n(A)$ for all $n \ge 3$, see \cite[Cor.~6.7]{Sus}. 
However, $t_n(A) =\infty$ for all $n \ge 2$, even for $d = 1$, 
if $\mF$ is of infinite transcendence degree over its prime field, 
see \cite[Prop.~(1.5)]{Kal}. It is worth noting that according to \cite[Thm.\,20]{DenVas},
if $A$ has finite Bass stable rank and $t_m(A)<\infty$ for some $m\ge 2$, 
then $\lim\limits_{n\to\infty}t_n(A)\le 6$.

\smallskip

Let $A$ be a unital associative topological algebra over $\mC$. 
Then every group $\text{SL}_n(A)$, $n\ge 2$, has the structure of a topological group 
induced from $A$. Let $(\text{SL}_n(A))_0$ be the path component of the unit matrix $I_n$. 
It is a normal subgroup of $\text{SL}_n(A)$. 
Obviously, $\text{E}_n(A)$ is a subgroup of $(\text{SL}_n(A))_0$. 
In general, however, these groups do not coincide. 
They coincide, e.g., if $A$ is a complex unital Banach algebra or
$A$ is the algebra $C(X)$ of complex continuous functions 
on a normal topological space $X$ of finite covering dimension 
(equipped  with the topology of pointwise convergence on $X$), see \cite{Vas}. 
Moreover, in the latter case, let $v(d)$ be the supremum of $t_n(C(X))$ 
taken over the sets of all $n\in\mN$ and all normal topological spaces $X$ 
of covering dimension $\le d$. It was shown in \cite[Thm.~4]{Vas} that $v(d)\in\mN$. 
Thus  for all such $X$,
\begin{equation}
\label{e2.14}
4\le t_n(C(X))\, \le v(d)\; \text{ for all } \; n\ge 2.
\end{equation}
Further, due to \cite[Thm.~2.3]{Phi}, and using \cite[Lm.~2.1]{Bru22}, 
we obtain that for $X$ a compact manifold of dimension $d$ and $n\ge 2$
(see the argument in Section~\ref{subsecti12.5}), 
\begin{equation}
\label{eq2.15}
\textstyle
t_n(C(X))\ge \max\left\{\left\lfloor{\scaleobj{1.2}{\frac{d-2}{n^2-1}}} \right\rfloor, 4\right\}.
\end{equation}
In particular, $\lim\limits_{d\to\infty}v(d)=\infty$.

In turn, for the algebras considered in this paper,  the following holds.

\begin{theorem}
\label{theorem2.7}
If  $A\in \{\mathscr W_G, \mathscr A_G,\mathscr O_{b,G}\},$ 
then the topological groups $\text{\em SL}_n(A)$ are path connected$,$ 
i.e.$,$ $\text{\em SL}_n(A)=(\text{\em SL}_n(A))_0$. 

\noindent Moreover$,$ $\text{\em SL}_n(A)=\text{\em E}_n(A)$ only in the following cases
\begin{itemize}[leftmargin=0.72cm]
\item[$(a)$] 
$A\in \{\mathscr W_G, \mathscr A_G\}$ $($the case of Banach algebras$);$ 
\item[$(b)$]
$A=\mathscr O_{b,G}$ with $|\mO|<\infty$. 
\end{itemize}
\noindent In addition$,$ for all $n\ge 2,$ 
\begin{equation}
\label{equ2.16}
\textstyle
\max\left\{4, \left\lfloor{\scaleobj{1.2}{\frac{|\mO|}{n^2-1}}} \right\rfloor\!-\!3\right\}
\le t_n(A)\le v(2|\mO|)\!+\!5,
\end{equation}
where $v(\cdot)$ is the bound from \eqref{e2.14} and we assume that $v(\infty)=\infty$.
\end{theorem}

\begin{remark}$\;$
\begin{enumerate}[leftmargin=*]
\item For the algebra $\mathscr H_G^\infty$, 
the left-hand inequality of \eqref{equ2.16} is also true. 
The proof is the same as for the algebras $A$ above. 

\smallskip

\item It was proved in \cite{IvaKut} that 
if $A=\mathcal O(X)$ is the algebra of holomorphic functions on a reduced Stein space $X$, 
then $(\text{SL}_n(A))_0=\text{E}_n(A)$ for all $n\in \mN$. 
In turn, Theorem~\ref{theorem2.7} shows that this is not true anymore 
for infinite-dimensional analogs of Stein spaces, 
e.g., for $X={\scaleobj{0.9}{\mD^\mO}}$ with $|\mO|=\infty$ and $A=\mathcal O_w(X)$.

\smallskip

\item For a reduced Stein space $X$ of complex dimension $d$ and $A=\mathcal O(X)$,  
Corollary~\ref{cor12.5} and inequality \eqref{equ12.51} of Section~\ref{subsecti12.5} 
give the following hitherto unknown lower bounds, analogous to that of \eqref{equ2.16}:
$$
\textstyle
t_n(A)\ge \max\left\{4, \left\lfloor{\scaleobj{1.2}{\frac{d}{n^2-1}}} \right\rfloor\!-\!3\right\}
\quad {\rm for\ all}\quad n\ge 2.
$$
\end{enumerate}
\end{remark}
\smallskip

In a forthcoming paper, we will present some results for the algebra $\mathcal O_u$. 
We mention that $\mathcal O_u$ is $G$-invariant 
only for certain subgroups $G\subset S_{\mN}$. 
It is also worth noting that for the Fr\'echet algebra $\mathcal O_u$ itself, 
one can prove that its maximal ideal space is contractible, 
and thus $\mathcal O_u^{-1}=e^{\mathcal O_u}$, and
the algebra is projective free. 
Moreover, $\text{sr}\;\!\mathcal O_u =\infty$, 
$\text{E}_n(\mathcal O_u)\ne \text{SL}_n(\mathcal O_u)$, 
and $t_n(\mathcal O_u)=\infty$ for all $n\ge 2$.

%============================================
\section{Proofs of Propositions~\ref{prop1.1}-\ref{prop3.1a}} 
%============================================

%--------------------------------
\subsection{Preliminaries}
%--------------------------------

We use the following notation.
\begin{itemize}[leftmargin=*]
\item[${\scaleobj{0.9}{\bullet}}$] $\mN_0:=\mN\cup \{0\}$; 
 
\item[${\scaleobj{0.9}{\bullet}}$] for $N\in \mN$ and 
$\balpha =(\alpha_1,\dots, \alpha_N)\in \mN_0^N$, 
$\bbz:=(z_1,\dots, z_N)\in\mC^N$, $\bbz^{\balpha}:=z_1^{\alpha_1}\cdots z_N^{\alpha_N}$;
  
\item[${\scaleobj{0.9}{\bullet}}$] $\mT:=\{z\in \mC: |z|=1\}$; 
  
\item[${\scaleobj{0.9}{\bullet}}$] for $x\in \mR$, $x^+:=\max\{x,0\}$; 
  
\item[${\scaleobj{0.9}{\bullet}}$] $\pi: [2,\infty)\to\mN$ denotes the {\em prime counting function}, 
i.e., for $x\ge 2$, $\pi(x)$ is equal to the cardinality of the set $\{p\in \mP: p\leq x\}$. 
\end{itemize}

Next, we use the following upper bound of the abscissa of uniform convergence 
$\sigma_u(D)$ of the Dirichlet series $D=\sum_{n=1}^\infty a_n n^{-s}$: 
\begin{equation}
\label{eq_sigma_u}
\sigma_u(D)^+
= \bigg(\limsup\limits_{N\rightarrow \infty} 
\frac{\log (\sup_{t\;\!\in\;\! \mR}|\sum_{\;\!n=1}^{N} a_n n^{-it}|)}{\log N}\bigg)^+,
\end{equation}
see, e.g., \cite[Prop.~1.6]{DGMS}.  
 
Recall that the Bohr fundamental lemma (see, e.g., \cite[Thm.~3.2,\,p.78]{DGMS}) 
states that for all $N\in \mN$ and $a_1,\dots, a_N\in \mC$, 
$$
\textstyle 
\sup\limits_{t\;\!\in \;\!\mR}\Big|\sum\limits_{n=1}^N a_n n^{it}\Big|
=
\sup\limits_{\bbw\;\! \in \;\!\mT^{\pi(N)}}
\bigg|\sum\limits_{{\scaleobj{0.9}{\substack{\balpha \;\!\in\;\! 
\mN_0^{\pi(N)}\\[0.06cm]1\leq\bbp^{\balpha} \leq N}}}} 
a_{\bbp^{\balpha}} \bbw^{\balpha}\bigg|.
$$

%--------------------------------
\subsection{Proof of Proposition~\ref{prop1.1}}
%--------------------------------

Let  $D=\sum_{n=1}^\infty a_n n^{-s}$ be such that $\sigma_u(D)\!\leq\! 0$, 
and let $r\!\in\!\mC$, $|r|\!\leq\! 1$. Let $D_N\!=\!\sum_{n=1}^N a_n n^{-s}$, $N\!\in\!\mN$, 
be the $N^{\text{th}}$ partial sum of $D$. 
Then $\mathfrak B (D_N)\in \mC\llbracket x_1,x_2,\dots\rrbracket$  
is a polynomial in $\pi(N)$ variables $x_1,\dots, x_{\pi(N)}$. 
Evaluating it at points $\bbz=(z_1,\dots, z_{\pi(N)})\in\mC^{\pi(N)}$, 
we obtain a holomorphic polynomial on $\mC^{\pi(N)}$ 
(which we also denote by the same symbol $\mathfrak B (D_N)$). 

Next, by \eqref{eq_2_10}, 
$$
\mathfrak{B} ((D_N)_r)(\bbz)
= (\mathfrak B (D_N))_{\texttt{r}}(\bbz)
:=(\mathfrak B (D_N))(r\bbz),
\quad \bbz\in\mC^{\pi(N)}.
$$
Since the map $\bbz\mapsto r\bbz$ maps the closed unit polydisc 
${\scaleobj{0.9}{\overline{\mD}^{\,\pi(N)}}}$ to itself, 
the maximum modulus principle implies  that 
$$
\sup_{\bbz\;\! \in \;\!\mT^{\pi(N)}} | (\mathfrak{B} (D_N)_r)(\bbz)|
\!=\!
\sup_{\bbz \;\!\in \;\!\mT^{\pi(N)}} |(\mathfrak{B} (D_N))_{\texttt{r}}(\bbz)|
\!\leq\! 
\sup_{\bbz \;\!\in \;\!\mT^{\pi(N)}} |(\mathfrak{B} (D_N))(\bbz)|.
$$
Applying Bohr's fundamental lemma to both sides of the above inequality now yields
$$
\textstyle
\sup\limits_{t\;\! \in\;\! \mR} \Big|  \sum\limits_{n=1}^N r^{\Omega(n)} a_n n^{-s} \Big| 
\leq \sup\limits_{t\;\! \in\;\!  \mR} \Big|  \sum\limits_{n=1}^N a_n n^{-s} \Big| 
\quad {\rm for\ all}\quad N\in \mN.
$$
From here using formula \eqref{eq_sigma_u} for  $\sigma_u^+$, 
we obtain $\sigma_u(D_r)^+\!\le\!\sigma_u(D)^+\!=\! 0$, as claimed.
\hfill $\Box$

%--------------------------------
\subsection{Proof of Proposition~\ref{prop1.2}} 
\label{subsectio3.3}
%--------------------------------

Let $D=\sum_{n=1}^\infty a_n n^{-s}$ be the Dirichlet series 
that converges uniformly on proper half planes $\mH$ of $\mC_+$ 
to the nonconstant function $g \in {\mathscr{O}}_b$. 
As in the proof of Proposition~\ref{prop1.1}, using the Bohr fundamental lemma, 
we have 
\begin{equation}
\label{e2.13}
\sup_{s\;\!\in\;\! \mH}|(D_r)_n(s)-(D_r)_m(s)|
\le 
\sup_{s\;\!\in \;\!\mH}|D_n(s)-D_m(s)|
\end{equation}
for all $n,m\in\mN$ and all $r\in\mC$, $|r|\le 1$. 
Here  $D_k=\sum_{n=1}^k a_n n^{-s}$, $k\in \mN$. As $g\in  {\mathscr{O}}_u$, 
$\{D_n\}_{n\in \mN}$ is a Cauchy sequence in the Fr\'echet algebra $\mathscr{O}_u$. 
Each term of the sequence $\{(D_{\bcdot})_n\}_{n\in \mN}$ 
is a holomorphic function in  $(r,s)\in \mC\times \mC_+$, and by inequality \eqref{e2.13}, 
this sequence converges uniformly on domains $\overline{\mD}\times \mH$, 
where $\mH$ is a proper half plane of $\mC_+$, 
to the function $ g_{\bcdot}$ on $\overline{\mD}\times \mC_+$, 
such that the Dirichlet series of $g_r$ is 
$D_r=\sum_{n=1}^\infty r^{\Omega(n)}a_nn^{-s}$, $(r,s)\in \mD\times\mC_+$.  
Thus  $g_{\bcdot}$ is a  continuous function on $\overline{\mD}\times\mC_+$,  
holomorphic on $\mD\times \mC_+$, and 
bounded on each set $\{r\}\!\times \mC_+$  for $r\!\in\! [0,1)$, 
by the definition of class $\mathscr{O}_b$. 

Using the Bohr fundamental lemma, it can be seen that 
$$
\sup_{s\;\!\in\;\!\mH}|(D_r)_n(s)| = \sup_{s\;\!\in\;\!\mH}|(D_{|r|})_n(s)|
$$
for each proper half plane $\mH$ of $\mC_+$, 
and each $n\in\mN$, $r\in\overline{\mD}$. This fact, 
together with the uniform convergence of the sequence $\{(D_{\bcdot})_n\}_{n\in\mN}$ 
on domains $\mD\times\mH$, and the definition of the class $\mathscr O_b$, imply that 
$$
P_r(g) = \sup_{s\;\!\in\;\!\mC_+}|g_r(s)|,\quad r\in\mD,
$$ 
is a well-defined  logarithmically-subharmonic function in $r\in\mD$, which satisfies 
$$
P_{r}(g) = P_ {|r|}(g) \quad {\rm for\ all}\quad r\in\mD.
$$
Hence, the function $\log P_{e^t}(g)$ is convex and nondecreasing in $t\in (-\infty, 0)$. 
In fact, let us show that it is increasing. 

Indeed, for otherwise, as $\log P_{e^{{\bf\cdot}}}(g)$ is convex and nondecreasing, 
the function $P_r (g)$, $r\in [0,1)$, must be constant on some interval $[0, r_0)$, 
i.e., it is equal to $|a_1|=P_0(g)$ there, implying that 
the function $g_r(s)$ in $(r,s)\in \mD\times \mC_+$, 
when restricted to $(r_0\mD)\times \mC_+$, 
attains its maximum modulus (which is $|a_1|$) at each point $(0,s) \in \mD\times \mC_+$. 
Hence, the restriction of $g_r(s)$ to $ (r_0\mD)\times \mC_+$ is constant, 
and so the function $g_r(s)$ on $\mD\times\mC_+$ is constant as well, and equals $a_1$. 
This  implies that $g=a_1$, a contradiction to the assumption that $g$ is nonconstant. 

Finally, 
$$
\lim_{r\to 0} P_r(g) = P_0(g) = |a_1| = \lim_{\textrm{Re}\;\!s\to \infty}|g(s)|,
$$
where the last equality follows from 
the uniform convergence on proper half planes of $\mC_+$ 
of partial sums of $D$ to $g$. Similarly, using \eqref{e2.13} 
and the uniform convergence, for each 
half plane $\{s\in\mC_+:\textrm{Re}\;\!s>\varepsilon\}\subset\mC_+$, $\varepsilon>0$, 
we obtain: 
$$
\sup\limits_{\textrm{Re}\;\!s>\varepsilon}|g_r(s)|
\le \sup\limits_{\textrm{Re}\;\!s>\varepsilon}|g(s)|,
\quad r\in\mD, 
$$
and, moreover,
$$
\lim\limits_{r\to 1}\,\sup\limits_{\textrm{Re}\;\!s
> \varepsilon}|g_r(s)|
= \sup\limits_{\textrm{Re}\;\!s
> \varepsilon}|g(s)|.
$$
From the first equation, passing to the limit as $\varepsilon\to 0$, we get
$$
P_r(g):=\sup\limits_{s\;\!\in\;\!\mC_+}|g_r(s)| 
\le \sup\limits_{s\;\!\in\;\!\mC_+}|g(s)|, \quad r\in\mD,
$$
which shows that $\lim_{r\to 1}P_r(g)\le \sup_{s\;\!\in\;\!\mC_+}|g(s)|$. 
In turn, from the second equation, by the definition of $P_r(g)$, we obtain
$$
\lim\limits_{r\to 1}P_r(g)
\ge \lim\limits_{r\to 1}\sup\limits_{\textrm{Re}\;\!s
> \varepsilon}|g_r(s)|
= \sup\limits_{\textrm{Re}\;\!s
> \varepsilon}|g(s)|.
$$
Thus, passing to the limit as $\varepsilon\to 0$, we get
$$
\lim\limits_{r\to 1}P_r(g)\ge \sup\limits_{s\;\!\in\;\!\mC_+}|g(s)|.
$$
Thus, $\lim\limits_{r\to 1}P_r(g)= \sup\limits_{s\;\!\in\;\!\mC_+}|g(s)|$, as required.
\hfill $\Box$

%--------------------------------
\subsection{Proof of Proposition~\ref{prop3.1a}}
%--------------------------------

First, we prove that each $S_\sigma\in {\rm Aut}(\mathscr D)$.
 
It is clear that $S_\sigma:\mathscr D\to\mathscr D$ is a linear bijection. 
Further, for Dirichlet series $D(s)=\sum_{n=1}^\infty a_n n^{-s}$ 
and $\widetilde{D}(s)=\sum_{n=1}^\infty \widetilde{a}_n n^{-s}$, 
using the substitution $d\!=\!\hat\sigma(m)$ 
and complete multiplicativity of $\hat\sigma$, we get
$$
\!\!\!\begin{array}{rcl}
(S_{\sigma} (D\!\cdot \!\widetilde{D}))(s)
\!\!\!\!&=&\!\!\!\!
S_\sigma \Big(\sum\limits_{n=1}^\infty  
\Big(\sum\limits_{d|n} a_d \;\! \widetilde{a}_{\frac{n}{d}} \Big) n^{-s}\Big)
\!=\!
\sum\limits_{n=1}^\infty  
\Big(\sum\limits_{d\mid\hat\sigma^{-1}(n)} 
a_d\;\! \widetilde{a}_{\frac{\hat\sigma^{-1}(n)}{d}} \Big) n^{-s}
\\
\!\!\!\!&=&\!\!\!\!
\sum\limits_{n=1}^\infty \! 
\Big(\!\sum\limits_{\hat\sigma^{-1}(m)\mid\hat\sigma^{-1}(n)} 
\!\!a_{\hat\sigma^{-1}(m)} \widetilde{a}_{\frac{\hat\sigma^{-1}(n)}{\hat\sigma^{-1}(m)}} 
\Big) n^{-s}\\
\!\!\!\!&=&\!\!\!\!
\sum\limits_{n=1}^\infty \!\Big(\sum\limits_{m \mid n}  
a_{\hat\sigma^{-1}(m)} \widetilde{a}_{\hat\sigma^{-1}(\frac{n}{m})} 
\Big) n^{-s}
\\
\!\!\!\!&=&\!\!\!\!
\Big(\sum\limits_{n=1}^\infty a_{\hat\sigma^{-1} (n)} n^{-s} \Big) \!
\Big(\sum\limits_{n=1}^\infty \widetilde{a}_{\hat\sigma^{-1} (n)} n^{-s} \Big)
\!=\!
(S_\sigma (D)\!\cdot\! S_\sigma (\widetilde{D}))(s),
\end{array}
$$
as required.

Next, we prove that $S: S_\mN\to {\rm Aut}(\mathscr D)$ is a monomorphism. 

It is clear that the map $S$ is injective and that $S_{{\rm id}_{\mN}}={\rm id}_{\mathscr D}$. 
Thus it remains to check that $S_{\sigma_1\circ \sigma_2}=S_{\sigma_1}\circ S_{\sigma_2}$ 
for all $\sigma_1,\sigma_2\in S_\mN$. Indeed, for $D=\sum_{n=1}^\infty a_n n^{-s}$, we have 
$$
\begin{array}{rcl}
((S_{\sigma_1}\circ S_{\sigma_2})(D))(s)
\!\!\!\!&=&\!\!\!\!
(S_{\sigma_1} (S_{\sigma_2}(D)))(s)
=S_{\sigma_1}\Big( \sum\limits_{n=1}^\infty a_{\hat{\sigma}_2^{-1} (n)} n^{-s} \Big)
\\
\!\!\!\!&=&\!\!\!\!
\sum\limits_{n=1}^\infty a_{\hat{\sigma}_2^{-1} (\hat{\sigma}_1^{-1}(n))}  n^{-s}
=\sum\limits_{n=1}^\infty a_{(\hat{\sigma}_1\circ \;\!\hat{\sigma}_2)^{-1}(n)}  n^{-s}
\\
\!\!\!\!&=&\!\!\!\!
S_{\sigma_1\circ\;\! \sigma_2}\Big( \sum\limits_{n=1}^\infty a_{n} n^{-s} \Big)
=(S_{\sigma_1\circ \;\!\sigma_2}(D))(s).
\end{array}
$$
This completes the proof of the proposition.
\hfill$\Box$

%=====================================================
\section{Holomorphic Functions on Open Balls of $c_0$ and $\ell^\infty$} 
\label{Section_3}
%=====================================================

To proceed with the proofs, we require some additional information 
about the structure of some subalgebras of holomorphic functions 
on open balls of $c_0$ and $\ell^\infty$. This  is collected in the present section. 
In what follows, we will repeat some of the definitions from Subsection~\ref{subsec_1.5}, 
for the case when the union of finite orbits $\mO$ coincides with $\mN$. 

Let $\ell^\infty$ be the Banach space of bounded sequences of complex numbers 
equipped with the supremum norm, and let $c_0$ be the closed subspace 
of $\ell^\infty$ of sequences that converge to zero. 
The open unit ball of $\ell^\infty$ is the countable Cartesian product 
of open unit discs $\mD\subset\mC$, denoted by $\mD^\infty$, 
and called {\em the infinite-dimensional polydisc}. 
In turn, the open unit ball $B_{c_0}$ of $c_0$ is $c_0\cap \mD^\infty$. 
Let $c_{00}$ be the subspace of all finitely supported sequences in $c_0$, 
and $B_{c_{00}}:=B_{c_0}\cap c_{00}$. 

By $\mathcal O(B_{c_0} )$ we denote the algebra of {\em holomorphic}  
(i.e., complex Fr\'echet differentiable) functions on $B_{c_0}$. 
It is known that a function $f$ belongs to $\mathcal O(B_{c_0} )$ 
if and only if $f$ is continuous and $f|_{B_{c_{00}}\cap V}$ 
is holomorphic for every finite dimensional subspace $V\subset c_{00}$, 
see, e.g.,  \cite[Thm.~15.35, and Cor.~2.22]{DGMS}. 
 
The next result follows, e.g., from \cite[Thm.~2.19]{DGMS}.
 
\begin{theorem*}
There is an algebra monomorphism 
$\mathfrak O: \mathcal O(B_{c_0} )\to\mathscr{D}$ such that 
for each $f\in \mathcal O(B_{c_0} )$ 
and every $\bbz =(z_1,z_2,\dots)\in B_{c_{00}},$
\begin{equation}
\label{eq_pg_4}
f(\bbz)= \mathfrak B(\mathfrak O(f))(\bbz).\footnote{Recall that 
$\mathfrak B$ stands for the Bohr transform, see \eqref{eq_pg_3}.}
\end{equation}
Moreover$,$  if $\mathfrak O(f)(s)=\sum_{n=1}^\infty a_n n^{-s},$ 
then for every integer $N\ge\pi(n)$ and real numbers $0<r_1,\dots, r_N<1,$ 
\begin{equation}
\label{eq2.8}
a_n=\frac{1}{(2\pi i)^N}
{\scaleobj{0.9}{\int_{|\zeta_1|=r_1}\! \!\cdots \int_{|\zeta_N|=r_N} 
\frac{f(\zeta_1,\dots, \zeta_N, 0,\dots)}{\zeta_1^{\nu_{p_1}(n)+1}\cdots \zeta_N^{\nu_{p_N}(n)+1}} \text{\em d}\zeta_N\cdots \text{\em d}\zeta_1}},
\end{equation}
and  for all $0<r<1,$
\begin{equation}\label{eq2.9}
|a_n|\cdot r^{\Omega(n)}
\leq \sup\limits_{ \underset{1\le i\le\pi(n)}{\max} |z_i|\le r} |f(z_1,\dots, z_{\pi(n)},0,\dots)|.
\end{equation}
\end{theorem*}

Let $c : \mC_+\rightarrow B_{c_0}$ be the holomorphic map given by 
\begin{equation}
\label{eq2.10}
c(s):= (2^{-s},3^{-s}, \dots ,p_n^{-s},\dots),\quad s\in\mC_+.
\end{equation}
The image $c(L_\sigma)$ of each vertical line 
$L_{\sigma}:=\{\sigma+it\in\mC\,:\, \, t\in\mR\}$, $\sigma>0$, in $\mC_+$, 
is a relatively compact subset of $B_{c_0}$. 
Hence, for every $f\in\mathcal O(B_{c_0})$, 
the holomorphic function $f\circ c$ on $\mC_+$ is bounded 
on each line $L_\sigma$, $\sigma>0$. From this, 
and the Bohr theorems applied to proper half-planes of $\mC_+$, 
see, e.g., \cite[Thms.~1.13,~3.8]{DGMS},  one obtains

\begin{theorem**}
The pullback map $c^*,$ $f\mapsto f\circ c,$ 
determines an algebra monomorphism 
from $\mathcal O(B_{c_0})$ to $\mathscr O_u,$ 
such that for each $f\in \mathcal O(B_{c_0}),$ 
the Dirichlet series of $c^*(f)$ is $\mathfrak O(f)$.
\end{theorem**}

Note that $c^*(\mathcal O(B_{c_0}))$ is a proper (dense) subalgebra of $\mathscr O_u$. 
Indeed, from equation \eqref{eq_sigma_u}, 
it follows that the limit function of a Dirichlet series $D=\sum_{n=1}^\infty a_n n^{-s}$ 
belongs to $\mathscr O_u$ if and only if 
\begin{equation}
\label{eq2.6}
\limsup\limits_{N\rightarrow \infty} 
\frac{\log (\sup_{t\;\!\in\;\! \mR}|\sum_{n=1}^N a_n n^{-it}|)}{\log N}\leq 0.
\end{equation} 
In particular, this implies that the limit function $f$ 
of the Dirichlet series $D=\sum_{n=1}^\infty  p_{2^n}^{-s}$ 
belongs to $\mathscr O_u$ because $p_{2^n}\sim 2^n\log 2^n$ as $n\to\infty$ 
by the prime number theorem, see, e.g., \cite[pp.129-130]{Cha}. 
However, $f\not\in c^*(\mathcal O(B_{c_0}))$, as otherwise according to Theorem~B, 
$f= c^*(F)$, where $F(\bbz)=\lim_{N\to\infty}\sum_{n=1}^N z_{2^n}$, 
$\bbz=(z_1,z_2,\dots)\in B_{c_0}$, must belong to $\mathcal O(B_{c_0})$, 
a contradiction, as the limit does not exist, 
e.g., at points $\bbz\in B_{c_0}$ 
with the coordinates $z_i=\frac{t}{\log(i+2)}$, $i\in\mN$, $t\in (0,1]$.

\smallskip

Let $H^\infty(B_{c_0} )\subset \mathcal O(B_{c_0} )$ 
be the complex Banach algebra of bounded holomorphic functions on $B_{c_0}$ 
equipped with the supremum norm, and 
$\calA_u(B_{c_0})$ be the subalgebra of functions from $H^\infty(B_{c_0})$ 
that are uniformly continuous on $B_{c_0}$. 
Recall that $\mathscr{H}^\infty$ denotes the Banach algebra 
of bounded holomorphic functions from $\mathscr O_u$ 
equipped with the supremum norm, and 
$\mathscr A $ denotes the subalgebra of functions from $\mathscr{H}^\infty$ 
that are uniformly continuous on $\mC_+$. The following is true: 
\begin{itemize}
\item[(F1)] 
The homomorphism $c^*$ of Theorem~B 
maps $H^\infty(B_{c_0})$ isometrically onto $\mathscr{H}^\infty$ (see \cite{HedLinSei}), 
and $\calA_u(B_{c_0})$ isometrically onto $\mathscr{A}$ (see \cite[Thm.~2.5]{ABG}).
\end{itemize}

Let $W(B_{c_0})\subset \calA_u(B_{c_0})$ be the subspace of functions 
that are pointwise limits of  power series of $\mC\llbracket z_1,z_2,\dots \rrbracket$ 
(see \eqref{eq_pg_3}) with finite sums of moduli of their coefficients. 
Clearly, $c^*$ maps  $W(B_{c_0})$ onto $\mathscr{W}$ 
(the Banach algebra of limit functions of Dirichlet series 
that are absolutely convergent at each point of $\overline{\mC}_+$, 
equipped with the norm equal to the $\ell^1$ norm 
of the vector of coefficients of the Dirichlet series, see \eqref{eq2.12}). 
We introduce a norm on $W(B_{c_0})$ so that 
the isomorphism $c^*: W(B_{c_0})\to\mathscr{W}$ is an isometry 
(i.e., if $f\in W(B_{c_0})$ is the limit function of 
a series $S\in\mC\llbracket z_1,z_2,\dots \rrbracket$, 
then $\|f\|_1$ is defined to be the sum of moduli of coefficients of $S$). 
Thus, with pointwise operations and the norm $\lVert\cdot\rVert_1$, 
$W(B_{c_0})$ is a complex Banach algebra isometrically isomorphic via $c^*$ 
to the algebra $\mathscr W$. 
 
Next, the Fr\'echet algebra $\mathcal O_b(B_{c_0})$ 
consists of functions from  $\mathcal O(B_{c_0})$ 
which are bounded on each proper ball $rB_{c_0}$, $0<r<1$, 
equipped with the topology defined by the system of seminorms 
\begin{equation}
\label{eq3.20}
P_r(f):=\sup_{\bbz\;\!\in\;\! B_{c_0}} |f(r\bbz)|,\quad 0<r<1.
\end{equation}
In what follows, for a Banach space $X$, its norm is denoted by $\lVert\cdot\rVert_X$.

\begin{proposition}
\label{prop3.1} 
The image of  $\mathcal O_{b}(B_{c_0})$ under $c^*$ is $\mathscr O_b$. 
Moreover$,$ the map $c^*: \mathcal O_b(B_{c_0}) \to \mathscr O_b$ 
is an isometric isomorphism of the Fr\'echet algebras. 
\end{proposition}

This means that $c^*$ maps $\mathcal O_b(B_{c_0})$ onto $\mathscr O_b$ 
isometrically with respect to the chosen systems of seminorms 
on these Fr\'echet spaces, see \eqref{eq3.20} and \eqref{equation_1.5}.

\begin{proof}
Let $f\in \mathcal O(B_{c_0})$, and $f_{\texttt{r}}(\bbz):=f(r\bbz)$, $r\in (0,1]$. 
According to Theorem~B, each $c^*(f_{\texttt{r}})\in\mathscr O_u$. 
Moreover, if $\mathfrak O(f_{\texttt{r}})=\sum_{n=1}^\infty a_{n,r} n^{-s}$ 
is the Dirichlet series of $c^*(f_{\texttt{r}})$, then due to Theorem~A, 
for every integer $N\ge\pi(n)$ and real numbers $0<r_1,\dots, r_N<1,$ we have
$$
\!\begin{array}{rcl}
a_{n,r}\!\!\!\!&=&\!\!\!\!\displaystyle\frac{1}{(2\pi i)^N}
{\scaleobj{0.9}{\int_{|\zeta_1|=r_1}\! \!\cdots \int_{|\zeta_N|=r_N} 
\frac{f_{\texttt{r}}(\zeta_1,\dots, \zeta_N, 0,\dots)}{\zeta_1^{\nu_{p_1}(n)+1}\cdots 
\zeta_N^{\nu_{p_N}(n)+1}} 
\text{d}\zeta_N\cdots \text{d}\zeta_1}}
\\[0.48cm]
\!\!\!\!&=&\!\!\!\!\displaystyle\frac{1}{(2\pi i)^N}
{\scaleobj{0.9}{\int_{|\zeta_1|=r_1}\! \!\cdots \int_{|\zeta_N|=r_N} 
\frac{f(r\zeta_1,\dots, r\zeta_N, 0,\dots)}{\zeta_1^{\nu_{p_1}(n)+1}\cdots 
\zeta_N^{\nu_{p_N}(n)+1}} 
\text{d}\zeta_N\cdots \text{d}\zeta_1}}
\\[0.48cm]
\!\!\!\!&=&\!\!\!\!\displaystyle\frac{1}{(2\pi i)^N}
{\scaleobj{0.9}{\int_{|\zeta_1|=r\cdot r_1}\! \!\cdots \int_{|\zeta_N|=r\cdot r_N} 
\frac{r^{\Omega(n)}f(\zeta_1,\dots, \zeta_N, 0,\dots)}{\zeta_1^{\nu_{p_1}(n)+1}\cdots \zeta_N^{\nu_{p_N}(n)+1}} 
\text{d}\zeta_N\cdots \text{d}\zeta_1}}
\\[0.48cm]
\!\!\!\!&=:&\!\!\!\! r^{\Omega(n)} a_{n,1}.
\phantom{\displaystyle\frac{1}{(2\pi i)^N}{\scaleobj{0.9}{\int_{|\zeta_1|=r\cdot r_1}}}}
\end{array}
$$
This implies, see Section~\ref{subsect1.2} for the definition of ${\bm{\cdot}}_r$ , 
that 
\begin{equation}
\label{eq3.21}
\mathfrak O(f_{\texttt{r}})=(\mathfrak O(f))_r
\;\text{ for all }\; r\in (0,1),
\end{equation}
and, hence, 
\begin{equation}
\label{eq3.22}
c^*(f_{\texttt{r}})=(c^*(f))_r\;\text{ for all }\; r\in (0,1).
\end{equation}

Now, assume that $f\in \mathcal O_b(B_{c_0})$. 
Then each function $f_{\texttt{r}}\in H^\infty(B_{c_0})$, 
and so by (F1) and \eqref{eq3.22}, 
$(c^*(f))_r\in \mathscr H^\infty$ and, see \eqref{equation_1.5}, \eqref{eq3.20}, 
for every $r\in (0,1)$,
$$
P_r(c^*(f))
:= \|(c^*(f))_r\|_{\mathscr H^\infty}
= \|c^*(f_{\texttt{r}})\|_{\mathscr H^\infty}
= \|f_{\texttt{r}}\|_{H^\infty(B_{c_0})}
=: P_r(f)
< \infty.
$$
This implies that $c^*$ maps $\mathcal O_b(B_{c_0})$ to $\mathscr O_b$ 
isometrically with respect to the chosen systems of seminorms on these Fr\'echet spaces. 

To complete the proof, we need to check that 
$c^*$ maps $\mathcal O_b(B_{c_0})$ onto $\mathscr O_b$. 
 
To this end, let $D=\sum_{n=1}^\infty a_n n^{-s}$ be a Dirichlet series 
converging to $g\in \mathscr O_b$. Then, by the definition of $\mathscr O_b$, 
for each $r\in (0,1)$, the series $D_r:=\sum_{n=1}^\infty a_n r^{\Omega(n)} n^{-s}$ 
converges to a function $g_r\in \mathscr{H}^\infty$. So by (F1) and Theorem B, 
there exists a (unique) function $G^{(r)}\in H^\infty(B_{c_0})$, 
such that $g_r=c^*(G^{(r)})$ and $\mathfrak{O}(G^{(r)})=D_r$, $r\in (0,1)$. 
Using the equivariance from \eqref{eq_2_10}, we get 
$$
(\mathfrak{B} (D))_{{\texttt r}}
=\mathfrak{B}(D_r) 
=\mathfrak{B}(\mathfrak{O} (G^{(r)})) . 
$$ 
From this, and equation \eqref{eq_pg_4} of Theorem~A, we get 
$$
(\mathfrak{B} (D))(r\bbz )
=: (\mathfrak{B} (D))_{{\texttt r}}(\bbz)
= G^{(r)}(\bbz) \textrm{ for all } \bbz \in B_{c_{00}}.
$$
Next, if $r_1\!=\!tr_2$ for some $t\!\in\! (0,1)$, 
then the latter implies for all $\bbz\!\in\! B_{c_{00}}$,
$$
G^{(r_1)} (\bbz)
= (\mathfrak{B}(D))_{{\texttt r}_1} (\bbz)
= (\mathfrak{B}(D_{r_2}))_{{\texttt t}} (\bbz)
= (\mathfrak{B}(D_{r_2}))(t \bbz)
= G^{(r_2)} (t\bbz).
$$
From this, using the density of $B_{c_{00}}$ in $B_{c_0}$,  
and the continuity of  $G^{(r)}$ on $B_{c_0}$ for each $r\in (0,1)$, 
we obtain 
\begin{equation}
\label{e3.23}
G^{(r_1)} (\bbz)=G^{(r_2)} (t\bbz)
\quad {\rm for\ all}\ z\in B_{c_0}. 
\end{equation}
We set for $r\in (0,1)$ 
$$
\widetilde{G}^{(r)} (\bbz):=G^{(r)} (\mbox{$\frac{1}{r} \bbz$}), 
\quad \bbz\in r B_{c_0}.
$$
Then \eqref{e3.23} implies that if $r_1=tr_2$, for $t\in (0,1)$, 
then for all $\bbz\in r_1 B_{c_0}$, 
$$
\widetilde{G}^{(r_1)}(\bbz)
= G^{(tr_2)} (\mbox{$\frac{1}{tr_2} \bbz$}) 
= G^{(r_2)} (\mbox{$\frac{t}{tr_2} \bbz$}) 
= \widetilde{G}^{(r_2)} (\bbz).
$$
Therefore, the function $G$ given by 
$$
G(\bbz)=\widetilde{G}^{(r)} (\bbz) 
\textrm{ if } \bbz \in rB_{c_0},\ r\in (0,1),
$$
is well-defined on $B_{c_0}$, and belongs to $\mathcal{O}(B_{c_0})$ 
(because each $\widetilde{G}^{(r)}$ is holomorphic on $rB_{c_0}$, $r\in (0,1)$). 
Moreover, since all functions $\widetilde{G}^{(r)}$ 
are bounded on balls $rB_{c_0}$, $r\in (0,1)$,  
the function  $G$ belongs to $ \mathcal{O}_b (B_{c_0})$. 

We have, for all $s\in\mC_+$ and $r\in (0,1)$, that
\begin{equation}
\label{e3.24}
g_r(s) = c^{*}(G^{(r)})(s) := G^{(r)}(c(s)) = G(rc(s)).
\end{equation}
Clearly, $\lim\limits_{r\to 1}\, rc(s)=c(s)$ in $c_{0}$. Hence, by continuity of $G$, 
$$
\lim\limits_{r\to 1}G(rc(s)) = G(c(s)) = (c^*(G))(s).
$$
Similarly, as was shown in the proof of Proposition~\ref{prop1.2} (see Section~\ref{subsectio3.3}), the function $g_r(s)$, $(r,s)\in\overline{\mD}\times\mC_+$, is continuous. So 
$$
\lim\limits_{r\to 1}g_r(s)=g(s).
$$
Thus, passing to the limit as $r\to 1$ in \eqref{e3.24}, we obtain
$$
g(s)=(c^*(G))(s)\;\text{ for all }\; s\in\mC_+.
$$
Hence, $c^*(G)=g$.

This shows that the image of $\mathcal O_b(B_{c_0})$ under $c^*$ is $\mathscr O_b$, 
and completes the proof of the proposition. 
\end{proof}

A function $P\!:\!c_0\!\to\mC$ is an {\em $m$-homogeneous polynomial} 
if there exists a continuous $m$-linear form $A\!:\! c_0^m\!\to\mC$, 
such that $P(\bbz) \!=\! A(\bbz,\dots,\bbz)$ for every $\bbz\in c_0$. 
For convenience, we assume that the $0$-homogeneous polynomials are constant functions. 
It is known, see, e.g., \cite[Prop.~2.28]{DGMS}, that for every $f\in H^\infty(B_{c_0})$, 
there is a unique sequence $\{P_m\}_{m=0}^\infty$ 
of $m$-homogeneous polynomials on $c_0$, such that
\begin{equation}
\label{eq2.13}
\textstyle
f=\sum\limits_{m=0}^\infty P_m\; 
\textrm{ uniformly on } rB_{c_0}\textrm{ for every } 0<r<1,
\end{equation}
and
\begin{equation}
\label{eq2.14}
\sup_{B_{c_0}}|P_m| \le \sup_{B_{c_0}}|f|.
\end{equation}
Clearly, each $P_m\in \calA_u(B_{c_0})$ 
(i.e., it is in $H^\infty(B_{c_0})$ and uniformly continuous on $B_{c_0}$). 
Moreover, by \cite{Bog}, every homogeneous polynomial $P$ on $c_0$ 
is weakly uniformly continuous when restricted to $B_{c_0}$. 
Hence, $P$ has a unique extension $\widehat P$ to 
the closed unit ball $\overline{\mD}^\infty$ of $\ell^\infty$
which is continuous with respect to the weak$^*$ topology on $\ell^\infty=\ell_1^*$. Clearly, 
$$
\sup_{\overline{\mD}^\infty}|\widehat P| = \sup_{B_{c_0}}| P|.
$$
This and \eqref{eq2.13}, \eqref{eq2.14} imply that 
$f$ has a unique extension $\widehat f$ to the unit ball $\mD^\infty$ of $\ell^\infty$ 
such that $\widehat f|_{r\overline{\mD}^\infty}$, $0\!<\!r\!<\!1$, 
is weak$^*$ continuous, and 
\begin{equation}
\label{eq2.14a}
\textstyle
\widehat f=\sum\limits_{m=0}^\infty \widehat{P}_m
\;\textrm{ uniformly on }r\overline{\mD}^\infty\textrm{ for\ every } 0<r<1.
\end{equation}
Here we use the fact that the weak$^*$ closure of $rB_{c_0}$ is $r\overline{\mD}^\infty$. 
Also, it is worth noting that the weak$^*$ topology coincides with the product topology 
on each weak$^*$ compact subset of $\overline{\mD}^\infty$.

We denote by $H_{w}^\infty(\mD^\infty)$ 
the algebra of the extended functions $\widehat f$ with $f\in H^\infty(B_{c_0})$, 
equipped with the supremum norm. Then the restriction to $B_{c_0}$ 
determines an isometric isomorphism between 
$H_{w}^\infty(\mD^\infty)$ and $H^\infty(B_{c_0})$. 
Let $A(\mD^\infty)$ be the infinite {\em polydisc algebra}, 
i.e., the algebra of all weak$^*$ uniformly continuous, 
separately holomorphic functions on $\mD^\infty$, 
equipped with the supremum norm. By \cite{ABG}, we have 
\begin{itemize}
\item[(F2)] 
The subalgebra of $H_{w}^\infty(\mD^\infty)$ 
of the extended functions $\widehat f$ with $f\in \calA_u(B_{c_0})$, 
coincides with $A(\mD^\infty)$.
\end{itemize}

Let $W(\mD^\infty)$ be the subspace of $A(\mD^\infty)$ 
of functions which can be represented as 
pointwise limits of power series from $\mC\llbracket z_1,z_2,\dots \rrbracket$ 
that converge absolutely at each point of $\overline{\mD}^\infty$. 
If $f\in W(\mD^\infty)$ is the limit function of a series $S$, 
then $\|f\|_1$ is defined to be the sum of moduli of the coefficients of $S$. 
With pointwise operations and the norm $\lVert\cdot\rVert_1$, 
$W(\mD^\infty)$ is a complex Banach algebra. 
Clearly, we have
\begin{itemize}
\item[(F3)] 
The restriction to $B_{c_0}$ determines 
a Banach algebra isometric isomorphism 
between  $W(\mD^\infty)$ and  $W(B_{c_0})$. 
\end{itemize}
Let $f\in  \mathcal O_b(B_{c_0})$. Applying the Taylor series expansion \eqref{eq2.13} 
to the functions $f_\texttt{r}:=f(r\,\cdot)\in H^\infty(B_{c_0})$, $0<r<1$, 
and using its uniqueness, we obtain a similar expansion for $f$, 
i.e.,  there is a unique sequence $\{P_m\}_{m=0}^\infty$ 
of $m$-homogeneous polynomials on $c_0$, such that
\begin{equation}
\label{eq2.16}
\textstyle
f=\sum\limits_{m=0}^\infty P_m
\; \textrm{ uniformly on } rB_{c_0}\textrm{ for every } 0<r<1.
\end{equation}
Thus, as in the case of $H^\infty(B_{c_0})$, 
every $f\in\mathcal O_b(B_{c_0})$ has a unique extension $\widehat f$ to $\mD^\infty$, 
such that $\widehat f|_{r\overline{\mD}^\infty}$, $0\!<\!r\!<\!1$, 
is weak$^{*}$ continuous, and 
\begin{equation}
\label{eq2.17}
\textstyle
\widehat f=\sum\limits_{m=0}^\infty \widehat{P}_m
\;\textrm{ uniformly on }r\overline{\mD}^\infty\textrm{ for\ every } 0<r<1.
\end{equation} 
We use symbol $\mathcal O_{w}(\mD^\infty)$ to denote 
the algebra of the extended functions $\widehat f$ with $f\in \mathcal O_b(B_{c_0})$, 
equipped with the metrisable locally convex topology 
defined by the system of seminorms
$$
P_r(\widehat f)
:=\sup\limits_{\bbz\;\!\in\;\!\mD^\infty}|\widehat f(r\bbz)|\, (=P_r(f)),
\quad 0<r<1.
$$ 
Then $\mathcal O_{w}(\mD^\infty)$ is a Fr\'echet algebra, and moreover:
\begin{itemize}
\item[(F4)] 
The restriction to $B_{c_0}$ determines a Fr\'echet algebra isomorphism 
between $\mathcal O_{w}(\mD^\infty)$ and $\mathcal O_b(B_{c_0})$. 
\end{itemize}
  
Next, the homomorphism $c^*$ given by \eqref{eq2.10} 
can be extended to $\mathcal O_{w}(\mD^\infty)$ by the formula 
$$
c^*(\widehat f):=c^*(f).  
$$
Combining all the above facts, we get
\begin{itemize}
\item[(F5)] 
The extended homomorphism $c^*$ maps the algebras 
$\mathcal O_{w}(\mD^\infty)$, $H^\infty_w(\mD^\infty)$, $A(\mD^\infty)$, and $W(\mD^\infty)$ isometrically onto the algebras 
$\mathscr O_b$, $\mathscr H^\infty$, $\mathscr A$, and $\mathscr{W}$, respectively.
\end{itemize}

%=========================================================
\section{Action of the Permutation Group of $\mN$ on $\mathcal O(B_{c_0})$}
\label{SubSect_1.3}
%=========================================================

In this section we prove Propositions \ref{prop1.4} and \ref{prop1.5}.

The group $S_\mN$ of permutations of $\mN$ 
acts on the set $\mP$ of primes  by the formula 
$$
\sigma(p_n):=p_{\sigma(n)}, 
\quad \sigma\in S_\mN.
$$
Also, there exists a faithful representation $A$ of $S_\mN$ 
on  $\ell^\infty$ by linear isometries $A_\sigma$, $\sigma\in S_\mN$:
$$
A_\sigma(\bbz) := (z_{\sigma(1)},z_{\sigma(2)},\dots)
\;\text{ for }\; \bbz=(z_1,z_2,\dots)\in\ell^\infty.
$$
In particular, $c_0$ is a subrepresentation of $A$, 
i.e., $A_\sigma(c_0)\subset c_0$ for all $\sigma\in S_\mN$. 
Note that $A_\sigma$ is the second adjoint of $A_\sigma|_{c_0}$; 
hence, each operator $A_\sigma$ is weak$^*$ continuous. 

Since each $A_\sigma$ maps $B_{c_0}$ biholomorphically onto itself, 
the adjoint (faithful) representation $A^*$ on $\mathcal O(B_{c_0})$,  
$A_\sigma^*(f):=f\circ A_\sigma$, is well-defined. 
By the definition, each $A_\sigma^*\in {\rm Aut}(\mathcal O(B_{c_0}))$, 
the group of automorphisms of the algebra $\mathcal O(B_{c_0})$.

Further, recall that the action of $\sigma$ on $\mP$ 
determines a completely multiplicative permutation $\hat\sigma\in S_{\mN}$ 
(i.e. such that $\hat\sigma(mn)=\hat\sigma(n)\hat\sigma(m)$ for all $m,n\in \mN$) 
via the fundamental theorem of arithmetic, by setting 
$$
\textstyle
\hat\sigma(n):=\prod\limits_{p \in   \mP} \bigl(\sigma(p)\bigr)^{\nu_{p}(n)}.
$$ 
Our next result connects the actions of $S_\mN$ introduced above 
with the homomorphism $\mathfrak O$ of Theorem~A (see Section~\ref{Section_3}). 
Recall that the action of the permutation group on $\mathscr D$ is given by 
$$
\textstyle
S_\sigma (D) := \sum\limits_{n=1}^\infty a_{\hat{\sigma}^{-1} (n)} n^{-s} 
\;\textrm{ for }\; D=\sum\limits_{n=1}^\infty a_n n^{-s}\in \mathscr{D},
$$
see  Proposition~\ref{prop3.1a}.
  
\begin{theorem}
\label{te3.1}
The homomorphism $\mathfrak O:\mathcal O(B_{c_0}) \to \mathscr D$ 
is equivariant with respect to the actions $A^*$ and $S,$ that is$,$
$$
\mathfrak O\circ A_\sigma^* = S_\sigma\circ\mathfrak O
\,\textrm{ for all }\, \sigma\in S_\mN.
$$
\end{theorem}
\begin{proof}
We use equation \eqref{eq_pg_4} of Theorem~A. 
For all $f\in\mathcal O(B_{c_0})$ and $\bbz\in B_{c_{00}}$, 
since each $A_\sigma$ maps $B_{c_{00}}$ to itself, we have
\begin{equation}
\label{te3.1eq1}
\mathfrak B(\mathfrak O(A_\sigma^*(f)))(\bbz)
= A_\sigma^*(f)(\bbz)
= f(A_\sigma(\bbz)) 
= \mathfrak B(\mathfrak O(f))(A_\sigma(\bbz)).
\end{equation} 
Suppose that $\mathfrak O(f)=\sum_{n=1}^\infty a_n n^{-s}$. 
Then, by the definition of the Bohr transform, see \eqref{eq_pg_3}, we have
$$
\textstyle 
\mathfrak B(\mathfrak O(f))(\bbw) 
= \sum\limits_{n\;\! \in\;\!  \mN} a_n\prod\limits_{i\;\! \in \;\! \mN} w_{i}^{\nu_{p_i}(n)}, 
\quad \bbw\in B_{c_{00}}.
$$
Due to inequality \eqref{eq2.9} of Theorem~A, 
the above series converges absolutely at each point of $B_{c_{00}}$. 
Thus its sum $\mathfrak B(\mathfrak O(f))(\bbw)$ 
does not depend on the rearrangement of the terms of the series. 
Using this, and \eqref{te3.1eq1}, in particular 
we have at $\bbw=A_\sigma(\bbz)$, $\bbz\in B_{c_{00}}$, that 
$$
\begin{array}{rcl}
A_\sigma^*(f)(\bbz)
\!\!\!\!&=&\!\!\!\!
\mathfrak B(\mathfrak O(A_\sigma^*(f)))(\bbz)
=
\mathfrak B(\mathfrak O(f))(A_\sigma(\bbz))
=
\sum\limits_{n\;\! \in\;\!  \mN} a_n
\prod\limits_{i\;\! \in \;\! \mN} z_{\sigma(i)}^{\nu_{p_i}(n)}
\\[0.3cm]
\!\!\!\!&=&\!\!\!\!
\sum\limits_{n\;\! \in\;\!  \mN} a_{\hat\sigma^{-1}(n)}
\prod\limits_{i\;\! \in\;\!  \mN} z_{i}^{\nu_{\sigma^{-1}(p_i)}(\hat\sigma^{-1}(n))} 
=
\sum\limits_{n\;\! \in\;\!  \mN} a_{\hat\sigma^{-1}(n)} 
\prod\limits_{i\;\! \in \;\! \mN} z_{i}^{\nu_{p_i}(n)}. 
\end{array}
$$
Here we use that $\nu_{\sigma^{-1}(p_i)}(\hat\sigma^{-1}(n))=\nu_{p_i}(n)$. 
 
Substituting this expression for $A_\sigma^*(f)$ in formula \eqref{eq2.8} 
to compute the coefficients of the Dirichlet series $\mathfrak O(A_\sigma^*(f))$, 
and using the fact that the series on the right in the above equation 
converges uniformly on each polydisc 
$r\mD^N := rB_{c_{00}}\cap  \{\bbz\in c_0 : z_{i} 
= 0\text{ for all }i> N\}$, $r\in (0,1)$, $N\in\mN$, we obtain that  
$$
\textstyle 
\mathfrak O(A_\sigma^*(f)) 
= \sum\limits_{n=1}^\infty a_{\hat\sigma^{-1}(n)}n^{-s}
= S_\sigma\circ(\mathfrak O(f)).
$$
This completes the proof of the theorem.
\end{proof} 
  
\begin{proof}[Proof of Proposition~\ref{prop1.4}] 
For  $\mA=\mathscr H^\infty$, $\mathscr A$, $\mathscr W$, $\mathscr O_b$, 
let $\mathbf{A}$ denote 
$H^\infty(B_{c_0})$, $\mathcal A_u(B_{c_0})$, $W (B_{c_0})$, $\mathcal O_b (B_{c_0})$, respectively. Let $D\in \mathscr D_\mA$ 
(the subalgebra of Dirichlet series converging to functions in $\mA$). 
Then by (F1), Proposition~\ref{prop3.1}, and Theorem~B, 
there exists an $f\in \mathbf{A}$ such that $D=\mathfrak O (f)$. 
Let $\sigma \in S_\mN$. Using Theorem~\ref{te3.1}, we obtain 
\begin{equation}
\label{eq4.31}
S_\sigma (D) = S_\sigma (\mathfrak O (f)) = \mathfrak O(A_\sigma^* (f)).
\end{equation}
From the definitions of the classes 
$H^\infty(B_{c_0})$, $\mathcal A_u(B_{c_0})$, $W (B_{c_0})$, $\mathcal O_b (B_{c_0})$, 
and  using the facts that  $A_\sigma: B_{c_0}\to B_{c_0}$ is a biholomorphic map 
and $A_\sigma (rB_{c_0})= rB_{c_0}$, for all $r\in [0,1]$, 
it follows that $A_\sigma^*( f)=f\circ A_\sigma \in \mathbf{A}$. 
By \eqref{eq4.31}, and Theorem~B, we get that 
$S_\sigma(D)$ is the Dirichlet series of $c^*(A_\sigma^* (f))$. 
Since $A_\sigma^* (f)\in \mathbf{A}$, it follows from (F1) and Proposition~\ref{prop3.1}, 
that $c^*(A_\sigma^*( f))\in \mA$, and so, by definition, 
its Dirichlet series  $S_\sigma(D)$ belongs to $\mathscr{D}_{\mA}$. 
\end{proof}

\begin{proof}[Proof of Proposition~\ref{prop1.5}] 
Propositions \ref{prop3.1a} and \ref{prop1.4} imply that 
the correspondence $\sigma\mapsto S_\sigma$ 
determines a monomorphism of $S_\mN$ 
in the group of algebraic automorphisms of 
$\mA\in\{\mathscr H^\infty, \mathscr A, \mathscr W, \mathscr O_b\}$. 
Thus, it remains to show that each $S_\sigma$ is an isometry. 
   
Indeed, if $\mA=\mathscr{W}$, then since $S_\sigma$, $\sigma\in S_\mN$, 
rearranges the coefficients of the Dirichlet series, 
it is isometric by the  definition of the norm in $\mathscr W$, see \eqref{eq2.12}. 
   
Next, in the case $\mA\in\{\mathscr H^\infty, \mathscr A, \mathscr O_b\}$, 
by Theorem~\ref{te3.1}, 
\begin{equation}
\label{4.32}
c^*\circ A_\sigma^*=S_\sigma\circ c^*.
\end{equation}
From this, and noting that  $c^*$ is an isometry (see (F1), and Proposition \ref{prop3.1}), 
it suffices to check that $A_\sigma^*$ is an isometry.  
In turn, the norms on $H^\infty(B_{c_0})$ and $\mathcal A_u(B_{c_0})$ 
are  supremum norms over $B_{c_0}$, 
and since $A^*_\sigma$ maps $B_{c_0}$ biholomorphically onto itself, 
it preserves these supremum norms, 
which implies  the required statement (about isometry) for these two algebras. 
Finally, since $A_\sigma^*$ maps each $rB_{c_0}$, $r\in (0,1)$, 
biholomorphically onto itself, 
it preserves the seminorms defining $\mathcal O_b(B_{c_0})$ (see \eqref{eq3.20}), 
which gives the required statement for $\mathcal O_b(B_{c_0})$, 
and completes the proof. 
\end{proof}

%===========================
\section{Proof of Theorem \ref{te1.6}} 
%===========================

%--------------------------------------------------------------------------------------   
\subsection{$G$-invariant subalgebras of $\mathcal O(B_{c_0})$}
%-------------------------------------------------------------------------------------- 
   
Given a subgroup $G$ of $S_\mN$, we define 
the subalgebra of $\mathcal O(B_{c_0})$ of $G$-invariant holomorphic functions 
with respect to the action $A^*$:
$$
\mathcal O_G(B_{c_0})
:= \{f\in \mathcal O(B_{c_0})\, :\, f=A_{\sigma}^* (f)\textrm{ for all }\sigma \in G\}.
$$ 
Theorem~\ref{te3.1} implies that $f\in \mathcal O_G(B_{c_0})$ 
if and only if $\mathfrak O(f)$ is $G$-invariant with respect to the action $S|_G$, i.e., 
$$
S_\sigma(\mathfrak O(f)) = \mathfrak O(f)
\text{ for all }\sigma\in G.
$$
In particular, if $\mathfrak O(f)=\sum_{n=1}^\infty a_n n^{-s}$, 
then $f\in \mathcal O_G(B_{c_0})$ if and only if
\begin{equation}
\label{eq5.34}
a_{\hat\sigma(n)}=a_n\,\textrm{ for all }\, n\in\mN,\,\, \sigma\in G.
\end{equation}
Here $\hat\sigma$  is the completely multiplicative permutation defined in \eqref{eq1.8b}. 

Recall (see Section~\ref{subsection_3.1}, p.\pageref{pagenumber_mO_mOinfty}) 
that $\mO$ and $\mO_\infty=\mN\setminus \mO$ 
are the unions of all finite and infinite orbits, respectively, 
of the (natural) action of $G$ on $\mN$. 
We have the direct sum decomposition $c_0= c_{0,\mO}\oplus c_{0,\mO_\infty}$, where
\begin{eqnarray*}
c_{0,\mO}\!\!\!\!&=&\!\!\!\!\{\bbz=(z_1,z_2,\dots)\in c_0\, :\, z_n=0,\ n\in\mO_\infty\},
\\
c_{0,\mO_\infty}\!\!\!\!&=&\!\!\!\!\{\bbz=(z_1,z_2,\dots)\in c_0\, :\, z_n=0, \ n\in\mO\}.
\end{eqnarray*}
By the definition, $c_{0,\mO}$ and $c_{0,\mO_\infty}$ 
are $G$-invariant subspaces of $c_0$ with respect to the action $A|_G$, 
i.e., $A_\sigma(c_{0,\mO})=c_{0,\mO}$ and 
$A_\sigma(c_{0,\mO_\infty})=c_{0,\mO_\infty}$ for all $\sigma\in G$. 
The corresponding projections $\pi_\mO: c_0\to c_{0,\mO}$ and 
$\pi_{\mO_\infty}: c_0\to c_{0,\mO_\infty}$ 
are equivariant with respect to the action $A|_G$, and they map 
$B_{c_0}$ surjectively onto the open balls 
$B_{c_{0,\mO}}:=B_{c_0}\cap c_{0,\mO}$ and 
$B_{c_{0,\mO_\infty}}:=B_{c_0}\cap c_{0,\mO_\infty}$ 
of $c_{0,\mO}$ and $c_{0,\mO_\infty}$, respectively. 
Let $\mathcal O_G(B_{c_{0,\mO}})$ be 
the algebra of holomorphic functions on $B_{c_{0,\mO}}$ 
that are $G$-invariant with respect to $A^*|_G$. 

\begin{theorem}
\label{te3.3}
If $f\in \mathcal O_G(B_{c_0}),$ then
$$
f(\bbz'+\bbz) = f(\bbz) 
\textrm{ for all } \bbz'\in B_{c_{0,\mO_\infty}},
\ \bbz\in B_{c_0}  \textrm{ such that } \bbz'+\bbz\in B_{c_0}.
$$ 
Thus$,$ the pullback map $\pi_\mO^*,$ $f\mapsto f\circ\pi_\mO,$ 
maps $\mathcal O_G(B_{c_{0,\mO}})$ isomorphically onto $\mathcal O_G(B_{c_0})$.
\end{theorem}
\begin{proof}
First, we will prove the theorem for $f\in H^\infty(B_{c_0})\cap \mathcal O_G(B_{c_0})$,  
i.e., show that such functions belong to the range of $\pi_\mO^*$. 
In this case, $c^*(f)\in\mathscr H^\infty$ (see \cite{HedLinSei}), 
is the limit function of the Dirichlet series $\mathfrak O(f)$ (by Theorem~B), 
and is invariant with respect to the action $S|_G$ (see Theorem~\ref{te3.1}). 
So, if $\mathfrak O(f)=\sum_{n=1}^\infty a_n n^{-s}$, then 
$$
\lim_{n\to\infty}a_n=0
$$ 
(see, e.g., \cite[Prop.~1.19]{DGMS}), 
and $a_{\hat{\sigma}(n)}=a_n$ for all $n\in\mN$ and $\sigma\in G$ (see \eqref{eq5.34}). 
In particular, if $n$ has a divisor $p_k\in\mP$ with $k\in\mO_\infty$, 
then by the definition of an infinite orbit of the action of $G$ on $\mN$, 
there is a sequence $\{\sigma_i\}_{i\in\mN}\subset G$ such that 
$\lim_{i\to\infty}\sigma_i(p_k)=\infty$ 
which implies that $\lim_{i\to\infty}\hat{\sigma}_i(n)=\infty$ as well. Hence,
$$
\textstyle 
a_n=\lim\limits_{i\to\infty} a_{\hat\sigma_i(n)}=0.
$$
Thus, we have proved that
$$
\begin{array}{lrcl}
\mathfrak O(f) = \sum\limits_{n=1}^\infty a_n n^{-s},
\textrm{ where}& a_n\!\!\!&=&\!\!\!0 \ \
\textrm{ for all } n\in\mN\setminus\mN_\mO\, \textrm{ and}
\\[-0.2cm]
& \!\!\! a_{\hat{\sigma}(n)}\!\!\!&=&\!\!\!a_n\ \!
\textrm{ for all } n\in\mN_\mO,\ \sigma\in G.
\end{array}
$$
From this, by~Theorem A, we obtain that for $\bbz=(z_1,z_2,\dots)\in B_{c_{00}}$,
$$
\textstyle
f(\bbz)
= \mathfrak B(\mathfrak O(f))(\bbz)
= \sum\limits_{n\;\! \in \;\! \mN_\mO} a_n
   \Big(\prod\limits_{i\in \mO} z_i^{\nu_{p_i}(n)}\Big).
$$
Hence, $f|_{B_{c_{00}}}$ depends only on variables $z_i$ with $i\in\mO$. 
This implies 
$$
\;{\scaleobj{0.97}{
f(\bbz'\!+\!\bbz)
\!=\!f(\bbz)\textrm{ for all }\bbz' \!\in \!c_{00}\!\cap \!B_{c_{0,\mO_\infty}}
\textrm{and } \bbz\!\in\! B_{c_{00}}\textrm{ such that } \bbz'\!+\!\bbz\!\in\! B_{c_{00}}.}}
$$
Note that $c_{00}\cap B_{c_{0,\mO_\infty}}$ 
is dense in $B_{c_{0,\mO_\infty}}$, and $B_{c_{00}}$ is dense in $B_{c_0}$; 
hence by continuity of $f$, the previous identity is also valid  
for all $\bbz'\in  B_{c_{0,\mO_\infty}}$ 
and $\bbz\in B_{c_{0}}$ satisfying $\bbz'\!+\!\bbz\in B_{c_{0}}$. 
Finally, using the fact that the projections 
$\pi_\mO: c_0\to c_{0,\mO}$ and $\pi_{\mO_\infty} :c_0\to c_{0,\mO_\infty}$ 
map $B_{c_0}$ surjectively onto 
$B_{c_{0,\mO}}$ and $B_{c_{0,\mO_\infty}}$, respectively, 
the latter implies that
$$
f(\bbz)
= f(\pi_{\mO_\infty}(\bbz)+\pi_\mO(\bbz))
= f(\pi_\mO(\bbz)) 
\;\text{ for all }\; \bbz\in B_{c_0}.
$$
We define $\tilde f:=f|_{B_{c_{0,\mO}}}$. 
Then $\tilde f\in\mathcal O_G(B_{c_{0,\mO}})$ 
and $f=\pi_{\mO}^{*}(\tilde f)$, as required.

Now, let us prove the theorem for general $f\in\mathcal O_G(B_{c_0})$.

Since $f$ is continuous, it is bounded on an open neighbourhood of $\mathbf{0} \in c_0$, 
i.e., there is some $r_0\in (0,1)$ such that $f$ is bounded on $r_0 B_{c_0}$. 
Equivalently, $f_{{\texttt{r}}_0}:=f(r_0\cdot)\in H^\infty(B_{c_0})$. 
Clearly, $f_{{\texttt{r}}_0}\in \mathcal O_G(B_{c_0})$ 
and, hence, by the first part of the proof, 
$f_{{\texttt{r}}_0}(\bbz) = f_{\texttt{r}_0}(\pi_{\mO}(\bbz))$ for all $\bbz\in B_{c_0}$. 
The latter implies that the function $f-\pi_{\mO}^*(f)\in \mathcal O(B_{c_0})$ 
is equal to zero on $r_0 B_{c_0}$. 
Then, due to analyticity, $f-\pi_{\mO}^*(f)=0$, as required. 
This completes the proof of the theorem. 
\end{proof}

%---------------------------------------------------
\subsection{Proof of Theorem \ref{te1.6}} 
%---------------------------------------------------

Let $f\in \mathscr{O}_{b,G}$ be the limit function 
of a Dirichlet series $\sum_{n=1}^\infty a_n n^{-s}\in \mathscr D_{\mathscr O_b}$. 
Then, by the definition of $\mathscr  O_{b,G}$, see Section~\ref{subsection_3.1}, 
\begin{equation}
\label{e5.35}
a_{\hat\sigma(n)} = a_n
\ \textrm{ for all }\, n\in\mN,\ \sigma\in G. 
\end{equation}
By Proposition~\ref{prop3.1}, there exists an $F\in \mathcal O_b(B_{c_0})$ 
such that $c^*(F)=f$. Moreover, it follows from Theorem~\ref{te3.1} that  
$c^*\circ A_\sigma^*=S_\sigma\circ c^*$ for all $\sigma\in S_\mN$.  
Applying this to $F$, and using the fact that $S_\sigma (f)=f$ for all $\sigma \in G$, 
we obtain 
\begin{equation}
\label{eq_5.36}
c^*(A_\sigma^* (F)) = c^*(F) 
\text{ for all }\sigma \in G.
\end{equation}
As $A_\sigma^*(F)\in \mathcal O_b(B_{c_0})$ for all $\sigma \in S_\mN$, 
and since $c^* :\mathcal O_b(B_{c_0})\to \mathscr O_b$ 
is a Fr\'echet algebra isomorphism by Proposition~\ref{prop3.1}, 
equation \eqref{eq_5.36} implies that $A_\sigma^* (F)=F$ for all $\sigma \in G$. 
Hence, $F\in \mathcal O_G(B_{c_0})$, and so, according to Theorem~\ref{te3.3}, 
there exists an $\widetilde{F}\in \mathcal O_G (B_{c_{0,\mO}})$ 
such that $F=\pi_{\mO}^*( \widetilde{F})$. This shows that in \eqref{e5.35}, 
$a_n=0$ for all $ n\in (\mN_\mO)^{\com}$, as stated in the theorem. 

The converse is obvious, see \eqref{e5.35}, 
i.e., if $f\in \mathscr O_b$ satisfies conditions (i) and (ii) of the theorem, 
then it is $G$-invariant.   
\hfill$\Box$

%============================
\section{Proof of Theorem \ref{teo1.7}}
%============================

%-------------------------------------
\subsection{Auxiliary Results}
\label{sec7.1}
%-------------------------------------

In the proof of the theorem, we use the following notation and results.

For a subset $\mathcal P\subset\mP$, 
we denote by $\langle \mathcal P\rangle$ 
the unital multiplicative semigroup contained in $\mN$ 
generated by the elements of $\mathcal P$. 
For a Dirichlet series $D=\sum_{n=1}^\infty a_n n^{-s}$ 
and a subset $\mathcal P\subset\mP$, 
we let $D_\mathcal P$ be the series obtained from $D$ 
by removing all terms with $n\not\in \langle\mathcal P\rangle$, 
i.e., $D_\mathcal P:=\sum_{n\;\! \in\;\! \langle\mathcal P\rangle} a_n n^{-s}$. 

\begin{proposition}
\label{prop6.1}
Let $D\in\mathscr D_{\mathscr H^\infty}$ and $f\in\mathscr H^\infty$ 
be the limit function of $D$. Then for every subset $\mathcal P\subset\mP$ 
the Dirichlet series $D_\mathcal P$ has a limit function $f_\mathcal P\in \mathscr H^\infty$ 
and$,$ moreover$,$ $\|f_{\mathcal P}\|_{\mathscr H^\infty}\le \|f\|_{\mathscr H^\infty}$. 
Analogously$,$ if $\mathcal P_1\subset\mathcal P_2\subset\mP,$ 
then $\|f_{\mathcal P_1}\|_{\mathscr H^\infty}\le \|f_{\mathcal P_2}\|_{\mathscr H^\infty}$.
\end{proposition}
\begin{proof}
The proof follows easily from statement (F1) of Section~3 (see \cite{HedLinSei}). 
Indeed, according to this statement, 
there exists a unique element $F\in H^\infty(B_{c_0})$ 
such that $c^*(F)=f$ and $\|F\|_{H^\infty(B_{c_0})}=\|f\|_{\mathscr H^\infty}$. 
Let $c_{0,\mathcal P}\subset c_0$ be the subspace of vectors $\bbz=(z_1,z_2,\dots)\in c_0$ 
such that $z_n=0$ for all $n$ with $p_n\not\in\mathcal P$. 
Let $\pi_\mathcal P: c_0\to c_{0,\mathcal P}$ be the natural projection 
which preserves coordinates $z_n$ 
of a point $\bbz$ for all $n$ with $p_n\in  \mathcal P$, 
and sets all other coordinates to zero. 
Consider the pullback function 
$F_\mathcal P:=\pi_\mathcal P^*(F|_{c_{0,\mathcal P}})\, (:=F\circ \pi_\mathcal P)$. 
Then it is easy to see, using e.g., Theorem~A of Section~\ref{Section_3}, 
that $\mathfrak O(F_\mathcal P)=D_\mathcal P$. 
Since $F_\mathcal P\in H^\infty(B_{c_0})$, 
statement (F1) and Theorem~B then imply that 
$D_\mathcal P\in\mathscr D_{\mathscr H^\infty}$, 
and its limit function $f_\mathcal P$ coincides with $c^*(F_\mathcal P)$. 
In particular, by statement (F1), 
$\|f_\mathcal P\|_{\mathscr H^\infty}=\|F_\mathcal P\|_{H^\infty(B_{c_0})}$. 
In turn, by construction, 
$$
\|F_\mathcal P\|_{H^\infty(B_{c_0})}
= \sup_{B(c_0)\cap c_{0,\mathcal P}}|F|
\le \|F\|_{H^\infty(B(c_0))}
= \|f\|_{\mathscr H^\infty}.
$$
Combining this inequality with the previous identity, 
we get the first inequality of the proposition. 
The second inequality is proved similarly.
\end{proof}

Since the algebra $\mathscr D$ is isomorphic via the Bohr transform $\mathfrak B$ 
to the algebra of formal power series 
$\mC\llbracket x_1,x_2,\dots \rrbracket$, see  \eqref{eq_pg_3}, 
the correspondence $\Phi_\mathcal P:\mathscr D\to\mathscr D$, 
$\Phi_\mathcal P(D):=D_\mathcal P$, 
is an algebra epimorphism 
and a linear projection onto a subalgebra $\mathscr D_\mathcal P$ 
of Dirichlet series of the form $\sum_{n\;\! \in\;\! \langle\mathcal P\rangle} a_n n^{-s}$. 
The kernel of $\Phi_\mathcal P$ is an ideal $\mathscr I_\mathcal P$ of $\mathscr D$ 
consisting of series of the form $\sum_{n\;\! \not\in\;\! \langle\mathcal P\rangle} a_n n^{-s}$, 
so that $\mathscr D$ is a semidirect product of $\mathscr D_\mathcal P$ 
and $\mathscr I_\mathcal P$. As a corollary of Proposition \ref{prop6.1}, we get
  
\begin{corollary}
\label{cor6.2}
The restriction of the homomorphism $\Phi_{\mathcal P_{\mO_n}}$ 
to $\mathscr D_{\mathscr H^\infty}$ determines a Banach algebra epimorphism 
and a continuous projection of norm one from $\mathscr H^\infty$ 
onto a closed subalgebra $\mathscr H^\infty_{\mathcal P}$ 
of limit functions of Dirichlet series 
from $\mathscr D_\mathcal P\cap\mathscr D_{\mathscr H^\infty}$.
\end{corollary}

In what follows, we also denote this epimorphism by $\Phi_\mathcal P$.

The following result is derived straightforwardly from \cite[Cor.~3.10]{DGMS} 
by means of Proposition \ref{prop6.1}.

\begin{proposition}
\label{prop6.3}
Let $\mathcal P_1\subset\mathcal P_2\subset\cdots $ 
be a sequence of finite sets of primes 
such that ${\scaleobj{0.9}{\bigcup_{\;\!n=1}^{\;\!\infty}}}\mathcal P_n=\mP$. 
A Dirichlet series $D$ converges to $f\in \mathscr H^\infty$ 
if and only if for all $n\in \mN,$ the Dirichlet series $D_{\mathcal P_n}$ 
converges to an element $f_{\mathcal P_n}\in \mathscr H^\infty,$ 
and $\sup_{n\;\!\in\;\! \mN} \| f_{\mathcal P_n}\|_{\mathscr H^\infty} <\infty$. 
Moreover$,$ $\sup_{n\;\!\in\;\! \mN} \| f_{\mathcal P_n}\|_{\mathscr H^\infty} 
=\|f\|_{\mathscr H^\infty}$. 
\end{proposition}

%------------------------------------------------------
\subsection{Proof of Theorem~\ref{teo1.7}} 
\label{subsect6.2}
%------------------------------------------------------

We will first prove the theorem for the algebra $\mathscr H^\infty$. 

We retain the notation of Section~\ref{sec1.4a}. 
Recall that the set of (finite) orbits of the action of $G$ on $\mO\subset\mN$ 
is labelled as $O_1, O_2, O_3,\dots$ 
(for instance, such that  if $i<j$, then $\min O_i <\min O_j$). We set 
$$
\textstyle 
\mO_n := \bigcup\limits_{i=1}^n O_i
\;\text{ and }\;
\mathcal P_{\mO_n} := \{p_i\in\mP,\ i\in\mO_n\},\quad n\in\mN.
$$
Let $G_{\mO_n}$ be the normal subgroup of $G$ 
consisting of all $\sigma\in G$ acting identically on $\mO_n$, 
i.e., $\sigma|_{\mO_n}=\text{id}_{\mO_n}$. 
Then the action of $G$ on $\mO_n$ induces an {\em effective} action 
of the (finite) quotient group $G/G_{\mO_n}$ on $\mO_n$. 
In turn, the latter defines an action $\tilde S$ of $G/G_{\mO_n}$ 
on the algebra $\mathscr D_{\mathcal P_{\mO_n}}$ 
of the Dirichlet series of the form 
$D=\sum_{n\;\! \in\;\! \langle\mathcal P_{\mO_n}\rangle} a_n n^{-s}$ 
given for $[\sigma]:=\sigma G_{\mO_n}\in G/G_{\mO_n}$, $\sigma\in G$, by  
$$
\tilde S_{[\sigma]}(D) := S_{\sigma}(D).
$$
Here $S\in {\rm Hom}(G, {\rm Aut}(\mathscr D))$ 
is the homomorphism of Proposition \ref{prop3.1a}.
 
\begin{lemma}
\label{lem6.4}
Let $\pi_G:\mathscr D\to\mathscr D$ be the linear map given by equation \eqref{eq1.5}. 
If $D\in \mathscr D_{\mathcal P_{\mO_n}} ,$ then 
$$
\textstyle
\pi_G(D)=
{\scaleobj{1.2}{\frac{1}{|G/G_{\mO_n}|}}}
\sum\limits_{g\;\! \in\;\!  G/G_{\mO_n}}  \tilde S_{g}(D).
$$
\end{lemma}
 
\noindent (Recall that $|S|$ stands for the cardinality of a finite set $S$.) 
 
\begin{proof} 
For each $g\in G/G_{\mO_n}$, we fix $\sigma_g\in g$. 
Let $D=\sum_{k=1}^\infty a_k k^{-s}\in \mathscr D_{\calP_{\mO_n}}$.  
Let $H_k:=\{g\in G/G_{\mO_n}\, :\, \widehat{\sigma_g}^{-1}(k)=k\}$ 
be the stabiliser subgroup of the action of $G/G_{\mO_n}$ 
on the orbit $\textsf{O}(k)$ of $k$. 
If $g$ runs through  $G/G_{\mO_n}$, 
then $\widehat{\sigma_g}^{-1}(k)$ runs through $\textsf{O}(k)$, 
so that each term of the orbit occurs $|H_k|$ times. 
By definition, $|\textsf{O}(k)|=\frac{|G/G_{\mO_n}|}{|H_k|}$. Therefore 
$$
\phantom{AAa}
\begin{array}{rcl}
{\scaleobj{1.2}{\frac{1}{|G/G_{\mO_n}|}}}
\! \sum\limits_{g\;\! \in\;\!  G/G_{\mO_n}}\! \tilde S_{g}(D)
\!\!\!&=&\!\!\! 
\sum\limits_{k=1}^{\infty}\!\Big(
{\scaleobj{1.2}{\frac{1}{|G/G_{\mO_n}|}}}
\!\sum\limits_{g\;\! \in \;\! G/G_{\mO_n}}\!a_{\widehat{\sigma_g}^{-1}(k)}\Big) k^{-s}
\medskip
\\
\!\!\!&=&\!\!\! 
\sum\limits_{k=1}^{\infty}\!\Big(
{\scaleobj{1.2}{\frac{1}{|\textsf{O}(k)|}}}
\!\sum\limits_{m\;\! \in\;\!  \textsf{O}(k)}\!a_{m}
\Big) k^{-s}
=\pi_G (D).
\phantom{A}\qedhere
\end{array}
$$
\end{proof}

\noindent Let $f\in \mathscr H^\infty$ be the limit function 
of the Dirichlet series $D\!=\!\sum_{n=1}^\infty a_n n^{-s}$. 
Then $(\pi_{G}(f))_{\mathcal P_{\mO_n}} :=(\Phi_{\mathcal P_{\mO_n}} \circ \pi_G)(f)$ 
is the limit function of the Dirichlet series 
$(\pi_G(D))_{\mathcal P_{\mO_n}}:=(\Phi_{\mathcal P_{\mO_n}} \circ \pi_G)(D)$. 
Since by the definition of $\mO_n$, 
$$
\Phi_{\mathcal P_{\mO_n}} \circ \pi_G = \pi_G\circ \Phi_{\mathcal P_{\mO_n}},
$$
Lemma~\ref{lem6.4} implies that
$$
\textstyle
(\pi_G(D))_{\mathcal P_{\mO_n}}
=
{\scaleobj{1.2}{\frac{1}{|G/G_{\mO_n}|}}}
\sum\limits_{g\;\! \in \;\! G/G_{\mO_n}}  \tilde S_{g}(\Phi_{\mathcal P_{\mO_n}}(D)).
$$
From this, Proposition~\ref{prop1.5}, Corollary \ref{cor6.2}, 
and using the notation of the proof of Lemma \ref{lem6.4}, we obtain 
$$
\begin{array}{rcl}
\|(\pi_{G}(f))_{\mathcal P_{\mO_n}}\|_{\mathscr H^\infty} 
\!\!\!& \le &\!\!\!
{\scaleobj{1.2}{\frac{1}{|G/G_{\mO_n}|} }}
\sum\limits_{g\;\! \in \;\! G/G_{\mO_n}}  
\| S_{\sigma_g}(\Phi_{\mathcal P_{\mO_n}}(f)) \|_{\mathscr H^\infty}
\smallskip \\
\!\!\!&=&\!\!\!
\| \Phi_{\mathcal P_{\mO_n}} (f) \|_{\mathscr H^\infty}
\le \|f\|_{\mathscr H^\infty}.
\end{array}
$$
Note that $\pi_G(D)\in\mathscr D_{\mathcal P_\mO}$, where 
$$
\textstyle 
\mathcal P_\mO 
:= \{p_i\in\mP\, :\, i\in\mO\}
 = \bigcup\limits_{n\ge 1}\mathcal P_{\mO_n} . 
$$  
Hence, as  the previous inequality holds for all $\mO_n$, 
we can apply Proposition~\ref{prop6.3} to get 
$\pi_G( f)\in \mathscr H^\infty$ and 
$\|\pi_G( f)\|_{\mathscr H^\infty}\le \|f\|_{\mathscr H^\infty}$. 
Since 
\begin{itemize}[leftmargin=*]
\item[${\scaleobj{0.9}{\bullet}}$]
$\pi_G(f)$ depends linearly on $f\in\mathscr H^\infty$, 

\item[${\scaleobj{0.9}{\bullet}}$]
$\pi_G(f)\in \mathscr H_G^\infty$,  
 
\item[${\scaleobj{0.9}{\bullet}}$] 
$\pi_G( f)=f$ for all $f\in\mathscr H_G^\infty$ 
(by Theorem \ref{te1.6} and equation \eqref{eq1.5}),
\end{itemize}   
\noindent the latter inequality implies that 
$\pi_G:\mathscr H^\infty\to \mathscr H_G^\infty$ 
is a continuous linear projection of norm one. 

To complete the proof of the theorem for the algebra $\mathscr H^\infty$, 
we must check that $\pi_G(fg)=f\pi_G(g)$ for all $f\in\mathscr H_G^\infty$ 
and $g\in\mathscr H^\infty$. In fact, due to Proposition~\ref{prop6.3}, 
it suffices to check for such $f$ and $g$ that
$$
\Phi_{\mathcal P_{\mO_n}}(\pi_G(fg)) = \Phi_{\mathcal P_{\mO_n}}(f\pi_G(g))
$$
for all $\mO_n$. In turn, using the fact that $\Phi_{\mathcal P_{\mO_n}}$ 
is an algebra homomorphism and  
$\Phi_{\mathcal P_{\mO_n}} \circ \pi_G=\pi_G\circ \Phi_{\mathcal P_{\mO_n}}$, 
we obtain
\begin{eqnarray}
\nonumber 
&&\!\!\!
\Phi_{\mathcal P_{\mO_n}}(\pi_G(fg))-\Phi_{\mathcal P_{\mO_n}}(f\pi_G(g))
\\
&=&\!\!\!
\pi_G(\Phi_{\mathcal P_{\mO_n}}(f)\Phi_{\mathcal P_{\mO_n}}(g))
-\Phi_{\mathcal P_{\mO_n}}(f)\pi_G(\Phi_{\mathcal P_{\mO_n}}(g)).
\label{e6.37}
\end{eqnarray}
Next, applying Lemma~\ref{lem6.4}, and the two facts that 
\begin{itemize}[leftmargin=*]
\item[${\scaleobj{0.9}{\bullet}}$] 
$S_\sigma $ is an algebra homomorphism by  Proposition~\ref{prop1.5}, and  

\item[${\scaleobj{0.9}{\bullet}}$] 
$\Phi_{\mathcal P_{\mO_n}}(f)\in\mathscr H_G^\infty$ by the definition of $\mO_n$, 
\end{itemize} 
we get
$$
\begin{array}{rcl}
\pi_G(\Phi_{\mathcal P_{\mO_n}}(f)\;\!\Phi_{\mathcal P_{\mO_n}}(g))
\!\!\!&=&\!\!\!
{\scaleobj{1.2}{\frac{1}{|G/G_{\mO_n}|}}} 
\sum\limits_{g\;\! \in \;\! G/G_{\mO_n}}  
S_{\sigma_g}(\Phi_{\mathcal P_{\mO_n}}(f)
\Phi_{\mathcal P_{\mO_n}}(g))
\\
\!\!\!&=&\!\!\!
{\scaleobj{1.2}{\frac{1}{|G/G_{\mO_n}|}}}
\sum\limits_{g\;\! \in \;\! G/G_{\mO_n}}  
S_{\sigma_g}(\Phi_{\mathcal P_{\mO_n}}(f))
S_{\sigma_g}(\Phi_{\mathcal P_{\mO_n}}(g))
\\
\!\!\!&=&\!\!\!
{\scaleobj{1.2}{\frac{1}{|G/G_{\mO_n}|} }}
\sum\limits_{g\;\! \in \;\! G/G_{\mO_n}}  
\Phi_{\mathcal P_{\mO_n}}(f)
S_{\sigma_g}(\Phi_{\mathcal P_{\mO_n}}(g))
\\
\!\!\!&=&\!\!\! 
\Phi_{\mathcal P_{\mO_n}}(f)\;\! \pi_G(\Phi_{\mathcal P_{\mO_n}}(g)).
\phantom{{\scaleobj{1.2}{\frac{1}{|G/G_{\mO_n}|} }}}
\end{array}
$$
Thus the expression in \eqref{e6.37} is equal to zero, and we are done.

Having proved the theorem for $\mathscr H^\infty$, 
we now use this result to deduce the theorem for the other listed algebras. 

Let $f\in \mathscr O_{b}$ be the limit of the Dirichlet series $\sum_{n=1}^\infty a_n n^{-s}$.  
As $f\in \mathscr O_{b}$, it follows that for $r\in (0,1)$, $f_r \in \mathscr H^\infty$. 
(Recall that $f_r$ stands for the limit function 
of the Dirichlet series $\sum_{n=1}^\infty r^{\Omega(n)}a_n n^{-s}$.) 
Thus the function $\pi_G(f_r)\in \mathscr H^\infty_G$ 
by the previous part of the proof. But $\pi_G(f_r)$ 
is the limit of the Dirichlet series $\sum_{n=1}^\infty b_n n^{-s}$, where
$$
b_n
:=\left\{ \begin{array}{ccl}
0
&{\rm if}&n\in (\mN_\mO)^{\com},
\\[0.09cm]
{\scaleobj{1.2}{ \frac{1}{|{\textsf O}(n)|} }}
\sum\limits_{k\;\! \in\;\!  {\textsf O}(n)} r^{\Omega(k)} a_k
 &{\rm if}&n\in\mN_\mO.
\end{array} \right.
$$
As $\Omega(k)$ has the constant value $\Omega(n)$ for all $k\in {\textsf O}(n)$, 
it follows that $(\pi_G (f))_r=\pi_G(f_r)\in \mathscr H^\infty_G$. 
As $r\in (0,1)$ was arbitrary, we conclude that $\pi_G(f)\in \mathscr O_{b,G}$. 
Also if $g\in  \mathscr O_{b}$, then for all $r\in (0,1)$, $g_r \in \mathscr H^\infty$, 
and so we have $\pi_G(f_r g_r)=f_r\pi_G(g_r)$ by the previous part of the proof. 
Thus, using the fact that  the correspondence $r\mapsto \cdot_r$ 
determines a monomorphism 
of the multiplicative group $\mC^*$ in the group ${\rm Aut}(\mathscr D)$, 
we obtain 
$$
\begin{array}{rcl}
(\pi_G (fg))_r 
\!\!\!&=&\!\!\! 
\pi_G ((fg)_r) 
= \pi_G(f_r g_r) 
= f_r \pi_G (g_r) 
= f_r (\pi_G(g))_r 
\\
\!\!\!&=&\!\!\! 
(f(\pi_G(g)))_r.
\end{array}
$$
It follows from here that $\pi_G (fg)=f(\pi_G(g))$, as wanted. 
Also, for $f\in \mathscr O_{b}$, and for $r\in (0,1)$, 
we have by the previous part of the proof,
$$
P_r(\pi_G(f))
:= \|(\pi_G(f))_r\|_{\mathscr H^\infty}
= \|(\pi_G(f_r))\|_{\mathscr H^\infty}
\leq \|f_r\|_{\mathscr H^\infty}
=: P_r(f).
$$
This implies that $\pi_G:\mathscr O_{b}\to\mathscr O_{b,G}$ 
is a continuous linear projection of norm one, 
and completes the proof of the theorem for the algebra $\mathscr O_b$.

Next, for the algebra $\mathscr A$, 
we use the fact that 
it is the uniform closure of the algebra $\mathscr W$, see \cite[Thm.~2.3]{ABG}. 
Thus, given  $f\in \mathscr A$ there is a sequence $\{f_n\}_{n=1}^\infty$ of Dirichlet series 
with finitely many terms such that $\lim_{n\to\infty}\|f-f_n\|_{\mathscr H^\infty}\!=\!0$. 
By the definition of the map $\pi_G$, the sequence $\{\pi_G(f_n)\}_{n=1}^\infty$ 
consists also of Dirichlet series with finitely many terms, 
and moreover, due to the proved part of the theorem for the algebra $\mathscr H^\infty$,
$$
\|\pi_G(f)-\pi_G(f_n)\|_{\mathscr H^\infty}
\le \|f-f_n\|_{\mathscr H^\infty}.
$$
These imply that the sequence $\{\pi_G(f_n)\}_{n=1}^\infty$ 
converges uniformly on $\overline{\mC}_+$ 
to the function $\pi_G(f)\in\mathscr H_G^\infty$. 
Hence, in fact, $\pi_G(f)\in \mathscr A_G$, as required. 
Moreover, $\|\pi_G(f)\|_{\mathscr A}\le\|f\|_{\mathscr A}$ 
for all $f \in \mathscr A$, and $\pi_G (fg)=f(\pi_G(g))$ for all $f\in \mathscr A_G$ 
and $g\in \mathscr A$, since $\mathscr A$ is a uniform subalgebra of $\mathscr H^\infty$, 
for which these facts are true by the first part of the proof of the theorem. 
  
Finally, for $\mathscr W$, the proof is straightforward, 
as we can rearrange terms arbitrarily in an absolutely converging series: 
If  $f\!\in\! \mathscr W$ is the limit of the Dirichlet series 
$D = \sum_{n=1}^\infty a_nn^{-s}$, 
then  $\|\pi_G (f)\|_{\mathscr W}\!\le\!\|f\|_{\mathscr W}$ 
because 
$$
\!\!
\begin{array}{rcl}
\|\pi_G(f)\|_{\mathscr W}
\!\!\!&=&\!\!\!
\sum\limits_{n\;\!\in\;\! \mN_\mO} 
{\scaleobj{1.2}{\frac{1}{|{\textsf O}(n)|} }}
\Big|\sum\limits_{k\;\!\in\;\! {\textsf O}(n)} a_k\Big| 
\\
\!\!\!&\le &\!\!\!
\sum\limits_{n\;\!\in \;\!\mN_\mO} 
{\scaleobj{1.2}{\frac{1}{|{\textsf O}(n)|} }}
\sum\limits_{k\;\!\in \;\!{\textsf O}(n)} |a_k|
= 
\sum\limits_{n\;\!\in\;\! \mN_\mO} 
\sum\limits_{k\;\!\in\;\! {\textsf O}(n)} |a_k| 
\le
\sum\limits_{n\;\!\in\;\! \mN} |a_n| 
=
\|f\|_{\mathscr W}.
\end{array}
$$
So $\pi_G (f)\in \mathscr W_G$, and the norm of $\pi_G$ is equal to $1$. 
Also, we have $\pi_G (fg)=f(\pi_G(g))$ for $f\in \mathscr W_G$ 
and $g\in \mathscr W$, since this holds 
when $f,g$ are viewed as elements of $\mathscr A$. 
This completes the proof for the algebra $\mathscr W$, 
and so of the theorem. 
\hfill$\Box$

%==============================================
\section{Proofs of Propositions~\ref{prop1.8} and \ref{prop1.9}}
\label{section_prop1.8_proof}
%==============================================

\begin{proof}[Proof of Proposition~\ref{prop1.8}]
We have to prove that the map \penalty-10000 
$Q_\mO\times ({\scaleobj{0.9}{\overline{\mD}^{\mO}}},\tau_p) 
\to ({\scaleobj{0.9}{\overline{\mD}^{\mO}}},\tau_p)$, 
$(\textsl{q} ,f)\mapsto \psi^*(\textsl{q})(f)$ is continuous. 
Since both spaces there are at most countable direct products of metrisable spaces, 
they are also metrisable. Thus, it suffices to prove that if the sequence 
$\{(\textsl{q}_n,f_n)\}_{n\in \mN}
\subset Q_\mO\times ({\scaleobj{0.9}{\overline{\mD}^{\mO}}},\tau_p)$ converges to 
$(\textsl{q},f)\in Q_\mO\times ({\scaleobj{0.9}{\overline{\mD}^{\mO}}},\tau_p)$, 
then the sequence $\{\psi^*(\textsl{q}_n)(f_n)\}_{n\in\mN}$ 
converges to $\psi^*(\textsl{q})(f)$ in $({\scaleobj{0.9}{\overline{\mD}^{\mO}}},\tau_p)$. 
 
We have $\textsl{q}_n=(q_{1,n}, q_{2,n},…)$ and $\textsl{q}=(q_1,q_2,…)$, 
where $q_{i,n}, q_i\in S_{O_i}$, $i=1,2,\dots$, and $n\in\mN$. 
Also, we set $f_{i,n}:=f_n|_{O_i}$ and 
$f_i:=f|_{O_i}$ for all $i=1,2,\dots$, and $n\in\mN$. 
By the definition of the (product) topologies 
on $Q_\mO$ and ${\scaleobj{0.9}{\overline{\mD}^{\mO}}}$,
$$
\begin{array}{ll}
\lim\limits_{n\to\infty}\psi^*(\textsl{q}_n)(f_n)
= 
\psi^*(\textsl{q})(f)
\text{ if and only if} 
\\
\lim\limits_{n\to\infty} f_{i,n}\circ q_{i,n}
=f_i\circ q_i 
\text{ for all }i=1,2,\dots.
\end{array}
$$
In turn, $\lim\limits_{n\to\infty} \textsl{q}_n = \textsl{q}$ 
means that for every $i$, $q_{i,n}=q_i$ for all sufficiently large $n$. 
Thus, for every $i$,
$$
\lim\limits_{n\to\infty}\psi^*(q_{i,n})(f_{i,n})
= \lim\limits_{n\to\infty}\psi^*(q_i)(f_{i,n}) 
= \lim\limits_{n\to\infty}f_{i,n}\circ q_i 
=f_{i}\circ q_i.
$$
The  last equality holds because $\{f_{i,n}\}_{n\in\mN}$ converges pointwise to $f_i$.
\end{proof}
 
\begin{proof}[Proof of Proposition~\ref{prop1.9}]
We have to prove that the map \penalty-10000 
$Q_\mO\times ( {\scaleobj{0.9}{\mD^{\mO}}},\tau_{hk}) 
\to ( {\scaleobj{0.9}{\mD^{\mO}}},\tau_{hk})$, 
$(\textsl{q} ,f)\mapsto \psi^*(\textsl{q})(f)$, is continuous. 
Since the space $({\scaleobj{0.9}{ \mD^{\mO}}},\tau_{hk})$ is Hausdorff 
and ${\scaleobj{0.9}{\mD^\mO}}/{\scaleobj{0.9}{\widehat G_\mO}}$ is hemicompact, 
it is sufficient to prove that the restriction of the above map 
to each set $Q_\mO\times r{\scaleobj{0.9}{\overline{\mD}^\mO}}$, $r\in (0,1)$, 
is continuous, see, e.g., \cite[Prop.~2.7]{Ste}. 
However, the topologies $\tau_{p}$ and $\tau_{hk}$ 
coincide on each $r{\scaleobj{0.9}{\overline{\mD}^\mO}}$, 
and so the continuity of the restricted map follows from Proposition~\ref{prop1.8}.
\end{proof}

%=================================================
\section{Proofs of Theorem~\ref{te1.9} and Corollary~\ref{cor1.10}}
%=================================================
 
Recall that $c_{0,\mO}$ is a closed subspace 
of points $\bbz=(z_1,z_2,\dots)\in c_0$ 
with $z_n=0$ for all $n\in\mO_\infty$, and $\pi_\mO: c_0\to c_{0,\mO}$ 
is the natural projection which preserves 
the coordinates $z_n$ of a point $\bbz$ for all $n\in\mO$, 
and sets all other coordinates to zero. 
The pullback map $\pi_\mO^*$, $f\mapsto f\circ\pi_\mO$, 
maps $\mathcal O(B_{c_{0,\mO}})$ onto a subspace of $\mathcal O(B_{c_0})$ 
of holomorphic functions constant along the fibres of the projection $\pi_\mO$. 
Then it is easy to see using Theorem~A of Section~\ref{Section_3}, 
that $\mathfrak O(\pi_\mO^*(f))\subset\mathscr D[\mO]$, 
the subalgebra of Dirichlet series $D\!=\!\sum_{n=1}^\infty a_n n^{-s}$ 
with $a_n\!=\!0$ for all $n\!\in\!\mO_\infty\!:=\! \mO^{\com}$. 
Also, Theorem~B of Section~\ref{Section_3} 
implies that the pullback map $c^*$ maps $\pi_\mO^*(f)$ 
to a holomorphic function in $\mathscr O_u$ 
whose Dirichlet series is $\mathfrak O(\pi_\mO^*(f))$. 
By definition, $c^*(\pi_\mO^*(f))=(\pi_\mO\circ c)^*(f)$. 
Here the continuous map $\tilde c_{\mO}:=\pi_\mO\circ c:\mC_+\to B(c_{0,\mO})$ 
is given by the formula
\begin{equation}
\label{eq8.38}
\begin{array}{l}
\tilde c_{\mO}(s) = ((\tilde c_{\mO}(s))_1,(\tilde c_{\mO}(s))_2,\dots ),
\quad s\in\mC_+,\quad {\rm where}\medskip
\\
(\tilde c_{\mO}(s))_i 
=
\left\{ \begin{array}{cl}
p_i^{-s} & \text{if }\;i\in\mO\smallskip \\
0&{\rm otherwise}.
\end{array}\right.
\end{array}
\end{equation}
 
Let us now consider $\ell^\infty=\ell_1^*$ as the second dual of $c_0$, 
so that $c_0$ is isometrically embedded in $\ell^\infty$. 
The weak$^*$ closure of $c_{0,\mO}$ in $\ell^\infty$ 
is the subspace $\ell^\infty_\mO$ of points $\bbz\!=\!(z_1,z_2,\dots)\!\in\! \ell^\infty$ 
with $z_n\!=\!0$ for all $n\!\in\!\mO_\infty$. 
The projection $\pi_\mO$ extends to a continuous and weak$^*$ continuous projection 
from $\ell^\infty$ onto $\ell^\infty_\mO$, which we denote also by $\pi_\mO$. 
 
Let $H^\infty(B_{c_{0,\mO}})$ be the complex Banach algebra 
of bounded holomorphic functions on $B_{c_{0,\mO}}$, 
equipped with the supremum norm, and 
$\calA_u(B_{c_{0,\mO}})$ be the subalgebra of functions 
from $H^\infty(B_{c_{0,\mO}})$ that are uniformly continuous on $B_{c_{0,\mO}}$. 
Then the pullback $\pi_\mO^*$ maps both algebras 
isometrically onto closed subalgebras of $H^\infty(B_{c_0 })$ 
and $\calA_u(B_{c_0 })$, respectively, 
of holomorphic functions constant along the fibres of $\pi_\mO$. 
Applying to these subalgebras statement (F1) of Section~\ref{Section_3} 
along with Corollary~\ref{cor6.2}, we obtain that 
{\em the pullback homomorphism $(\tilde c_\mO)^*,$ 
with $\tilde c_\mO$ given by \eqref{eq8.38}$,$ 
maps $H^\infty(B_{c_{0,\mO}})$ isometrically onto $\mathscr{H}^\infty [\mO],$ 
and $\calA_u(B_{c_{0,\mO}})$ isometrically onto $\mathscr{A}[\mO]$}. 
(Recall that for 
$\mA\in\{\mathscr O_b,\mathscr{W}, \mathscr{A}, \mathscr{H}^\infty\}$, 
we denote by $\mA[\mO]\subset\mA$ 
the closed subalgebra of limit functions 
of Dirichlet series in $\mathscr D[\mO]\cap\mathscr D_\mA$.) 

Similarly, let $W(B_{c_{0,\mO}})\subset \calA_u(B_{c_{0,\mO}})$ 
consist of all functions $f$ such that $\pi_\mO^*(f)\in W(B_{c_0})$, 
with norm $\|f\|_{W(B_{c_{0,\mO}})}:=\|\pi_\mO^*(f)\|_{W(B_{c_0})}$. 
Then {\em $W(B_{c_{0,\mO}})$ is a complex Banach algebra, 
isometrically isomorphic to the algebra $\mathscr W[\mO]$ via $(\tilde c_\mO)^*$}.
 
Finally, let $\mathcal O_b(B_{c_{0,\mO}})\subset\mathcal O(B_{c_{0,\mO}})$ 
be the Fr\'echet algebra of holomorphic functions 
which are bounded on each ball $rB_{c_{0,\mO}}$, $0<r<1$, 
equipped with the topology defined by the system of seminorms 
\begin{equation}
\label{eq3.20b}
P_r(f):=\sup_{\bbz\;\!\in \;\!B_{c_{0,\mO}}} |f(r\bbz)|, 
\quad 0<r<1.
\end{equation}
Then the above arguments together with Proposition~\ref{prop3.1} 
imply that {\em the homomorphism $(\tilde c_\mO)^*$ maps 
$\mathcal O_{b}(B_{c_{0,\mO}})$ onto $\mathscr O_b[\mO]$ 
isometrically with respect to the chosen systems of seminorms on these Fr\'echet spaces}.

Further, each function from $\pi_\mO^*(H^\infty(B_{c_{0,\mO}}))\subset H^\infty(B_{c_{0}})$ 
admits a weak$^*$ continuous extension to $\mD^\infty$, see Section~\ref{Section_3}. 
The set of the extended functions forms a closed subalgebra of $H_{w}^\infty(\mD^\infty)$ 
of holomorphic functions constant along the fibres of $\pi_\mO$. 
This implies that each function from $H^\infty(B_{c_{0,\mO}})$ 
admits a weak$^*$ continuous extension to the open ball $B_{\ell^\infty_\mO}$ 
of the space $\ell^\infty_\mO$. We denote by $H_w^{\infty}(B_{\ell^\infty_\mO})$ 
the algebra of the extended functions. 
Clearly, $\pi_\mO^*(H_w^{\infty}(B_{\ell^\infty_\mO}))\subset H_{w}^\infty(\mD^\infty)$, 
and {\em the restriction to $B_{c_{0,\mO}}$ 
determines an isometric isomorphism 
between $H^{\infty}_w(B_{\ell^\infty_\mO})$ and $H^\infty(B_{c_{0,\mO}})$}.

Let $A(B_{\ell^\infty_\mO})$ and $W(B_{\ell^\infty_\mO})$ 
be subalgebras of $H^{\infty}_w(B_{\ell^\infty_\mO})$ 
of functions whose pullbacks by $\pi_\mO$ belong to $A(\mD^\infty)$ and $W(\mD^\infty)$, 
equipped with the supremum norm and the norm induced from $W(\mD^\infty)$, respectively. 
Using statements (F2) and (F3) of Section~\ref{Section_3}, 
we obtain that {\em the restriction to $B_{c_{0,\mO}}$ 
determines a Banach algebra isometric isomorphism 
between $A(B_{\ell^\infty_\mO})$ and $\calA_u(B_{c_0, \mO}),$ 
and between $W(B_{\ell^\infty_\mO})$ and  $W(B_{c_{0,\mO}}),$ respectively}.

Similarly, each function from $\mathcal O_b(B_{c_0,\mO})$ 
admits a weak$^*$ continuous extension to $B_{\ell^\infty_\mO}$. 
By $\mathcal O_{w}(B_{\ell^\infty_\mO})$ 
we denote the algebra of the extended functions, 
equipped with the metrisable locally convex topology 
defined by the system of seminorms
$$
P_r( f):=\sup_{\bbz\;\!\in \;\!B_{\ell^\infty_\mO}}| f(r\bbz)|, 
\quad 0<r<1.
$$
Then $\mathcal O_{w}(B_{\ell^\infty_\mO})$ is a Fr\'echet algebra, 
and {\em the restriction to $B_{c_{0,\mO}}$ determines a 
Fr\'echet algebra isomorphism between 
$\mathcal O_{w}(B_{\ell^\infty_\mO})$ and $\mathcal O_b(B_{c_0,\mO}),$  
isometric  with respect to the chosen systems of seminorms on these Fr\'echet spaces}.

Let us consider  $\tilde c_\mO$ as a map 
from $\mC_+$ into $B_{\ell^\infty_\mO}$. 
Then combining all the above facts, we get: 

\begin{proposition}
\label{prop8.1}
The extended homomorphism $(\tilde c_\mO)^*$ 
maps the algebras 
$\mathcal O_{w}(B_{\ell^\infty_\mO})$, 
$H^\infty_w(B_{\ell^\infty_\mO}),$ 
$A(B_{\ell^\infty_\mO}),$ and 
$W(B_{\ell^\infty_\mO})$ 
isometrically onto the algebras 
$\mathscr O_b[\mO],$ 
$\mathscr H^\infty[\mO],$ 
$\mathscr A[\mO],$ and 
$\mathscr{W}[\mO],$ respectively.
\end{proposition}

Next, recall that $A$ is a faithful representation of $S_\mN$ on  $\ell^\infty$ 
by linear isometries $A_\sigma$, $\sigma\in S_\mN$, defined by the formula
$$
A_\sigma(\bbz):=(z_{\sigma(1)},z_{\sigma(2)},\dots)
\quad{\rm for}\quad \bbz=(z_1,z_2,\dots)\in\ell^\infty,
$$
see Section~\ref{SubSect_1.3}.

Also, recall that $S$ is a monomorphism of $S_\mN$ in the group ${\rm Aut}(\mathscr D)$, 
given by the formula
$$
\textstyle 
S_{\sigma}(D):=\sum\limits_{n=1}^\infty a_{\hat\sigma^{-1}(n)} n^{-s}
\quad {\rm for}\quad  D=\sum\limits_{n=1}^\infty a_n n^{-s},
$$
see Proposition~\ref{prop3.1a}, 
and the correspondence  $\sigma\mapsto S_\sigma$ 
determines a monomorphism of $S_\mN$ in the group ${\rm Aut_I}(\mA)$ 
of  isometric algebra isomorphisms of 
$\mA\in \{ \mathscr{O}_b,\mathscr{W}, \mathscr{A}, \mathscr{H}^\infty\}$, 
denoted by the same symbol, see Proposition~\ref{prop1.5}.

Let $H \subset S_{\mN}$ be a subgroup of permutations which maps $\mO$ into itself. 
Then the ball $B_{\ell^\infty_\mO}$ is invariant 
with respect to the representation $A|_H$, and the algebras 
$\mathcal A\in 
\{\mathcal O_{w}(B_{\ell^\infty_\mO}), H^\infty_w(B_{\ell^\infty_\mO}), 
A(B_{\ell^\infty_\mO}), W(B_{\ell^\infty_\mO})\}$
and 
$\mA[\mO]
\in \{ \mathscr{O}_b[\mO], \mathscr{H}^\infty[\mO], 
\mathscr{A}[\mO], \mathscr{W}[\mO]\}$ 
are invariant with respect to the actions 
$A^*|_H$ and $S|_H$, respectively. 
Moreover, Theorem \ref{te3.1} implies easily that for each $f\in\mathcal A$,
\begin{equation}
\label{eq8.40}
( (\tilde c_\mO)^* \circ A^*_\sigma)(f)= (S_{\sigma} \circ  (\tilde c_\mO)^*)(f)
\,\textrm{ for all }\, \sigma\in H.
\end{equation}

Now we are ready to prove Theorem~\ref{te1.9}.

\begin{proof}
Let us define the isometric isomorphism $\iota:\ell^\infty_\mO\to \ell^\infty(\mO)$ by 
$$
(\iota (\bbz))(i) = z_i 
\text{ for all } i\in \mO, \text{ and all } \bbz=(z_1,z_2,\dots)\in \ell^\infty_\mO.
$$
Then it is easy to see that the pullback $\iota^*$, $f\mapsto f\circ \iota$, 
is an isometric isomorphism from 
$\mathcal O_w({\scaleobj{0.9}{\mD^\mO}}), H_w^\infty({\scaleobj{0.9}{\mD^\mO}}), 
A({\scaleobj{0.9}{\mD^\mO}}), W({\scaleobj{0.9}{\mD^\mO}})$ 
onto the algebras 
$\mathcal O_{w}(B_{\ell^\infty_\mO}),
H^\infty_w(B_{\ell^\infty_\mO})$, 
$A(B_{\ell^\infty_\mO}),
W(B_{\ell^\infty_\mO})$, respectively. 
(In the case of $\mathcal O_{w}(B_{\ell^\infty_\mO})$ 
and $O_w({\scaleobj{0.9}{\mD^\mO}})$, 
the isometry means that for each $r\in (0,1)$, the seminorm $P_r$ is preserved.)  
We have $\iota(B_{\ell^\infty_\mO})\!=\!{\scaleobj{0.9}{\mD^\mO}}$, 
$\iota \circ \tilde{c}_\mO\!=\!c_\mO$, and so $c_\mO^*=\tilde{c}_\mO^* \circ \iota^*$. Proposition~\ref{prop8.1} now implies that 
if $(\mathcal A,\mA)$ is one of the pairs 
$(\mathcal O_w({\scaleobj{0.9}{\mD^\mO}}),\mathscr O_b[\mO])$, 
$(H_w^\infty({\scaleobj{0.9}{\mD^\mO}}), \mathscr H^\infty [\mO])$, 
$(A({\scaleobj{0.9}{\mD^\mO}}), \mathscr A[\mO])$, 
$(W({\scaleobj{0.9}{\mD^\mO}}),\mathscr W[\mO])$, 
then $c_\mO^*:\mathcal A\to\mA$ is an isometric isomorphism of algebras.

Further, $S_\mO$ is naturally regarded as a subgroup of the group $S_\mN$ 
consisting of all permutations that fix points of $\mO_\infty:=\mN\setminus\mO$. 
Then $\psi(Q_\mO)$ is a subgroup of $S_\mO$ (see Section~\ref{sec1.4a}), 
and we can apply equation \eqref{eq8.40} to $H:=\psi(Q_\mO)$, 
observing that $\iota\circ A_{\psi(\textsl{q})}=\psi^*(\textsl{q})\circ\iota$,  
to get for all $\textsl{q}\in Q_\mO$ and $f\in\mathcal A$:
$$
\begin{array}{rcl}
(c_\mO^* \circ \Psi_{\textsl{q}})(f)
\!\!\!&=&\!\!\! 
(\tilde{c}_\mO^* \circ \iota^* \circ (\psi^*(\textsl{q}))^*)(f)
= (\psi^*(\textsl{q})\circ\iota\circ\tilde{c}_\mO)^*(f)
\\
\!\!\!&=&\!\!\!
(\iota\circ A_{\psi(\textsl{q})}\circ\tilde{c}_\mO)^*(f) 
= ( (\tilde c_\mO)^* \circ A^*_{\psi(\textsl{q})})(\iota^*(f))
\\
\!\!\!&=&\!\!\! 
(S_{\psi(\textsl{q})} \circ  (\tilde c_\mO)^*\circ\iota^*)(f)
= (S_{\psi(\textsl{q})} \circ c_{\mO}^*)(f).
\end{array}
$$
This gives the equivariance relation of Theorem \ref{te1.9}, and completes the proof.
\end{proof}
 
\begin{proof}[Proof of Corollary~\ref{cor1.10}] 
This follows directly from Proposition~\ref{prop1.8} and Theorem~\ref{te1.9}. 
We leave the details to the reader. 
\end{proof}
   
\begin{proof}[Proof of Proposition~\ref{prop1.11}]
It is sufficient to prove the proposition 
for the algebra $\mathcal O_w({\scaleobj{0.9}{\mD^\mO}})$ 
because it contains as subalgebras all other algebras of the proposition. 
Furthermore, since the operators on the left and on the right in the formula of the proposition 
are continuous with respect to the topology of $\mathcal O_w({\scaleobj{0.9}{\mD^\mO}})$ 
(which is the topology of the uniform convergence 
on subsets $r{\scaleobj{0.9}{\overline{\mD}^\mO}}$ for all $r\in (0,1)$), 
and the algebra of polynomials in coordinates $z_i$, $i\in\mO$, 
of $\bbz\in\ell^\infty(\mO)$ is dense in $\mathcal O_w({\scaleobj{0.9}{\mD^\mO}})$, 
it is sufficient to prove the result for monomials. 

Let us consider a monomial $f(\bbz):=z_{i_1}^{\alpha_1}\cdots z_{i_k}^{\alpha_k}$, 
$\bbz\in\ell^\infty(\mO)$, $k\in \mN$,  $i_1,\dots ,i_k\in  \mO$, 
and $\alpha_1,\dots,\alpha_k \in \mN$. Applying the map $c_\mO^*$, 
we get $c_\mO^*(f)=(p_{i_1}^{\alpha_1}\cdots p_{i_k}^{\alpha_k})^{-s}=:D$. 
Let $p_{i_1},\dots ,p_{i_k}\in \mathcal P_{\mO_n}:=\{p_i\in\mP,\ i\in\mO_n\}$, 
where $\mO_n:={\scaleobj{0.9}{\bigcup_{\;\!i=1}^{\;\! n}}}  O_i$ 
(we use the notation of Section~\ref{subsect6.2}). 
Let $G_{\mO_n}$ be the normal subgroup of $G$ 
consisting of all $\sigma\in G$ acting identically on $\mO_n$. 
By Lemma~\ref{lem6.4}, 
$$
\textstyle
(\pi_G|_{\mA} \circ c_{\mO}^*) (f)
= {\scaleobj{1.2}{\frac{1}{|G/G_{\mO_n}|}}} 
  \sum\limits_{g\;\! \in\;\!  G/G_{\mO_n}} \tilde{S}_g(D).
$$
Here if $g=\sigma G_{\mO_n}$ for an element $\sigma \in G$, 
then $\tilde{S}_g(D)=S_\sigma (D)$, where $S\in \text{Hom}(G, \text{Aut}(\mathscr D))$ 
is the homomorphism of Proposition~\ref{prop3.1a}. 
Let ${\scaleobj{0.9}{\widehat G_{\mO_n}}}$ be the closure in $Q_\mO$ 
of the group $\rho(G_{\mO_n})$ (see equation \eqref{eq1.6} for the definition of $\rho$). 
Then ${\scaleobj{0.9}{\widehat G_\mO}}/{\scaleobj{0.9}{\widehat G_{\mO_n}}}$ 
is the same as $G/G_{\mO_n}$ (and in particular, finite).  
For each $g\in {\scaleobj{0.9}{\widehat G_{\mO}}}/{\scaleobj{0.9}{\widehat G_{\mO_n}}}$, 
we fix $\textsl{q\,}_{\!g}\in g$. Using Theorem~\ref{te1.9}, we have
$$
\textstyle
\sum\limits_{g\;\! \in\;\!  G/G_{\mO_n}} (c_\mO^*)^{-1}  \tilde{S}_g(D)
=
\sum\limits_{g\;\! \in\;\!  {\scaleobj{0.9}{\widehat G_\mO}}/{\scaleobj{0.9}{\widehat G_{\mO_n}}}} \Psi_{\textsl{q\,}_{\!g}}(f).
$$
So the left-hand side in the claim of the proposition can be written as
\begin{equation}
\label{e8.41}
\textstyle
P_{{\scaleobj{0.9}{\widehat G_\mO}}}(f)
= ( (c_\mO^*)^{-1} \circ \pi_G|_{\mA} \circ c_{\mO}^*)(f)
= {\scaleobj{1.2}{\frac{1}{| G/G_{\mO_n} |}}}
\sum\limits_{ g\;\! \in\;\!  {\scaleobj{0.9}{\widehat G_\mO}}/{\scaleobj{0.9}{\widehat G_{\mO_n}}} } \Psi_{\textsl{q\,}_{\!g}}(f).
\end{equation}    
 
Next, the kernel of the quotient map 
${\scaleobj{0.9}{\widehat G_\mO}} \to 
{\scaleobj{0.9}{\widehat G_\mO}}/{\scaleobj{0.9}{\widehat G_{\mO_n}}}$, 
namely the group ${\scaleobj{0.9}{\widehat G_{\mO_n}}}$, 
is a clopen subgroup of ${\scaleobj{0.9}{\widehat{G}_{\mO}}}$, 
and in particular, it is $m_{{\scaleobj{0.9}{\widehat G_\mO}}}$ measurable. 
Since $m_{{\scaleobj{0.9}{\widehat G_\mO}}}$ is normalised, 
$m_{{\scaleobj{0.9}{\widehat G_\mO}}}(\textsl{q\,}{\scaleobj{0.9}{\widehat G_{\mO_n}}})
= \frac{1}{|G/G_{\mO_n}|}$ for all $\textsl{q}\in {\scaleobj{0.9}{\widehat G_{\mO}}}$. 
Also, by the definition of the action $\Psi$, see \eqref{eq1.7}, 
the monomial $f$  satisfies 
$\Psi_{\textsl{q}\textsl{q}'}(f)=\Psi_{\textsl{q}}(f)$ 
for all $\textsl{q}\in {\scaleobj{0.9}{\widehat G_{\mO}}}$, 
$\textsl{q\;\!}'\in {\scaleobj{0.9}{\widehat G_{\mO_n}}}$. 
Then 
$$
\begin{array}{rcl}
\displaystyle 
{\scaleobj{0.9}{\int\limits_{{\scaleobj{0.9}{\widehat{G}_\mO}}}}}\;
 (\Psi_{\textsl{q}}(f))(\bbz) \;\!\text{d}m_{{\scaleobj{0.9}{\widehat{G}_\mO}}}(\textsl{q})
\!\!\!&=&\!\!\!
\displaystyle 
{\scaleobj{0.9}{
\sum\limits_{g\;\! \in \;\! {\scaleobj{0.9}{\widehat G_\mO}}/{\scaleobj{0.9}{\widehat G_{\mO_n}}}} \;
\int\limits_{{\scaleobj{0.9}{\widehat{G}_{\mO_n}}}} }}
(\Psi_{\textsl{q\,}_{\!g} \textsl{q\;\!}'} (f))(\bbz)\;\! 
\text{d}m_{{\scaleobj{0.9}{\widehat{G}_{\mO}}}}(\textsl{q\;\!}')
\\[0.1cm]
\!\!\!&=&\!\!\!
\displaystyle 
{\scaleobj{0.9}{
\sum\limits_{g\;\! \in\;\!  {\scaleobj{0.9}{\widehat G_\mO}}/{\scaleobj{0.9}{\widehat G_{\mO_n}}}}
\frac{1}{| G/G_{\mO_n} |}}} 
(\Psi_{\textsl{q\,}_{\!g}} (f))(\bbz), 
\;\;\quad \bbz\in{\scaleobj{0.9}{\mD^\mO}},
\end{array}
$$
which is exactly the averaging sum \eqref{e8.41}, as required.
\end{proof}

%============================= 
\section{Proof of Theorem~\ref{teo1.11}} 
\label{section_9}
%=============================

We start with the following analog of the classical Wiener lemma \cite[p.91]{Wie}, 
proved in \cite[Thm.~1]{HW} (see also \cite{GN} for an elementary proof).

\begin{lemma}
\label{lem9.1}
The algebra $\mathscr{W}_G$ is {\em inverse closed}$,$ 
i.e.$,$if $f\in \mathscr{W}_G$ is  such that $\inf_{z\;\!\in\;\!\mC_+}|f(z)|>0,$ 
then $1/f\in \mathscr{W}_G$. 
\end{lemma}
\begin{proof}
Since the algebra $\mathscr{W}$ is inverse closed (see, e.g., \cite{GN}), 
there exists a function $g\in\mathscr W$ such that $gf=1$. 
Applying to this equation the projection $\pi_G:\mathscr{W}\to\mathscr W_G$ 
of Theorem~\ref{teo1.7}, we get that $\pi_G(g) f=1$, 
that is, $1/f=g=\pi_G(g)\in \mathscr W_G$, as required.
\end{proof}

$\;\;$Since the uniform closure of $\mathscr W$ is $\mathscr A$, 
applying the projection $\pi_G:\mathscr A\to \mathscr A_G$ 
implies that the uniform closure of $\mathscr W_G$ is $\mathscr A_G$. 
In turn, this and Lemma \ref{lem9.1} imply that the maximal ideal spaces 
$M(\mathscr{W}_G)$ and $M(\mathscr{A}_G)$ of $\mathscr{W}_G$ 
and $\mathscr{A}_G$ coincide, see, e.g., \cite[Prop.~3]{Roy}. 
Thus, by Corollary~\ref{cor1.10}, the maximal ideal spaces of 
$\mathscr W_G$, $\mathscr A_G$,
$W({\scaleobj{0.9}{\overline{\mD}^{\mO}}}/{\scaleobj{0.9}{\widehat G_\mO}})$ 
and $A({\scaleobj{0.9}{\overline{\mD}^{\mO}}}/{\scaleobj{0.9}{\widehat G_\mO}})$ 
are homeomorphic, which proves the first statement of Theorem~\ref{teo1.11}.

To prove the second statement, first, observe that 
$M(A({\scaleobj{0.9}{\mD^\mO}}))={\scaleobj{0.9}{\overline{\mD}^\mO}}$. 
Indeed, let us present the compact set  ${\scaleobj{0.9}{\overline{\mD}^\mO}}$ 
as the inverse limit of the sequence of finite-dimensional closed polydiscs 
${\scaleobj{0.9}{\overline{\mD}^{\mO_i}}}$, where $ \mO_i=\bigcup_{k\le i}O_k$ 
is the union of the first $i$ orbits of the action of $G$ on $\mN$. 
Then the algebra $A({\scaleobj{0.9}{\mD^\mO}})$ 
is the uniform closure of the pullbacks 
of the polydisc algebras $A(\mD^{\mO_i})$ on $\mD^{\mO_i}$ 
under the (inverse limit) projections 
$\pi_i:{\scaleobj{0.9}{\overline{\mD}^{\mO}}} 
\to {\scaleobj{0.9}{\overline{\mD}^{\mO_i}}}$ 
onto the first $|\mO_i|$ coordinates. 
From this, and  using the fact that $M(A(\mD^{\mO_i}))=\overline{\mD}^{\mO_i}$, 
we conclude that $M(A({\scaleobj{0.9}{\mD^\mO}}))={\scaleobj{0.9}{\overline{\mD}^\mO}}$, 
as stated, see, e.g., \cite[Prop.~9]{Roy}.

Next, let us show that functions in 
$A({\scaleobj{0.9}{\mD^\mO}})_{{\scaleobj{0.9}{\widehat G_\mO}}}$ 
separate the orbits of the action of ${\scaleobj{0.9}{\widehat G_\mO}}$ 
on $M(A({\scaleobj{0.9}{\mD^\mO}}))={\scaleobj{0.9}{\overline{\mD}^\mO}}$.\footnote{Here 
and also below, we assume without loss of generality 
that functions of $A({\scaleobj{0.9}{\mD^\mO}})$ 
are extended by continuity to ${\scaleobj{0.9}{\overline{\mD}^\mO}}$.}

Let  $F_1$ and $F_2$ be two distinct orbits 
of the action of ${\scaleobj{0.9}{\widehat G_\mO}}$ 
on \penalty-10000 $M(A({\scaleobj{0.9}{\mD^{\mO}}}))
={\scaleobj{0.9}{\overline{\mD}^{\mO}}}$. 
Since the group $ {\scaleobj{0.9}{\widehat G_\mO}}$ 
is compact and acts continuously on ${\scaleobj{0.9}{\overline{\mD}^{\mO}}}$ 
by Proposition \ref{prop1.8}, 
$F_1$ and $F_2$ are compact disjoint subsets 
of ${\scaleobj{0.9}{\overline{\mD}^{\mO}}}$. 
Hence, since the inverse limit of an inverse family of nonvoid compact spaces 
is a nonvoid compact space, there exists some $i$ 
such that $\pi_i(F_1)$ and $\pi_i(F_2)$ are compact disjoint subsets 
of ${\scaleobj{0.9}{\overline{\mD}^{\mO_i}}}$. 
Since by definition, see Section~\ref{sec1.4a}, 
${\scaleobj{0.9}{\widehat G_\mO}}$ acts via permutations of coordinates 
on ${\scaleobj{0.9}{\overline{\mD}^{\mO}}}$ by permuting, for all $i$, 
the elements of a subset of the first $|O_i|$ coordinates among themselves, 
$\pi_i$ maps each orbit of the action of 
${\scaleobj{0.9}{\widehat G_\mO}}$ on 
${\scaleobj{0.9}{\overline{\mD}^{\mO}}}$ 
to the orbit of the action of ${\scaleobj{0.9}{\widehat G_\mO}}$ 
on ${\scaleobj{0.9}{\overline{\mD}^{\mO_i}}}$ 
via permutations of the first $|O_i|$ coordinates of ${\scaleobj{0.9}{\overline{\mD}^{\mO}}}$. 
The kernel of this action is the group ${\scaleobj{0.9}{\widehat G_{\mO_i}}}$, 
the closure in $Q_\mO$ of the group $\rho(G_{\mO_i})$ (see equation \eqref{eq1.6}), 
where $G_{\mO_i}$ is the normal subgroup of $G$ 
consisting of all $\sigma\in G$ acting identically on  $\mO_i$, see Section~\ref{subsect6.2}. 
Thus the action of ${\scaleobj{0.9}{\widehat G_\mO}}$ 
on ${\scaleobj{0.9}{\overline{\mD}^{\mO_i}}}$ reduces to the effective action 
of the finite quotient group $\widetilde G_i
:={\scaleobj{0.9}{\widehat G_\mO}}/{\scaleobj{0.9}{\widehat G_{\mO_i}}}$ 
via permutation of coordinates on ${\scaleobj{0.9}{\overline{\mD}^{\mO_i}}}$. 
In particular, $\pi_i(F_1)$ and $\pi_i(F_2)$ are finite disjoint subsets 
of ${\scaleobj{0.9}{\overline{\mD}^{\mO_i}}}$, 
and therefore their union is an interpolating set 
for the algebra of holomorphic polynomials on $\mC^{\mO_i}$. 
This and \cite[Thm.~1.1]{Bjo} imply that there exists a function 
$f\in A(\mD^{\mO_i})$ invariant with respect to the action of $\widetilde G_i$, 
such that $f|_{\pi_i(F_1)}=0$ and $f|_{\pi_i(F_2)}=1$. 
Then the pullback $\pi_i^*(f)$ of $f$ is a function of 
$A({\scaleobj{0.9}{\mD^\mO}})_{{\scaleobj{0.9}{\widehat G_\mO}}}$ 
such that $\pi^*_{i}(f)|_{F_1}=0$ and $\pi^*_{i}(f)|_{F_2}=1$. 
Thus, functions in $A({\scaleobj{0.9}{\mD^\mO}})_{{\scaleobj{0.9}{\widehat G_\mO}}}$ 
separate the orbits of the action of ${\scaleobj{0.9}{\widehat G_\mO}}$, 
or equivalently, functions in 
$A({\scaleobj{0.9}{\overline{\mD}^\mO}}/{\scaleobj{0.9}{\widehat G_\mO}})$ 
separate the points of 
${\scaleobj{0.9}{\overline{\mD}^\mO}}/{\scaleobj{0.9}{\widehat G_\mO}}$.

The proof ends by a well-known argument. 
Specifically, if, on the contrary, there exists an $x\in 
M(A({\scaleobj{0.9}{\overline{\mD}^\mO}}/{\scaleobj{0.9}{\widehat G_\mO}})) 
\setminus  ({\scaleobj{0.9}{\overline{\mD}^\mO}}/{\scaleobj{0.9}{\widehat G_\mO}})$, 
then by the definition of the Gelfand topology on 
$M(A({\scaleobj{0.9}{\overline{\mD}^\mO}}/{\scaleobj{0.9}{\widehat G_\mO}}))$, 
there exist elements $f_1,\dots, f_k 
\in A({\scaleobj{0.9}{\overline{\mD}^\mO}}/{\scaleobj{0.9}{\widehat G_\mO}})$ 
such that 
\begin{equation}
\label{eq9.42}
\hat f_1(x)=\cdots =\hat f_k(x)=0
\;\text{ and }\;
\max\limits_{1\le i\le k}|f_i(z)|>0
\;\text{  for all }\; z \in 
{\scaleobj{0.9}{\overline{\mD}^\mO}}/{\scaleobj{0.9}{\widehat G_\mO}}.
\end{equation}
This implies that 
$$
\max\limits_{1\le i\le k}   |(\pi^*(f_i))(z)|>0
\;\text{ for all }\; z\in {\scaleobj{0.9}{\overline{\mD}^\mO}}.
$$
Here $\pi: {\scaleobj{0.9}{\overline{\mD}^\mO}}  \to 
{\scaleobj{0.9}{\overline{\mD}^\mO}}/{\scaleobj{0.9}{\widehat G_\mO}}$ 
is the quotient map.

Since $\pi^*(f_i)\in A({\scaleobj{0.9}{\mD^\mO}})$, $1\le i\le k$, 
and $M(A({\scaleobj{0.9}{\mD^\mO}}))={\scaleobj{0.9}{\overline{\mD}^\mO}}$, 
the previous condition implies that the family $\pi^*(f_1),\dots, \pi^*(f_k)$ 
does not belong to a maximal ideal of $A({\scaleobj{0.9}{\mD^\mO}})$, 
i.e., there exist $\widetilde g_1,\dots, \widetilde g_k\in A({\scaleobj{0.9}{\mD^\mO}})$ 
such that
$$
\widetilde g_1\pi^*(f_1)+\cdots +\widetilde g_k \pi^*(f_k)=1.
$$
Applying the projection $P_{{\scaleobj{0.9}{\widehat G_\mO}}}: 
A({\scaleobj{0.9}{\mD^{\mO}}}) \to 
A({\scaleobj{0.9}{\mD^\mO}})_{{\scaleobj{0.9}{\widehat G_\mO}}}$ 
of Remark \ref{rem1.11} to this equation, we get
\begin{equation}
\label{eq9.43}
P_{{\scaleobj{0.9}{\widehat G_\mO}}}(\widetilde g_1)\pi^*(f_1)
+\cdots 
+P_{{\scaleobj{0.9}{\widehat G_\mO}}}(\widetilde g_k)\pi^*(f_k)
=1.
\end{equation}
Next, by definition, there exist $g_i\!\in\! 
A({\scaleobj{0.9}{\overline{\mD}^\mO\!}}/{\scaleobj{0.9}{\widehat G_\mO}})$ 
such that $\pi^*(  g_i)\!=\!P_{{\scaleobj{0.9}{\widehat G_\mO}}}\!(\widetilde g_i)$, 
$1\le i\le k$. This and \eqref{eq9.43} imply that
\begin{equation}
\label{eq9.44}
g_1f_1+\cdots +  g_k f_k=1.
\end{equation} 
Applying the character $x$ to \eqref{eq9.44}, we get, due to \eqref{eq9.42}, that
$$
1
=\hat 1(x)
=\hat g_1(x)\hat f_1(x)+\cdots +\hat g_k(x)\hat f_k(x)
=0,
$$ 
a contradiction which, together with the fact that the algebra 
$A({\scaleobj{0.9}{\overline{\mD}^\mO}}/{\scaleobj{0.9}{\widehat G_\mO}})$ 
separates the points of 
${\scaleobj{0.9}{\overline{\mD}^\mO}}/{\scaleobj{0.9}{\widehat G_\mO}}$, 
shows that 
$M(A({\scaleobj{0.9}{\overline{\mD}^\mO}}/{\scaleobj{0.9}{\widehat G_\mO}})) 
={\scaleobj{0.9}{\overline{\mD}^\mO}}/{\scaleobj{0.9}{\widehat G_\mO}}$, 
as required.
 
This completes the proof of Theorem~\ref{teo1.11}. 

%=============================
\section{Proof of Theorem~\ref{teo1.12}} 
%=============================

From Corollary~\ref{cor1.10}, it follows that $M((\mathscr O_{b}[\mO])_G)$ 
and $M(\mathcal O_w({\scaleobj{0.9}{\mD^{\mO}}})_{{\scaleobj{0.9}{\widehat G_\mO}}})$ 
are homeomorphic. Moreover, as $\pi^*: 
\calO_w({\scaleobj{0.9}{\mD^\mO}}/{\scaleobj{0.9}{\widehat{G}_\mO}}) 
\to \calO_w({\scaleobj{0.9}{\mD^\mO}})_{{\scaleobj{0.9}{\widehat G_\mO}}}$ 
is an isometric algebra isomorphism, we conclude that 
the maximal ideal spaces of $\mathscr O_{b,G}\, (=(\mathscr O_{b}[\mO])_G)$ 
and $\mathcal O_w({\scaleobj{0.9}{\mD^{\mO}}}/{\scaleobj{0.9}{\widehat G_\mO}})$ 
are homeomorphic. This proves the first statement of Theorem~\ref{teo1.12}. 

Now, let us prove the second statement. 
To this end, we first show that 
$M(\mathcal{O}_w({\scaleobj{0.9}{\mD^\mO}}))
=({\scaleobj{0.9}{\mD^\mO}},\tau_{hk})$ 
(see Subsection~\ref{subsec_1.5} for the definition of the topology 
on ${\scaleobj{0.9}{\mD^\mO}}$). Let $\mathcal A_r$ 
be the uniform closure of the restriction of the algebra 
$\mathcal O_w({\scaleobj{0.9}{\mD^\mO}})$ 
to $r{\scaleobj{0.9}{\overline{\mD}^\mO}}$, $r\in (0,1)$. 

\begin{lemma}
\label{lemma10.1}
The pullback of the dilation map 
$\delta_{\mathfrak{r}},$ $\bbz\mapsto r\bbz,$ $\bbz\in {\scaleobj{0.9}{\mD^\mO}},$ 
maps $\mathcal A_r$ isometrically onto $A({\scaleobj{0.9}{\mD^\mO}})$.
\end{lemma}
\begin{proof} 
If $|\mO|<\infty$, then the statement is obvious. 
Otherwise, the set $\mO$ is countable, 
and any bijection between $\mO$ and $\mN$ 
results in an isometric isomorphism between $c_0(\mO)$ and $c_0$, 
and, in particular, the Fr\'echet spaces $\mathcal O_b(B_{c_0(\mO)})$ 
and  $\mathcal O_b(B_{c_0})$ are isomorphic. 
Thus, all facts formulated for $\mathcal O_b(B_{c_0})$ in Section~\ref{Section_3} 
are also valid for $\mathcal O_b(B_{c_0(\mO)})$. 
In particular, the restriction of each $f\in \mathcal O_b (B_{c_0(\mO)})$ 
to $rB_{c_0(\mO)}$ is uniformly continuous, see \eqref{eq2.13}. 
This implies that the pullback of $\mathcal A_r|_{rB_{c_0(\mO)}}$ 
under the dilation map $\delta_{\mathfrak{r}}$ 
belongs to $\mathcal A_u (B_{c_0(\mO)})$, 
the Banach subalgebra of functions  from $H^\infty(B_{c_0(\mO)})$ 
that are uniformly continuous on $B_{c_0(\mO)}$. 
Since it contains the algebra of polynomials 
in coordinates $z_i$, $i\in\mO$, of $\bbz\in B_{c_0(\mO)}$, 
the uniform closure of the image of the pullback of 
$\mathcal A_r|_{rB_{c_0(\mO)}}$ under the dilation map 
is $ \mathcal A_u (B_{c_0(\mO)})$, as required 
(because  $A({\scaleobj{0.9}{\mD^\mO}})|_{B_{c_0(\mO)}} 
= \mathcal A_u (B_{c_0(\mO)})$, 
and each function of $\mathcal A_u (B_{c_0(\mO)})$ 
extends by weak$^*$ continuity to ${\scaleobj{0.9}{\overline{\mD}^\mO}}$). 
\end{proof}

Next, the dilation map sends ${\scaleobj{0.9}{\overline{\mD}^\mO}}$ 
to $r{\scaleobj{0.9}{\overline{\mD}^\mO}}$, and hence, 
by Lemma~\ref{lemma10.1} and the previous argument 
(in the proof of Theorem~\ref{teo1.11}), 
the maximal ideal space of $\mathcal A_r$ is $r{\scaleobj{0.9}{\overline{\mD}^\mO}}$. 
Using \cite[Thm.~3.2.8, p.80]{Gol}, 
we conclude that the union $\bigcup_{r\in (0,1)} r{\scaleobj{0.9}{\overline{\mD}^\mO}}$ 
considered as the injective limit of the maximal ideal spaces 
of the algebras $\mathcal A_r$, adjoint to the projective limit of these algebras 
(where the projective limit system is defined by restriction homomorphisms), 
is the maximal ideal space of $\mathcal O_w ({\scaleobj{0.9}{\mD^\mO}})$ 
equipped with the corresponding hemicompact $k$-space topology. 

Since $A({\scaleobj{0.9}{\mD^\mO}})_{{\scaleobj{0.9}{\widehat G_\mO}}}$ 
is a subalgebra of 
$\mathcal O_w({\scaleobj{0.9}{\mD^\mO}})_{{\scaleobj{0.9}{\widehat G_\mO}}}$ 
which separates the orbits of the action of 
${\scaleobj{0.9}{\widehat G_\mO}}$ on ${\scaleobj{0.9}{\mD^\mO}}$ 
(see the proof of Theorem~\ref{teo1.11}), 
the algebra 
$\mathcal O_w({\scaleobj{0.9}{\mD^{\mO}}}/{\scaleobj{0.9}{\widehat G_\mO}}) 
\subset C(({\scaleobj{0.9}{\mD^\mO}}/{\scaleobj{0.9}{\widehat G_\mO}},\tau_{hk}))$ 
separates the points of ${\scaleobj{0.9}{\mD^\mO}}/{\scaleobj{0.9}{\widehat G_\mO}}$. 

Using this fact, we obtain that 
$M(\calO_w({\scaleobj{0.9}{\overline{\mD}^\mO}}/{\scaleobj{0.9}{\widehat G_\mO}})) 
= ({\scaleobj{0.9}{\overline{\mD}^\mO}}/{\scaleobj{0.9}{\widehat G_\mO}},\tau_{hk})$ 
by the same standard argument given at the end of Section~\ref{section_9}, 
by applying the projection $P_{{\scaleobj{0.9}{\widehat G_\mO}}}: 
\mathcal O_w({\scaleobj{0.9}{\mD^{\mO}}}) 
\to \mathcal O_w({\scaleobj{0.9}{\mD^\mO}})_{{\scaleobj{0.9}{\widehat G_\mO}}}$ 
of Remark \ref{rem1.11} to the corresponding B\'ezout equation, 
using that $M(\calO_w({\scaleobj{0.9}{\overline{\mD}^\mO}})) 
= ({\scaleobj{0.9}{\mD^\mO}}, \tau_{hk})$ (cf. \eqref{eq9.43}). 
We leave the details to the reader.
  
%=====================================
\section{Proofs of Results of Section~\ref{Sec2}}
%=====================================

%---------------------------------------------------------------------
\subsection{Proof of Theorem~\ref{teo_contractibility}}
%---------------------------------------------------------------------

In light of Theorems~\ref{teo1.11} and \ref{teo1.12}, 
we need to show that the topological spaces 
${\scaleobj{0.9}{\overline{\mD}^\mO}}/{\scaleobj{0.9}{\widehat{G}_\mO}}
\,(:=({\scaleobj{0.9}{\overline{\mD}^\mO}},\tau_p)/{\scaleobj{0.9}{\widehat{G}_\mO}}$), 
and \penalty-10000 $({\scaleobj{0.9}{\mD^\mO}}/{\scaleobj{0.9}{\widehat{G}_\mO}}, \tau_{hk})$ 
are contractible. 

In what follows, $I\!:=\![0,1]\!\subset\!\mR$ is equipped with the standard topology. 
First, we observe that the space ${\scaleobj{0.9}{\overline{\mD}^\mO}}$ 
equipped with the product topology $\tau_p$ is contractible.  
Indeed, the map $F: I\times ({\scaleobj{0.9}{\overline{\mD}^\mO}}, \tau_p) 
\to ({\scaleobj{0.9}{\overline{\mD}^\mO}},\tau_p)$, 
$F(t,\bbz):=t\cdot \bbz$ for all $t\in I$, 
$\bbz\in {\scaleobj{0.9}{\overline{\mD}^\mO}}$, is continuous, 
and determines a homotopy between 
$F(1,\cdot)={\rm id}_{{\scaleobj{0.9}{\overline{\mD}^\mO}}}$, 
and the constant map $F(0,\cdot)=\mathbf 0$. 
Next, since ${\scaleobj{0.9}{ \widehat{G}_\mO}}$ 
acts by linear transformations on $\ell^\infty(\mO)$, 
for all $g\in {\scaleobj{0.9}{\widehat{G}_\mO}}$ 
and $(t,\bbz)\in I\times {\scaleobj{0.9}{\overline{\mD}^\mO}}$, we have
\begin{equation}
\label{12.45}
\;\;F(t,g(\bbz))=g(t\cdot \bbz)=g(F(t,\bbz)). 
\end{equation}
This implies that there exists a map 
$H:I\times ({\scaleobj{0.9}{\overline{\mD}^\mO}}/{\scaleobj{0.9}{\widehat{G}_\mO}}) 
\to {\scaleobj{0.9}{\overline{\mD}^\mO}}/{\scaleobj{0.9}{\widehat{G}_\mO}}$, 
such that the diagram
\begin{center}
\begin{tikzcd}
I\!\times\! {\scaleobj{0.9}{\overline{\mD}^\mO}}  
\arrow{r}[name=U]{F} 
& {\scaleobj{0.9}{\overline{\mD}^\mO}}  \\
I \!\times\! ({\scaleobj{0.9}{\overline{\mD}^\mO\!}}/{\scaleobj{0.9}{\widehat{G}_\mO}}) 
\arrow[leftarrow]{u}{\text{id}_I\times \pi } 
\arrow{r}[name=U]{H} 
& {\scaleobj{0.9}{\overline{\mD}^\mO\!}}/{\scaleobj{0.9}{\widehat{G}_\mO}} 
\arrow[leftarrow,swap]{u}{\pi}
\end{tikzcd}
\end{center}
is commutative. 

By the definition, 
$H(1,\cdot)={\rm id}_{{\scaleobj{0.9}{\overline{\mD}^\mO\!/\widehat{G}_\mO}}}$ 
and $H(0,\cdot)=\pi({\mathbf 0})$. Thus, to prove that the space 
${\scaleobj{0.9}{\overline{\mD}^\mO\!}}/{\scaleobj{0.9}{\widehat{G}_\mO}}$ 
is contractible, it suffices to show that the map $H$ is continuous.

To this end, let $S$ be a closed subset of 
${\scaleobj{0.9}{\overline{\mD}^\mO\!}}/{\scaleobj{0.9}{\widehat{G}_\mO}}$. 
Then the above diagram implies that 
$(\text{id}_I\times \pi)^{-1}(H^{-1}(S)) = F^{-1}(\pi^{-1}(S)) =: \widetilde{S}$. 
Since $F$ and $\pi$ are  continuous maps, 
$\widetilde{S}$ is a closed subset of $I\times {\scaleobj{0.9}{\overline{\mD}^\mO}}$, 
and so it is compact, as $I\times{\scaleobj{0.9}{ \overline{\mD}^\mO}}$ is a compact space. 
Thus, since $\text{id}_I\times \pi$ is a continuous map, 
$H^{-1}(S)= (\text{id}_I\times \pi)(\widetilde{S})$ is a compact subset of 
${\scaleobj{0.9}{\overline{\mD}^\mO\!}}/ {\scaleobj{0.9}{\widehat{G}_\mO}}$, 
and in particular, is closed. This proves that $H$ is a continuous contraction of 
${\scaleobj{0.9}{\overline{\mD}^\mO\!}}/{\scaleobj{0.9}{\widehat{G}_\mO}}$ 
to $\pi(\mathbf{0})$.
   
Next, to show the contractibility of 
$({\scaleobj{0.9}{\mD^\mO}}/{\scaleobj{0.9}{\widehat G_\mO}}, \tau_{hk})$, 
we use the same maps $F$ and $H$ as above, 
restricted to $I\times{\scaleobj{0.9}{\mD^\mO}}$ and 
$I\times {\scaleobj{0.9}{\mD^\mO}}/{\scaleobj{0.9}{\widehat G_\mO}}$, respectively. 
Then the images of the restricted maps are ${\scaleobj{0.9}{\mD^\mO}}$ 
and ${\scaleobj{0.9}{\mD^\mO}}/{\scaleobj{0.9}{\widehat G_\mO}}$, 
respectively, as required. Assuming that ${\scaleobj{0.9}{\mD^\mO}}$ 
is given the hemicompact $k$-topology $\tau_{hk}$, 
we must prove that $H: 
I\times  ({\scaleobj{0.9}{\mD^\mO}}/{\scaleobj{0.9}{\widehat G_\mO}}, \tau_{hk}) 
\to ( {\scaleobj{0.9}{\mD^\mO}}/{\scaleobj{0.9}{\widehat G_\mO}}, \tau_{hk})$ 
is continuous. To this end, since the space 
$( {\scaleobj{0.9}{\mD^\mO}}/{\scaleobj{0.9}{\widehat G_\mO}}, \tau_{hk})$ 
is Hausdorff and the space 
$I\times  ({\scaleobj{0.9}{\mD^\mO}}/{\scaleobj{0.9}{\widehat G_\mO}}, \tau_{hk})$ 
is hemicompact, it is sufficient to prove that $H$ is continuous 
on each compact set $I\times \pi(r{\scaleobj{0.9}{\overline{\mD}^{\mO}}})$, $r\in (0,1)$, 
see, e.g., \cite[Prop.~2.7]{Ste}. 
In this case, the image of $H$ is $\pi(r{\scaleobj{0.9}{\overline{\mD}^{\mO}}})$, 
and the map $F$ sends $I\times r{\scaleobj{0.9}{\overline{\mD}^{\mO}}}$ 
to $r{\scaleobj{0.9}{\overline{\mD}^{\mO}}}$, 
so that the continuity of $H$ follows as in the above considered case, 
with $r{\scaleobj{0.9}{\overline{\mD}^{\mO}}}$ 
in place of ${\scaleobj{0.9}{\overline{\mD}^{\mO}}}$.

%----------------------------------------------------
\subsection{Proof of Theorem \ref{te2.2}}
%----------------------------------------------------

Let $A\in \{\mathscr W_G, \mathscr A_G, \mathscr H^\infty_G\}$. 
If $f\in A^{-1}$, then $0$ does not belong to the compact set $\hat f(M(A))\subset\mC$. 
Hence there exists a $\delta\!>\!0$ such that $\min_{M(A)}\! |\hat f|\!\ge\!\delta$. 
Since $\mC_+\!\subset \!M(A)$, this implies $\inf_{ \mC_+} |f|>0$. 
Conversely, if $f\in \mathscr H_G^\infty$ is such that $\inf_{\mC_+} |f|>0$, 
then $1/f$ is a bounded holomorphic  function on $\mC_+$, 
which by  \cite[Thm.~2.6]{Bon}, belongs to $\mathscr O_u$. 
Thus, $1/f\in\mathscr H^\infty$. In fact, since $f$ is $G$-invariant, 
as in the proof of Lemma~\ref{lem9.1}, we get $1/f\in\mathscr H_G^\infty$. 
Moreover, if $f\in\mathscr W_G$, then $1/f\in \mathscr W_G$ (see Lemma~\ref{lem9.1}), 
and if $f\in \mathscr A_G$, then also $1/f\in \mathscr A_G$, 
because $\mathscr A_G$ is the uniform closure of $\mathscr W_G$. 

So we have proved $A^{-1}\!=\!\{f\in A : \inf\limits_{\mC_+}|f|>0\}$ 
if $A\!\in\! \{\mathscr W_G, \mathscr A_G,\mathscr H_G^\infty\}$. 

Next, suppose that $f\in\mathscr O_{b,G}$ satisfies 
$f_r\in (\mathscr H_G^\infty)^{-1}$  for all $r\in (0,1)$. 
Let $F\in (H^\infty_w({\scaleobj{0.9}{\mD^\mO}}))_{{\scaleobj{0.9}{\widehat G_\mO}}}$ 
be such that $c_\mO^*(F)=f$, see Corollary \ref{cor1.10}.  
Then, $f_r=c_{\mO}^*(F_{\texttt{r}})$ for all $r\in (0,1)$, cf. \eqref{eq3.22}. 
This and the assumption $f_r\in (\mathscr H_G^\infty)^{-1}$ imply that 
$1/F_{\texttt{r}} 
\in (H^\infty_w({\scaleobj{0.9}{\mD^\mO}}))_{{\scaleobj{0.9}{\widehat G_\mO}}}$, 
i.e., $|1/F|$ is  bounded from above on each $r{\scaleobj{0.9}{\mD^\mO}}$, $r\in (0,1)$. 
Hence, $1/F\in  (H^\infty_w({\scaleobj{0.9}{\mD^\mO}}))_{{\scaleobj{0.9}{\widehat G_\mO}}}$. 
Thus $c_{\mO}^*(1/F)\in \mathscr O_{b,G}$ and $f\cdot c_{\mO}^*(1/F)=1$. 
This shows that $f\in (\mathscr O_{b,G})^{-1}$.

Conversely, if $f\in (\mathscr O_{b,G})^{-1}$, 
then $f\cdot(1/f)=1$ and $1/f\in \mathscr O_{b,G}$. 
Thus, $f_r\cdot (1/f)_r=1$ for all $r\in (0,1)$. 
Since $(1/f)_r\in\mathscr H^\infty$ by the definition of $\mathscr O_b$, 
this shows that each $f_r\in (\mathscr H_G^\infty)^{-1}$, 
and completes the proof of the first part of the theorem.

Now, let us prove the second statement of the theorem. 
If $A$ is  a unital commutative complex Banach algebra, 
then according to the Arens-Royden theorem (see, e.g., \cite[Thm.,\,p.295]{Roy}), 
the group $A^{-1}/e^{A}$ is isomorphic 
to the first \v{C}ech cohomology group $H^1(M(A), \mZ)$ of $M(A)$ with integer coefficients.  
(For background on \v{C}ech cohomology, see, e.g., \cite{EilSte}.) 
An analogous result due to Brooks also holds 
for unital commutative complex Fr\'echet algebras 
whose maximal ideal spaces are $k$-spaces, 
see, e.g., \cite[Sect.\,6.2,\,Cor.]{Gol}. 
In addition, for a contractible space, 
all cohomology groups are trivial (see, e.g., \cite[IX,\,Thm.\,3.4]{EilSte}). 
Applying these to the algebras 
$A\in \{\mathscr W_G, \mathscr A_G, \mathscr O_{b,G}\}$, 
we get $A^{-1}=e^{A}$. 

%-------------------------------------------------------
\subsection{Proof of Theorem~\ref{thm2.4}} 
\label{subsekt2.3}
%-------------------------------------------------------

If $|\mO|<\infty$, then ${\scaleobj{0.9}{\widehat G_{\mO}}}=G$ is a finite group, 
and for $\mathcal A\in
\{\mathcal O_w({\scaleobj{0.9}{\mD^\mO}}), 
H_w^\infty({\scaleobj{0.9}{\mD^\mO}}), 
A({\scaleobj{0.9}{\mD^\mO}}), 
W({\scaleobj{0.9}{\mD^\mO}})\}$, 
the projection $P_{G}: \mathcal A\to \mathcal A_{G}$ of Proposition \ref{prop1.11} 
has the form
$$
\textstyle
(P_{ G}(f))(\bbz)
= {\scaleobj{1.2}{\frac{1}{|G|}}} \sum\limits_{\sigma\;\! \in \;\! G}f(\sigma(\bbz)), 
\  \bbz\in{\scaleobj{0.9}{\mD^\mO}}, 
\  f\in\mathcal A,
$$
where if $\mO=\{i_1,\dots, i_{|\mO|}\}\subset\mN$, 
then $\sigma(\bbz):= (z_{\sigma(i_1)},\dots, z_{\sigma(i_{|\mO|})})$ 
for $\bbz=(z_{i_1},\dots, z_{i_{|\mO|}})\in\mC^\mO$.
Clearly, $P_{ G}$ maps the algebra $\mC[z_{i_1},\dots,z_{i_{|\mO|}}]$ 
of holomorphic polynomials on $\mC^\mO$ 
onto a subalgebra $\mathcal P_G$ of polynomials $p$ 
invariant with respect to the action of $G$, 
i.e., such that $p(\sigma(\bbz))=p(\bbz)$, $\bbz\in\mC^\mO$, for all $\sigma\in G$. 
By a result due to Chevalley (see \cite[Thm.~A]{Che}), 
there exist $|\mO|$ homogeneous  polynomials 
$h_1,\dots, h_{|\mO|}\in  \mC[z_{i_1},\dots,z_{i_{|\mO|}}]$, 
algebraically independent over $\mC$, 
that are invariant with respect  to the action of $G$, 
such that $\mathcal P_G$ is the free polynomial algebra 
over $\mC$ in $h_1,\dots, h_{|\mO|}$.

Further, the uniform closure of the algebra 
$(\mC[z_{i_1},\dots,z_{i_{|\mO|}}])|_{{\scaleobj{0.9}{\mD^\mO}}}$ 
is the polydisc algebra $A({\scaleobj{0.9}{\mD^\mO}})$. 
Then, applying the projection $P_{ G}$, 
we obtain that the uniform closure of the algebra 
$(\mathcal P_G)|_{{\scaleobj{0.9}{\mD^\mO}}}$ 
is the algebra $A({\scaleobj{0.9}{\mD^\mO}})_{G}$. 
Hence, $A({\scaleobj{0.9}{\mD^\mO}})_{G}$ is generated by functions 
$h_1|_{{\scaleobj{0.9}{\mD^\mO}}},\dots,h_{|\mO|}|_{{\scaleobj{0.9}{\mD^\mO}}}$. 
In turn, there exist functions 
$g_1,\dots, g_{|\mO|}\in A({\scaleobj{0.9}{\overline{\mD}^\mO}}/G)$ 
such that $\pi^*(g_i)=h_i$, $1\le i\le |\mO|$.\footnote{Recall that 
$\pi: {\scaleobj{0.9}{\overline{\mD}^\mO}} 
\to {\scaleobj{0.9}{\overline{\mD}^\mO}}/G$ 
is the quotient map, see Section~\ref{sec1.4a}.} 
Thus, the algebra $A({\scaleobj{0.9}{\overline{\mD}^\mO}}/G)$ 
is generated by $g_1,\dots, g_{|\mO|}$ 
and, in particular, these functions 
separate the points of ${\scaleobj{0.9}{\overline{\mD}^\mO}}/G$, 
the maximal ideal space of $A({\scaleobj{0.9}{\overline{\mD}^\mO}}/G)$.
 
Let us consider the holomorphic map 
$H\!=\!(h_1,\dots, h_{|\mO|})\!:\!\mC^\mO\!\to\!\mC^{|\mO|}$. 
Since $H$ is proper, by Remmert's proper mapping theorem 
(see, e.g., \cite[Ch.~3,\,Sect.~2,\,p.395]{GriHar}), 
$H$ is surjective. Next, let $\partial{\scaleobj{0.9}{\mD^\mO}}$ 
be the boundary of ${\scaleobj{0.9}{\mD^\mO}}$, 
and $\Gamma:=H(\partial{\scaleobj{0.9}{\mD^\mO}})$. 
(It is worth noting that since $\Gamma\subset\mC^{|\mO|}$ 
is the image of a finite analytic map restricted 
to the compact subanalytic set 
$\partial{\scaleobj{0.9}{\mD^\mO}}\subset\mC^{\mO}$, 
it is a $(2|\mO|-1)$-dimensional subanalytic set.)  
The holomorphic map of complex manifolds 
$H|_{\mC^\mO\setminus \partial{\scaleobj{0.9}{\mD^\mO}}}: 
\mC^\mO\setminus \partial{\scaleobj{0.9}{\mD^\mO}}
\to \mC^{|\mO|}\setminus\Gamma$ is also proper, 
so that Remmert's proper mapping theorem implies that 
the map $H|_{\mC^\mO\setminus \partial{\scaleobj{0.9}{\mD^\mO}}}$ is surjective, 
and the open set $\mC^{|\mO|}\setminus\Gamma$ 
has two connected components: 
the bounded component $H({\scaleobj{0.9}{\mD^\mO}})$, 
and the unbounded component 
$H(\mC^{\mO}\setminus{\scaleobj{0.9}{\overline{\mD}^{\mO}}})$. 
Moreover, $\Gamma$ is the boundary of both components. 
 
Finally, the joint spectrum of the generators 
$g_1,\dots, g_{|\mO|}$ of the algebra 
$A({\scaleobj{0.9}{\overline{\mD}^\mO}}/G)$, 
i.e., the set $\{(g_1(x),\dots, g_{|\mO|}(x)) 
\in \mC^{|\mO|}\,:\, x\in {\scaleobj{0.9}{\overline{\mD}^\mO}}/G\}$, 
is a compact subset of $\mC^{|\mO|}$, 
homeomorphic to ${\scaleobj{0.9}{\overline{\mD}^\mO}}/G$, 
and coincident with $H({\scaleobj{0.9}{\mD^\mO}})\cup\Gamma$. 
As the latter has a nonempty interior, 
by \cite[Cor.~3.13]{CorSua} (see also Thms~3.4 and 3.12 there), it follows that 
$\text{sr} (A({\scaleobj{0.9}{\mD^\mO}}/{\scaleobj{0.9}{\widehat{G}_\mO}}))
\!=\!\lfloor \frac{|\mO|}{2}\rfloor\!+\!1$. 
By Theorem~\ref{teo1.11}, we also get the same stable rank 
for $\mathscr W_G$ and $\mathscr A_G$.  
 
In the remaining case of the algebra $(\mathcal O_w({\scaleobj{0.9}{\mD^\mO}}))_G$, 
its maximal ideal space is ${\scaleobj{0.9}{\mD^\mO}}/G_{\mO}$. 
According to our construction, it can be naturally identified 
with the open subset $H({\scaleobj{0.9}{\mD^\mO}})\subset\mC^{|\mO|}$, 
so that by this identification, $(\mathcal O_w({\scaleobj{0.9}{\mD^\mO}}))_G$ 
becomes the algebra  $\mathcal O(H({\scaleobj{0.9}{\mD^\mO}}))$ 
of holomorphic functions on $H({\scaleobj{0.9}{\mD^\mO}})$. 
According to the definition of ${\scaleobj{0.9}{\mD^\mO}}/G_{\mO}$ 
as the maximal ideal space of the Fr\'{e}chet algebra 
$(\mathcal O_w({\scaleobj{0.9}{\mD^\mO}}))_G$, 
see the proof of Theorem~\ref{teo1.12}, 
the set $H({\scaleobj{0.9}{\mD^\mO}})$ admits 
an exhaustion by compact holomorphically convex subsets 
$H(r{\scaleobj{0.9}{\overline{\mD}^\mO}})$, $r\in (0,1)$. 
This implies that the complex manifold $H({\scaleobj{0.9}{\mD^\mO}})$ 
is holomorphically convex and, in particular, 
it is Stein (see, e.g., \cite[Ch.~4]{GraRem}, 
which also contains background on Stein spaces). 
Thus, $\text{sr}(\mathcal O({\scaleobj{0.9}{\mD^\mO}}/{\scaleobj{0.9}{\widehat G_{\mO}}}))
=\text{sr}(\mathcal O(H({\scaleobj{0.9}{\mD^\mO}}))
= \lfloor \frac{|\mO|}{2}\rfloor+1$ 
by \cite[Thm.~1.1]{Bru19}. 
In light of Theorem~\ref{teo1.12}, 
this is also the stable rank of $\mathscr O_{b,G}$.

Next, we give the proof of Theorem~\ref{thm2.4} in the case when $|\mO|=\infty$. 
For the convenience of the reader, 
we recall some relevant notations and definitions 
from Sections~\ref{sec7.1} and \ref{subsect6.2}.  

For a subset $\mathcal P\subset\mP$, 
we denote by $\langle \mathcal P\rangle$ 
the unital multiplicative semigroup contained in $\mN$ 
generated by the elements of $\mathcal P$.  
For a Dirichlet series $D=\sum_{n=1}^\infty a_n n^{-s}$ 
and a subset $\mathcal P\subset\mP$, 
we let $D_\mathcal P$ be the series obtained from $D$ 
by removing all terms with $n\not\in \langle\mathcal P\rangle$, 
that is, $D_\mathcal P:=\sum_{n\;\! \in\;\! \langle\mathcal P\rangle} a_n n^{-s}$. 
The correspondence $\Phi_\mathcal P:\mathscr D\to\mathscr D$, 
$\Phi_\mathcal P(D):=D_\mathcal P$, 
is an algebra epimorphism and a linear projection 
onto a subalgebra $\mathscr D_\mathcal P$ 
of Dirichlet series of the form $\sum_{n\;\! \in\;\! \langle\mathcal P\rangle} a_n n^{-s}$. 
Let $\mA \in\{ \mathscr W , \mathscr A , \mathscr H^\infty,\mathscr O_{b}\} $, 
and $\mathscr D_{\mA}$ be the set of Dirichlet series converging to functions in $\mA$. 
Then we have

\begin{proposition}
\label{prop12.1}
$\Phi_{\mathcal P}$ maps $\mathscr D_{\mA}$ to itself$,$ 
and thus determines a homomorphism $\mA\to\mA,$ 
denoted by the same symbol $\Phi_{\mathcal P}$. 
Moreover$,$ $\Phi_{\mathcal P}$ is a continuous projection of norm one of $\mA$ 
onto a closed subalgebra $\mA_{\mathcal P}$ 
of limit functions of Dirichlet series of $\mathscr D_{\mathcal P}\cap\mathscr D_{\mA}$.
\end{proposition}
\begin{proof}
The statement is true for the algebra $\mathscr W\!$, 
since obviously $\Phi_{\mathcal P}$ does not increase the norm of $\mathscr W$. 
Next, for the algebra $\mathscr H^\infty$, 
the required statement is the content of Corollary~\ref{cor6.2}. 
In turn, the previous two cases imply that 
the statement of the proposition holds 
for the algebra $\mathscr A$, 
since  $\mathscr W$ is dense in $\mathscr A$ in the uniform norm. 
Finally, noting that the projection $\Phi_{\mathcal P}$ 
commutes with the map $\cdot_r:\mathscr D\to\mathscr D $, see Section~\ref{subsect1.2}, 
by the definition of $\mathscr O_b$ and Corollary~\ref{cor6.2}, 
it sends $\mathscr O_b$ to itself, 
and does not increase the seminorms $P_r$, $r\in (0,1)$, 
determining its topology, see \eqref{equation_1.5}. 
Thus we obtain the result also for the algebra $\mathscr O_b$.
\end{proof}

As before, the (finite) orbits of the action of $G$ on $\mO\subset\mN$ 
are labelled as $O_1, O_2, O_3,\dots$  
(e.g., such that if $i<j$, then $\min O_i <\min O_j$). 
Also, for $n\in \mN$, $\mO_n:={\scaleobj{0.9}{\bigcup_{\;\!i=1}^{\;\!n}}\;\!} O_i$, 
and $\mathcal P_{\mO_n}:=\{p_i\in\mP: i\in\mO_n\}$. 

For $\mA$ as above, let $\mA_G$ 
and $\mA_{G,\mathcal P_{\mO_n}}:=\mA_G\cap \mA_{\mathcal P_{\mO_n}}$ 
denote the algebras of limit functions of invariant Dirichlet series 
with respect to the action $S$ of $G$ on $\mA$ and $\mA_{\mathcal P_{\mO_n}}$, 
respectively, see Section~\ref{subsection_3.1}. 
(Note that since $\mO_n$ is invariant with respect to the action of $G$, 
the algebra $\mA_{\mathcal P_{\mO_n}}$ is invariant 
with respect to the action $S$ of $G$ on $\mA$.)

Now, from Proposition~\ref{prop12.1}, we obtain:

\begin{corollary}\label{cor12.2}
The restriction of the homomorphism 
$\Phi_{\mathcal P_{\mO_n}}\!\!:\!\mA\!\to\! \mA_{\mathcal P_{\mO_n}}$ to $\mA_G$ 
determines an algebra epimorphism 
and a continuous projection of norm one onto $\mA_{G,\mathcal P_{\mO_n}}$.
\end{corollary}
\begin{proof} 
The result follows straightforwardly from Proposition \ref{prop12.1} 
by the fact that the action $S$ of $G$ on $\mA$ 
commutes with the homomorphism $\Phi_{\mathcal P_{\mO_n}}$. 
\end{proof}

Now, we are ready to prove that if $|\mO|=\infty$, then ${\rm sr}\,\mA_G=\infty$.

Assume on the contrary that ${\rm sr}\,\mA_G\!=\!r\!<\!\infty$. 
Since $|\mO|\!=\!\infty$, the set of orbits $\{O_i\}$ of the action of $G$ on $\mN$ is infinite, 
and thus there is some $\mO_n$ such that $|\mO_n|\!>\!2r$. 
Thus each $r\!+\!1$ tuple 
$f_1,\dots,f_{r+1}\!\in\! U_{r+1}(\mA_{G,\mathcal P_{\mO_n}})$ 
can be reduced over $\mA_{G}$ to an $r$ tuple in $U_{r}(\mA_{G})$, 
and so there exist $h_1,\dots , h_r, g_1,\dots, g_r \in \mA_G$ such that 
$$
(f_1\!+\!h_1 f_{r+1})g_1+\cdots+ (f_r \!+\!h_r f_{r+1})g_r =1.
$$
Applying the homomorphism $\Phi_{\mathcal P_{\mO_n}}$, we get 
$$
(f_1\!+\!\Phi_{\mathcal P_{\mO_n}}\!(h_1) f_{r+1})\Phi_{\mathcal P_{\mO_n}}\!(g_1)
+ \cdots + 
(f_r \!+\!\Phi_{\mathcal P_{\mO_n}}\!(h_r) f_{r+1})\Phi_{\mathcal P_{\mO_n}}\!(g_r) 
=1.
$$
As $\Phi_{\mathcal P_{\mO_n}}(g_1),\dots, \Phi_{\mathcal P_{\mO_n}}(g_r)
\in \mA_{G,\mathcal P_{\mO_n}}$ by Corollary \ref{cor12.2}, 
we conclude from the above that the $r+1$ tuple 
$f_1,\dots,f_{r+1}$  is reducible over $\mA_{G,\mathcal P_{\mO_n}}$, 
which implies that
\begin{equation}
\label{eq12.46}
{\rm sr}\,\mA_{G,\mathcal P_{\mO_n}}\le r.
\end{equation} 

However, $\mA_{G,\mathcal P_{\mO_n}}$ 
is the algebra of limit functions of invariant Dirichlet series 
with respect to the action $\tilde S$ of the finite group $G/G_{\mO_n}$ 
on $\mA_{\mathcal P_{\mO_n}}$, 
where $G_{\mO_n}$ denotes the normal subgroup of $G$ 
consisting of all $\sigma\in G$ acting identically on $\mO_n$, see Section~\ref{subsect6.2}. 
Then we can apply to $\mA_{G,\mathcal P_{\mO_n}}$ 
the previous case (of the finite set of orbits) 
with $G/G_{\mO_n}$ in place of $G$, 
and $\mO_n$ in place of $\mO$. 
Consequently, we obtain for $\mA\in \{\mathscr W , \mathscr A ,\mathscr O_{b}\}$ 
that $\text{sr}\,\mA_{G,\mathcal P_{\mO_n}} = \lfloor\frac{|\mO_n|}{2}\rfloor +1 > r+1$,  
contradicting \eqref{eq12.46}. 
This shows that if $|\mO|=\infty$, then $\text{sr}\,\mA_{G}=\infty$ for these algebras. 

Further, for $\mA=\mathscr H^\infty$, 
the Banach algebra $\mA_{G,\mathcal P_{\mO_n}}$ 
is isometrically isomorphic to the Banach algebra 
$H^\infty(\mD^{\mO_n}/(G/G_{\mO_n}))$ 
of bounded holomorphic functions 
on the Stein manifold $\mD^{\mO_n}/(G/G_{\mO_n})$. 
As in the proof of the first part of the theorem, 
we identify $\mD^{\mO_n}/(G/G_{\mO_n})$ 
with a bounded Stein domain in $\mD^{|\mO_n|}$ 
containing ${\bf 0}\in\mC^{|\mO_n|}$. 
Expanding the latter if necessary, 
we may assume that 
the open unit Euclidean ball $\mB$ of $\mC^{|\mO|}$ centered at ${\bf 0}$ 
lies in $\mD^{\mO_n}/(G/G_{\mO_n})$. For $s:=\lfloor\frac{|\mO_n|}{2}\rfloor +1$,  
consider the element $u=(f_1,\dots, f_{s-1},p) 
\in U_{s}(H^\infty(\mD^{\mO_n}/(G/G_{\mO_n})))$, where 
$$
\begin{array}{ll}
f_1(\bbz) \!=\! z_1, 
\; f_2(\bbz)\!=\!z_3, 
\; \dots, 
\; f_{s-1}(\bbz)\!=\! z_{2s-3},\\
p(\bbz) = z_1z_2 + z_3z_4 + \cdots + z_{2s-3}z_{2s-2} -1; 
\end{array}
\big(\bbz\!=\! (z_1,\dots,z_{|\mO|}) \!\in\! \mC^{|\mO|}\big).
$$
Assuming that $\text{sr}\,H^\infty(\mD^{\mO_n}/(G/G_{\mO_n}))<s$, 
we obtain that the $s$-tuple $u$ is reducible over 
$H^\infty(\mD^{\mO_n}/(G/G_{\mO_n}))$, and therefore its restriction to $\mB$ 
is reducible over the algebra $A(\mB)$ of holomorphic functions on $\mB$
which extend continuously to its boundary. 
However, it was shown in the proof of \cite[Thm.~3.12]{CorSua} 
that $u|_{\mB}\in U_s(A(\mB))$ is irreducible. 
This contradiction shows that 
${\rm sr}\,\mathscr H^\infty_{G,\mathcal P_{\mO_n}}
\!\ge\!\lfloor\frac{|\mO_n|}{2}\rfloor \!+\!1
\!>\!r\!+\!1$, which in turn contradicts \eqref{eq12.46} 
and shows that if $|\mO|\!=\!\infty$, then ${\rm sr}\,\mathscr H^\infty_{G}\!=\!\infty$ as well.

\begin{remark}
\label{rem13.3}
The transposes of the embedding homomorphism \penalty-10000 
$(\mathscr O_{b})_{G,\mathcal P_{\mO_n}}\hookrightarrow \mathscr O_{b,G}$ 
and the homomorphism 
$\Phi_{\mathcal P_{\mO_n}} : 
\mathscr O_{b,G}\to (\mathscr O_{b})_{G,\mathcal P_{\mO_n}}$ 
determine the following continuous maps of the maximal ideal spaces 
$$
\;\textstyle
\mathscr I_n: M(\mathscr O_{b,G}) \to M((\mathscr O_{b})_{G,\mathcal P_{\mO_n}})
\;\text{ and }\;
\mathscr P_n: M((\mathscr O_{b})_{G,\mathcal P_{\mO_n}}) \to M(\mathscr O_{b,G}), 
$$
respectively, such that 
$\mathscr I_n\circ\mathscr P_n 
= {\rm id}|_{M((\mathscr O_{b})_{G,\mathcal P_{\mO_n}})}$ 
(because $\Phi_{\mathcal P_{\mO_n}}$ 
is a projection onto $(\mathscr O_{b})_{G,\mathcal P_{\mO_n}}$). 
Similarly, if $\mO_i\subset\mO_j$, 
then the transpose of the embedding homomorphism 
$(\mathscr O_{b})_{G,\mathcal P_{\mO_i}} 
\hookrightarrow  (\mathscr O_{b})_{G,\mathcal P_{\mO_j}}$ 
determines a continuous map 
$\mathscr I_{ij}: M((\mathscr O_b)_{G,\mathcal P_{\mO_i}}) 
\to M((\mathscr O_{b})_{G,\mathcal P_{\mO_j}})$. 
It is easy to check that  $M(\mathscr O_{b,G})$ 
is the inverse limit of the inverse system
$\bigl((M((\mathscr O_b)_{G,\mathcal P_{\mO_i}}))_{i\in\mN}, 
(\mathscr I_{ij})_{i\le j\in\mN}\bigr)$, 
and that the set 
${\scaleobj{0.9}{\bigcup_{\;\!n=1}^{\;\!\infty}}} 
\mathscr P_n(M((\mathscr O_{b})_{G,\mathcal P_{\mO_n}}))$ 
is dense in $M(\mathscr O_{b,G})$. As before, each 
$M((\mathscr O_{b})_{G,\mathcal P_{\mO_n}}) 
= \mD^{\mO_n}/(G/G_{\mO_n})$ 
can be identified with a bounded Stein domain $\Omega_n$ in $\mD^{|\mO_n|}$ 
containing ${\bf 0}\in\mC^{|\mO_n|}$, 
such that the Fr\'{e}chet algebra $(\mathscr O_{b})_{G,\mathcal P_{\mO_n}}$ 
is identified with the Fr\'{e}chet algebra $\mathscr O(\Omega_n)$ 
of holomorphic functions on $\Omega_n$.
\end{remark}

%-------------------------------------------------------
\subsection{Proof of Theorem~\ref{teo2.5a}} 
%-------------------------------------------------------

By Theorem~\ref{teo_contractibility},  $M(\mathscr W_G)$ and $M(\mathscr A_G)$ are contractible, and so the projective freeness of $\mathscr W_G$ and $\mathscr A_G$ follows from  \cite[Thm.~4.1]{BruSas23}. 
It remains to show that $\mathscr O_{b,G}$ is projective free. 
By Corollary \ref{cor1.10} (see also Section~\ref{subsec1.7}), it suffices to show that the algebra $\mathcal O_w({\scaleobj{0.9}{\mD^\mO}}/{\scaleobj{0.9}{\widehat G_\mO}})$ is projective free. 

To this end, recall that the set ${\scaleobj{0.9}{\mD^\mO}}/{\scaleobj{0.9}{\widehat G_\mO}}$ 
is the union of compact sets $\pi(r{\scaleobj{0.9}{\overline{\mD}^\mO}})$, $r\in (0,1)$, 
where $\pi$ is the quotient map 
$\pi: {\scaleobj{0.9}{\mD^\mO}} 
\to {\scaleobj{0.9}{\mD^\mO}}/{\scaleobj{0.9}{\widehat G_\mO}}$. 
Let $\mathcal B_r$ denote the uniform closure of the algebra 
$\mathcal O_w({\scaleobj{0.9}{\mD^\mO}} / {\scaleobj{0.9}{\widehat G_\mO}}) 
|_{{\scaleobj{0.9}{\pi(r{\scaleobj{0.9}{\overline{\mD}^\mO}})}}}$. 
Since ${\scaleobj{0.9}{ \widehat{G}_\mO}}$ 
acts by linear transformations on $\ell^\infty(\mO)$, 
the dilation map $\delta_{\mathfrak{r}},$ $\bbz\mapsto r\bbz,$ 
$\bbz\in {\scaleobj{0.9}{\mD^\mO}},$ 
commutes with the action of ${\scaleobj{0.9}{\widehat{G}_\mO}}$. 
Then there exists a continuous map 
$\delta_{\mathfrak{r};{\scaleobj{0.9}{\widehat{G}_\mO}}} : 
{\scaleobj{0.9}{\overline{\mD}^{\mO}}}\to \pi(r{\scaleobj{0.9}{\overline{\mD}^{\mO}}})$ 
such that the following diagram is commutative:
\begin{center}
\begin{tikzcd}
{\scaleobj{0.9}{\overline{\mD}^\mO }} 
\arrow{r}[name=U]{\delta_{\mathfrak{r}}} 
& r {\scaleobj{0.9}{\overline{\mD}^\mO}}  
\\ 
{\scaleobj{0.9}{\overline{\mD}^\mO}}/{\scaleobj{0.9}{\widehat{G}_\mO}} 
\arrow[leftarrow]{u}{ \pi } 
\arrow{r}[]{\delta_{\mathfrak{r}; {\scaleobj{0.9}{\widehat{G}_\mO}}}}  
& \pi(r{\scaleobj{0.9}{\overline{\mD}^{\mO}}}) 
\arrow[leftarrow,swap]{u}{\pi}
\end{tikzcd}
\end{center}
Moreover, by Lemma \ref{lemma10.1}, 
the pullback of $\delta_{\mathfrak{r};{\scaleobj{0.9}{\widehat{G}_\mO}}}$ 
maps $\mathcal B_r$ isometrically onto the algebra 
$A({\scaleobj{0.9}{\mD^\mO}}/{{\scaleobj{0.9}{\widehat{G}_\mO}}})$. 
This implies
 
\begin{lemma}
\label{lem12.3}
By identifying the points with the evaluation functionals at these points$,$ 
the maximal ideal space $M(\mathcal B_r)$ of the algebra $\mathcal B_r$ 
is coincident with $(\pi(r{\scaleobj{0.9}{\overline{\mD}^\mO}}),\tau_p),$ 
and homeomorphic to 
$({\scaleobj{0.9}{\overline{\mD}^\mO}}/{\scaleobj{0.9}{\widehat{G}_\mO}},\tau_p)$. 
\hfill$\Box$
\end{lemma}
  
Let $P \in 
\text{M}_n(\mathcal O_w({\scaleobj{0.9}{\mD^\mO}}/{\scaleobj{0.9}{\widehat G_\mO}}))$ 
satisfy $P^2=P$. Given $r\in (0,1)$, consider the idempotent 
$P_r:=P|_{{\scaleobj{0.84}{\pi(r{\scaleobj{0.9}{\overline{\mD}^\mO}})}}} 
\in \text{M}_n(\mathcal B_r)$. 
Due to Lemma~\ref{lem12.3} and Theorem~\ref{teo_contractibility}, 
the maximal ideal space $M(\mathcal B_r)$ is contractible, 
and so, by \cite[Thm.~4.1]{BruSas23}, 
there exists $F_r\in \text{GL}_n(\mathcal B_r)$ such that 
\begin{equation}
\label{eq12.47}
F_r^{-1} P_r F_r 
= [\begin{smallmatrix}
I_k & 0\\
0 & 0
\end{smallmatrix}].
\end{equation}
(Note that $k$ does not depend on $r$. Indeed, since the space 
$({\scaleobj{0.9}{\mD^\mO}}/{\scaleobj{0.9}{\widehat G_\mO}},\tau_{hk})$ 
is connected, the set $P({\scaleobj{0.9}{\mD^\mO}}/{\scaleobj{0.9}{\widehat G_\mO}})$ 
belongs to a connected component 
of the complex manifold of idempotents of $\text{M}_n(\mC)$.)  

In what follows, we will use the following notation. 
For a matrix $\alpha \!\in\! \text{M}_n(\mC)$, by $\|\alpha\|_{\text{M}_n(\mC)}$, 
we mean the operator norm of the matrix multiplication map $\alpha: \mC^n\to \mC^n$. 
If $X$ is a compact topological space, 
and $A$ is a Banach subalgebra of $C(X)$ with the supremum norm, 
then for $\alpha \in  \text{M}_n(A)$, we define 
$\|\alpha\|_{\text{M}_n(A)}=\sup_{x\;\!\in\;\! X} \|\alpha(x)\|_{\text{M}_n(\mC)}$.

\begin{lemma}
\label{lem12.4}
Let $0<r_1<r_2<1$. Given $\varepsilon>0,$  
there exists a matrix $D_{r_2,\varepsilon} \in \text{\em GL}_n(\mathcal B_{r_2})$ 
such that for $F_{r_2,\varepsilon}:=F_{r_2}D_{r_2,\varepsilon} 
\in \text{\em GL}_n(\mathcal B_{r_2}),$ 
we have
$$
\|
F_{r_1}^{\pm 1} - F_{r_2,\varepsilon}^{\pm 1}|_{\pi(r_1{\scaleobj{0.9}{\overline{\mD}^{\mO}}})}
\|_{\text{\em M}_n(\mathcal B_{r_1})} < \varepsilon,
$$
and
$$
F_{r_2,\varepsilon}^{-1} P_{r_2} F_{r_2,\varepsilon}
= [\begin{smallmatrix}
I_k & 0\\
0 & 0
\end{smallmatrix}].
$$
\end{lemma}
\begin{proof}
Note  that due to \eqref{eq12.47},
$$
F_{r_1}^{-1}(F_{r_2}|_{\pi(r_1{\scaleobj{0.9}{\overline{\mD}^{\mO}}})})
= [\begin{smallmatrix}
G_{r_1r_2} & 0\\
0 & H_{r_1r_2} 
\end{smallmatrix}]
$$
for some $G_{r_1r_2}\in \text{GL}_k(\mathcal B_{r_1})$, 
$H_{r_1r_2}\in \text{GL}_{n-k}(\mathcal B_{r_1})$. 
Since by Lemma \ref{lem12.3} $M(\mathcal B_r)$ is contractible, 
the Banach Lie groups $\text{GL}_k(\mathcal B_{r_1})$ 
and $\text{GL}_{n-k}(\mathcal B_{r_1})$ are connected. 
Thus, $G_{r_1r_2}$ and $H_{r_1r_2}$ can be joined by paths 
in these groups with their identities. 
This implies that there exist $p,q\in\mN$ 
such that $G_{r_1r_2}=\textstyle\prod_{i=1}^p e^{M_i}$ 
for some $M_i\in \text{M}_k(\mathcal B_{r_1})$, $1\le i\le p$, 
and $H_{r_1r_2}=\textstyle\prod_{i=1}^q e^{N_i}$ 
for some $N_i\in \text{M}_{n-k}(\mathcal B_{r_1})$, $1\le i\le q$. 
By the definition of $\mathcal B_{r_1}$, each such $M_i$ and $N_i$ 
can be uniformly approximated on $\pi(r{\scaleobj{0.9}{\overline{\mD}^\mO}})$ 
by matrices of the same size with entries in 
$\mathcal O_w({\scaleobj{0.9}{\mD^\mO}}/{\scaleobj{0.9}{\widehat G_\mO}})$. 
Thus, by replacing the matrices $M_i$ and $N_i$ in the above finite products 
with their appropriate approximations $M_{i,\varepsilon}$ and $N_{i,\varepsilon}$, 
we obtain the desired matrix $D_{r_2,\varepsilon}$ of the form 
$$
D_{r_2,\varepsilon}
= [\begin{smallmatrix}
G_{r_2,\varepsilon}& 0\\
0 &H_{r_2,\varepsilon}
\end{smallmatrix}],
$$
where $G_{r_2,\varepsilon} := 
\textstyle\prod_{i=1}^p e^{M_{i,\varepsilon}} 
\in \text{GL}_k(\mathcal B_{r_2})$,  
$H_{r_2,\varepsilon} := 
\textstyle\prod_{i=1}^q e^{N_{i,\varepsilon}} 
\in \text{GL}_{n-k}(\mathcal B_{r_2})$, 
such that if $F_{r_2,\varepsilon} := F_{r_2}D_{r_2,\varepsilon}$, then
$$
\|
F_{r_1}^{\pm 1}-F_{r_2,\varepsilon}^{\pm 1}|_{\pi(r_1{\scaleobj{0.9}{\overline{\mD}^{\mO}}})}
\|_{\text{M}_n(\mathcal B_{r_1})} < \varepsilon.
$$
Also,
$$
\begin{array}{rcl}
F_{r_2,\varepsilon}^{-1} P_{r_2} F_{r_2,\varepsilon}
\!\!\!&=&\!\!\!
D_{r_2,\varepsilon}^{-1} F_{r_2}^{-1}P_{r_2} F_{r_2}D_{r_2,\varepsilon}
=
D_{r_2,\varepsilon}^{-1}
[\begin{smallmatrix}
I_k & 0\\
0 & 0
\end{smallmatrix}]
D_{r_2,\varepsilon}
\\[0.15cm]
\!\!\!&=&\!\!\!
[\begin{smallmatrix}
G_{r_2,\varepsilon}^{-1} & 0\\
0 & H_{r_2,\varepsilon}^{-1}
\end{smallmatrix}]
[\begin{smallmatrix}
I_k & 0\\
0 & 0
\end{smallmatrix}]
[\begin{smallmatrix}
G_{r_2,\varepsilon} & 0\\
0 & H_{r_2,\varepsilon}
\end{smallmatrix}]
=
[\begin{smallmatrix}
I_k & 0\\
0 & 0
\end{smallmatrix}].
\end{array}
$$
This completes the proof of the lemma.
\end{proof}

Next, we set $r_m:=\frac{m}{m+1}$ 
and $\varepsilon_m:=\frac{1}{2^m}$, $m\ge 1$. 
Let $F_{m}\in \text{GL}_{n}(\mathcal B_{r_m})$ be such that
$$
F_{m }^{-1} P_{r_m} F_{m}
= [\begin{smallmatrix}
I_k & 0\\
0 & 0
\end{smallmatrix}].
$$
We construct $G_{m}\in \text{GL}_{n}(\mathcal B_{r_{m}})$, $m\in\mN$, 
inductively as follows:
\begin{quote}
$G_1:=F_{1}\in \text{GL}_n(\mathcal B_{r_1})$. 
If $G_1,\dots, G_{m}$ are defined, 
then we apply Lemma~\ref{lem12.4} 
with $F_{r_1}:=G_{m}$, $F_{r_2}:=F_{m+1}$ and $\varepsilon:=\varepsilon_{m+1}$. 
Then $G_{m+1}$ is defined to be the resulting function  $F_{r_2,\varepsilon}$, and so on.
\end{quote}
Then, given $k\!\in\!\mN$,  for all sufficiently large $m_2\!>\!m_1\!>\!k$, 
we have on $\pi(r_{k}{\scaleobj{0.9}{\overline{\mD}^\mO}})$:
$$
\textstyle 
\|G_{m_2}^{\pm 1}-G_{m_1}^{\pm 1}\|_{\text{M}_n(\mathcal B_{r_{k}})}
\le 
\sum\limits_{i=m_{1}}^{m_{2}-1} 
\|G_{i+1}^{\pm 1}-G_{i}^{\pm 1}\|_{\text{M}_n(\mathcal B_{r_{k}})} 
\le \sum\limits_{i=m_{1}}^{m_{2}-1}\varepsilon_{i+1} 
\le {\scaleobj{1.2}{\frac{1}{2^{m_1}}}}.
$$
This shows that for each $k\in\mN$, 
the sequences $\{G_m\}_{m\ge k}$ and $\{G_m^{-1}\}_{m\ge k}$ 
are Cauchy in the  Banach space $\text{M}_n(\mathcal B_{r_{k}})$. 
Thus, the pointwise limits of  these sequences for all $k\in\mN$ 
determine the matrix functions $G$ and $H$ 
on ${\scaleobj{0.9}{\mD^{\mO}}}/{\scaleobj{0.9}{\widehat G_{\mO}}}$. 
Since $G_m\cdot G_m^{-1}=I_n$, we get that $H=G^{-1}$. 
Hence, $G\in 
\text{GL}_n(C({\scaleobj{0.9}{\mD^{\mO}}}/{\scaleobj{0.9}{\widehat G_{\mO}}},\tau_{hk}))$. 
Moreover, since 
$$
G_{m}^{-1} P G_{m}
= [\begin{smallmatrix}
I_k & 0\\
0 & 0
\end{smallmatrix}],
$$
the same is true with $G$ replacing $G_m$. 

Finally, by the definition of the algebras $\mathcal B_{r_m}$, 
every  $\pi^*(G_m)|_{r_m{\scaleobj{0.9}{\mD^\mO}}}$ 
is a continuous G\^{a}teaux holomorphic matrix function on $r_m{\scaleobj{0.9}{\mD^\mO}}$. 
Then, since $\pi^*(G)$ is the uniform limit on compact subsets 
of $({\scaleobj{0.9}{\mD^{\mO}}},\tau_{hk})$ 
of the sequences $\{\pi^*(G_m)\}_{m\ge k}$, $k\in\mN$, 
the matrix function $\pi^*(G)$ 
is continuous G\^{a}teaux holomorphic on ${\scaleobj{0.9}{\mD^{\mO}}}$, 
i.e., it belongs to $\text{GL}_n(\mathcal O_w({\scaleobj{0.9}{\mD^{\mO}}}))$, 
which implies that $G\in 
\text{GL}_n(\mathcal O_{w}({\scaleobj{0.9}{\mD^{\mO}}}/{\scaleobj{0.9}{\widehat G_{\mO}}}))$, 
and  completes the proof of the projective freeness of the algebra $\mathscr O_{b,G}$. 

%-------------------------------------------------------
\subsection[]{Proof of Theorem~\ref{theorem2.7}} 
\label{subsecti12.5}
%-------------------------------------------------------

%.......................
\subsubsection{}
%.......................

In this subsection, we formulate some auxiliary results used in the proof of the theorem.

Let $A$ be an associative algebra  with $1$ over a field of characteristic zero. 
Let ${\rm Exp}_n(A)\subset \text{GL}_n( A)$ 
denote the subgroup generated by exponents of matrices from $\text{M}_n(A)$. 
Clearly, $\text{E}_n(A)\subset {\rm Exp}_n(A)$. 
The {\em exponential rank} $e_n(A)$ of $\text{M}_n(A)$ 
is defined as the least $k\in\mN$ such that 
every matrix in ${\rm Exp}_n(A)$ is a product of $k$ exponents. 
We write $e_n(A)=\infty$ if such a natural number $k$ does not exist. 
 
In our proofs, we use the following result. 

\begin{Th}[\mbox{\cite[Thm.\,2.3]{Phi}}]
Let $n\!\ge\! 2$ and $X$ be any compact manifold 
of dimension at least $m=2(n^2-1)\ell+2$. 
Then $e_n(C(X)) \ge \ell+1$.
\end{Th}

As a  corollary of this result, we obtain

\begin{corollary}\label{cor12.5}
Let $X$ be a reduced Stein space of dimension\footnote{Hereinafter, 
the dimension means the dimension of the complex analytic space.} 
$m\ge  2 (n^2-1)\ell +2,$ where $n\ge 2$. 
Let $S\subset X$ be a set with a nonempty interior $\mathring S$ of dimension $m,$ 
and let $\mathcal O(S)$ be the algebra of continuous functions on $S$ 
that are holomorphic on $\mathring{S}$. 
Then $e_n(\mathcal O(S))\ge\ell$.
\end{corollary}
\begin{proof} We consider two cases. 

\noindent $\bf{1.}$ Suppose that $X=\mC^m$ 
and $S\subset\mC^m$ is a set with nonempty interior 
containing a compact $m$-dimensional cube $K\subset \mR^m$. 

By the above theorem \cite[Thm.~2.3]{Phi}, 
there exists an $F\in {\rm Exp}_n(C(K))$ 
that cannot be represented 
as a product of $\ell$ exponents of matrices in $\text{M}_n(C(K))$. 
Since the set $K$ is contractible, 
there exists $g\in C(K)$ such that $e^{g}=\det F$. 
Thus, replacing $F$ by $e^{-\frac{g}{n}}\cdot F$ if necessary, 
we can assume, without loss of generality, that $F\in \text{SL}_n(C(K))$. 
Also, the contractibility of $K$ implies that 
$F$ can be represented as a finite product of unitriangular matrices 
$F=E_1\cdots E_s$ over $C(K)$ for some $s\in\mN$. 
By the Stone-Weierstrass approximation theorem, 
we can approximate the off-diagonal nonzero entries of these matrices 
to obtain unitriangular matrices $G_1,\dots, G_s\in \text{M}_n(\mC[z_1,\dots, z_m])$ 
(i.e., with holomorphic polynomial entries), 
such that for $G:=G_1\cdots G_s\, (\in {\rm Exp}_n(\mC[z_1,\dots, z_m]))$,
$$
\textstyle
\|FG^{-1}-I_n\|_{\text{M}_n(C(K))} < \frac{1}{2}.
$$
This inequality shows that $FG^{-1}$ has a logarithm in $\text{M}_n(C(K))$ 
given by the formula 
$$
\textstyle
H:=-\sum\limits_{i=1}^\infty{\scaleobj{1.2}{\frac{(I_n-FG^{-1})^i}{i}}}.
$$
Since $F=Ge^{H}$, by our choice of $F$ and because $K\subset S$, 
the holomorphic matrix function $G\in {\rm Exp}_n(\mC[z_1,\dots, z_m])$ restricted to $S$ 
cannot be represented 
as a product of $\ell-1$ exponents of matrices in $\text{ M}_n(C(S))$. 
Obviously the same is true for $\mathcal O(S)\, (\subset C(S))$, 
i.e., $e_n(\mathcal O(S))\ge\ell$.

\smallskip

\noindent ${\bf 2.}$ Now let us consider the general case.

Recall that  by our hypothesis, 
the interior $\mathring S$ of $S$ is of maximal dimension $m$. 
Then the set of regular points of $X$ 
with embedding dimension $m$ contained in $\mathring S$ is nonvoid. 
Let $x\in \mathring S$ be one of these points. 
Since $X$ has embedding dimension $m$ at $x$, 
due to Cartan’s Theorem \cite[Ch.\,IV]{GraRem}, 
there exist holomorphic functions $f_1, \dots , f_m$ on $X$, 
and an open neighbourhood $U\subset \mathring S$ of $x$, 
such that the holomorphic map $\mathscr F = (f_1,\dots,f_m) : X \to \mC^m$ 
is one-to-one on the closure $\overline{U}$ of $U$, maps $x$ to $\mathbf{0}$ 
and $U$ biholomorphically onto the open unit polydisc $\mD^m \subset \mC^m$. 
By our construction, the algebra $\mathcal O(U)$ (of holomorphic functions on $U$) 
is isomorphic by means of the pullback of $\mathscr F|_{U}$ 
to the  algebra $\mathcal O(\mD^m)$. 
Let $K$ be an $m$-dimensional compact cube in $(-1,1)^m=\mD^m\cap\mR^m$, 
and let $G\in \text{Exp}_n(\mC[z_1,\dots, z_m])$ be the matrix, 
as in case $\mathbf{1}$ above, 
that approximates on $K$ some $F\in \text{SL}_n(C(K))$ 
that is not a product of $\ell$ exponents of matrices in $\text{M}_n(C(K))$. 
In particular, $G|_{\mD^n}$ cannot be represented 
as a product of $\ell-1$ exponents of matrices in $\text{M}_n(C(\mD^n))$. 
Then, since $U\subset S$, 
the restriction of the pullback $\mathscr F^*(G)\in \text{SL}_n(\mathcal O(X))$ to $S$ 
cannot be represented as a product 
of less than $\ell$ exponents of matrices in $\text{Exp}_n(\mathcal O(S))$, 
i.e., $e_n(\mathcal O(S))\ge\ell$.
\end{proof}

%.......................
\subsubsection{}
%.......................

In this subsection, we prove the path connectedness 
of the topological groups of Theorem~\ref{theorem2.7}.

Since the algebras $A=\mathscr W_G, \mathscr A_G, \mathscr O_{b,G}$ 
are isometrically isomorphic to algebras 
$\mathcal A
=W({\scaleobj{0.9}{\overline{\mD}^\mO}}/{\scaleobj{0.9}{\widehat{G}_\mO}}), 
 A({\scaleobj{0.9}{\overline{\mD}^\mO}}/{\scaleobj{0.9}{\widehat{G}_\mO}}), 
\mathcal O_w({\scaleobj{0.9}{\mD^\mO}}/{\scaleobj{0.9}{\widehat G_\mO}})$, respectively, 
the corresponding topological groups 
$\text{SL}_n(A)$ and $\text{SL}_n(\mathcal A)$ are isomorphic. 
So we have to prove the path connectedness of groups $\text{SL}_n(\mathcal A)$. 
But ${\scaleobj{0.9}{\overline{\mD}^\mO}}/{\scaleobj{0.9}{\widehat G_\mO}}$ 
and $({\scaleobj{0.9}{\mD^\mO}}/{\scaleobj{0.9}{\widehat G_\mO}},\tau_{hk})$ 
are contractible spaces. Thus the transpose of a homotopy 
between the identity map and a constant map of value $x_*$ (i.e., a contractibility map) 
applied to each entry of $F\in \text{SL}_n(\mathcal A)$ 
determines a path in $\text{SL}_n(\mathcal A)$ 
that connects $F$ with the matrix function of the constant value $F(x_*)$. 
Thus, since the group $\text{SL}_n(\mC)$ is path connected, 
each element of $\text{SL}_n(\mathcal A)$ 
can be connected by a path in $\text{SL}_n(\mathcal A)$ with $I_n$, as required.

%.......................
\subsubsection{}
%.......................

In this subsection, we prove inequality \eqref{equ2.16}.

First,  we prove inequality \eqref{equ2.16} for $|\mO|<\infty$. 
We use results of Section~\ref{subsekt2.3}. 

For $|\mO|<\infty$, the algebra $\mathscr O_{b,G}$ 
is isometrically isomorphic to the algebra $\mathcal O(\Omega)$ 
of holomorphic functions on a Stein domain $\Omega\Subset\mC^{|\mO|}$. 
Then due to \cite{IvaKut} (see also \cite{IvaKut2}), 
the connectedness of $\text{SL}_n(\mathscr O_{b,G})$, 
and a remark in \cite[Sect.\,1.1.]{Bru18}, 
$\text{E}_n(\mathscr O_{b,G})=\text{SL}_n(\mathscr O_{b,G})$, 
and $t_n(\mathscr O_{b,G})\le v(2|\mO|)+3$. 
In this case, the left-hand side inequality in \eqref{equ2.16}  
follows from Corollary~\ref{cor12.5} with $S=X:=\Omega$ 
by an  analog of \cite[Lm.\,2.1]{Bru22} (with the same proof), 
which asserts 
in the case of an algebra $A\in\{\mathscr W_G, \mathscr A_G, \mathscr O_{b,G}\}$ 
that 
\begin{equation}
\label{equ12.51}
\textstyle 
e_n(A) 
\le \left\lfloor {\scaleobj{1.2}{\frac{t_n(A)}{2}}} \right\rfloor +1.
\end{equation}
Thus, we get for $n\ge 2$,   
$$
\textstyle
t_n(\mathscr O_{b,G})
\ge 2\;\!e_n(\mathscr O_{b,G})\!-\!2
= 2\;\!e_n(\mathcal O(\Omega))\!-\!2
\ge {\scaleobj{1.2}{\frac{|\mO|-2}{n^2-1}}} \!-\!2
\ge \left\lfloor{\scaleobj{1.2}{\frac{|\mO|}{n^2-1}}}\right\rfloor\!-\!3. 
$$
As $t_n(\mathscr O_{b,G})\ge 4$ for all $n\ge 2$, 
the latter gives the left-hand side inequality in \eqref{equ2.16} in this case.

Next, we consider the algebras $\mathscr W_G$ and $\mathscr A_G$ 
(isometrically isomorphic to algebras 
$W({\scaleobj{0.9}{\overline{\mD}^\mO}}/{\scaleobj{0.9}{\widehat{G}_\mO}})$ 
and $A({\scaleobj{0.9}{\overline{\mD}^\mO}}/{\scaleobj{0.9}{\widehat{G}_\mO}})$, respectively) 
when $|\mO|<\infty$. In this case, according to results of Section~\ref{subsekt2.3}, 
$\mathscr A_G$  is isometrically isomorphic to the algebra $A(\Omega)$, 
that is, the uniform closure of the restriction 
to the above Stein domain $\Omega\Subset\mC^{|\mO|}$ 
of the algebra of polynomials $\mC[z_1,\dots, z_{|\mO|}]$. 
Moreover, the maximal ideal space $M(\mathscr A_G)$ 
is naturally identified with the closure $\overline{\Omega}$ of $\Omega$. 
Since the covering dimension of $\overline{\Omega}$ is $2|\Omega|$, 
due to \cite[Thm.\,1.1]{Bru18}, $t_n(\mathscr A_G)\le v(2|\mO|)+5$. 
The lower bound is the same 
as in the case of the algebra $\mathscr O_{b,G}$ 
(to obtain it, we use Corollary~\ref{cor12.5} with $X=\mC^{|\mO|}$ 
and $S=\overline{\Omega}$, and inequality \eqref{equ12.51}).  

In turn, the algebra $\mathscr W_G$ is isometrically isomorphic 
to a dense subalgebra of $A(\Omega)$, 
and its maximal ideal space is homeomorphic to $\overline{\Omega}$. 
Hence, the proof of inequality \eqref{equ2.16} in this case 
is the same as that for the algebra $\mathscr A_G$.

Now, we prove \eqref{equ2.16} for $|\mO|=\infty$, 
i.e., we prove that in this case 
for an algebra $\mA_G\in\{\mathscr W_G, \mathscr A_G, \mathscr O_{b,G}\}$, 
$t_n(\mA_G)=\infty$ for each $n\ge 2$. 

We use the notation of Corollary~\ref{cor12.2}. 
Thus, $\mO_m:={\scaleobj{0.9}{\bigcup_{i=1}^m}\;\!} O_i$, 
where $O_1, O_2, O_3,\dots$ are the (finite) orbits 
of the action of  $G$ on $\mO\subset\mN$, 
and $\mathcal P_{\mO_m}:=\{p_i\in\mP,\ i\in\mO_m\}$, $m\in\mN$. 
Also, $\langle \mathcal P_{\mO_m}\rangle\subset\mN$ 
is the unital multiplicative semigroup 
generated by the elements of $\mathcal P_{\mO_m}$, 
and $\mA_{\mathcal P_{\mO_m}}$ is a subalgebra 
of functions of $\mA \in \{\mathscr W, \mathscr A, \mathscr O_{b}\}$ 
with Dirichlet series of the form 
$\sum_{k\;\! \in\;\! \langle\mathcal P_{\mO_m}\rangle} a_k k^{-s}$. 
In turn, $\mA_{G,\mathcal P_{\mO_m}} = \mA_G\cap \mA_{\mathcal P_{\mO_m}}$ 
is the algebra of limit functions of Dirichlet series 
that are invariant with respect to the action $\tilde S$ 
of the finite group $G/G_{\mO_m}$ on $\mA_{\mathcal P_{\mO_m}}$, 
where $G_{\mO_m}$ is the normal subgroup of $G$ 
consisting of all $\sigma\in G$ acting identically on  $\mO_m$, 
see Section~\ref{subsect6.2}. 
Therefore we can apply the previous case (of the finite set of orbits) 
to $\mA_{G,\mathcal P_{\mO_m}}$ 
with $G/G_{\mO_m}$ in place of $G$, and $\mO_m$ in place of $\mO$. 
In particular, in this setting, as in the proof of the previous case, 
we obtain for $n\ge 2$,  
\begin{equation}
\label{equ12.52}
\textstyle
t_n(\mA_{G,\mathcal P_{\mO_m}})
\ge  2\;\!e_n(\mA_{G,\mathcal P_{\mO_m}})\!-\!2
\ge  \left\lfloor{\scaleobj{1.2}{\frac{|\mO_m|}{n^2-1}}}\right\rfloor\!-\!3. 
\end{equation}
Assume, on the contrary, that $t_n(\mA_G)=:t<\infty$ for some $n\ge 2$. 
Since $|\mO|\!=\!\infty$, the set of orbits $\{O_i\}$ of the action of $G$ on $\mN$ is infinite, 
and so there is some $\mO_m$ 
such that $\bigl\lfloor\frac{|\mO_m|}{n^2-1}\bigr\rfloor-3 > 2t$. 
Therefore according to \eqref{equ12.52}, 
$t_n(\mA_{G,\mathcal P_{\mO_m}})
\ge 2\;\!e_n(\mA_{G,\mathcal P_{\mO_m}})-2 
> 2t$. Let $F\in \text{Exp}_n(\mA_{G,\mathcal P_{\mO_m}})$ 
be an element that cannot be represented 
as a product of less than $e_n(\mA_{G,\mathcal P_{\mO_m}})\,(\ge t\!+\!2)$ 
exponents of matrices in $\text{M}_n(\mA_{G,\mathcal P_{\mO_m}})$. 
However, since 
${\rm SL}_n(\mA_{G,\mathcal P_{\mO_m}}) 
= \text{E}_n(\mA_{G,\mathcal P_{\mO_m}}) 
\subset \text{E}_n(\mA_G)$, by our assumption, we have 
$F=e^{F_1}\cdots e^{F_t}$ for some $F_i\in \text{M}_n(\mA_G)$, $1\le i\le t$. 
(We have used the fact that each unitriangular matrix in $\text{M}_n(\mA_G)$ 
is an exponent of some other matrix in $\text{M}_n(\mA_G)$.) 
Applying the continuous homomorphism 
$\Phi_{\mathcal P_{\mO_m}}:\mA\to \mA_{G,\mathcal P_{\mO_m}}$ 
of Corollary \ref{cor12.2} to this equation, 
we get $F=e^{H_1}\cdots e^{H_t}$, 
where $H_i:=\Phi_{\mathcal P_{\mO_m}}\!(F_i) 
\in \text{M}_n(\mA_{G,\mathcal P_{\mO_m}})$. 
This contradicts the choice of $F$, 
and shows that $t_n(\mA_G)=\infty$ for each $n\ge 2$.

\begin{remark}
\label{rem12.6}
The above argument also shows that if $|\mO|=\infty$, 
then for an algebra 
$\mA_G\in\{\mathscr W_G, \mathscr A_G, \mathscr O_{b,G}\}$,  
$e_n(\mA_G)=\infty$  for each $n\ge 2$.
\end{remark}

%.......................
\subsubsection{}
%.......................

In this subsection we prove that if $|\mO|=\infty$, 
then we have 
${\scaleobj{0.98}{\text{E}_n(\mathscr O_{b,G}) 
\!\subsetneq\! 
\text{SL}_n(\mathscr O_{b,G})}}$ 
(equivalently, 
${\scaleobj{0.98}{\text{E}_n( 
\mathcal O_w({\scaleobj{0.9}{\mD^{\mO}}}\!/{\scaleobj{0.9}{\widehat G_{\mO}}}))
\!\subsetneq\!
\text{SL}_n(
\mathcal O_w({\scaleobj{0.9}{\mD^{\mO}}}\!/{\scaleobj{0.9}{\widehat G_{\mO}}}))}}$) 
for each $n\ge 2$.

First, let us consider the function 
$f\in \mathcal O_w({\scaleobj{0.9}{\mD^{\mO}}}/{\scaleobj{0.9}{\widehat G_{\mO}}})$ 
such that 
$$
\textstyle
(\pi^*(f))(\bbz) := {\scaleobj{1.2}{\frac{1}{|O_1|}}}\sum\limits_{i\;\! \in \;\! O_1}z_i.
$$
(As before, 
$\pi:{\scaleobj{0.9}{\mD^{\mO}}} 
\to {\scaleobj{0.9}{\mD^{\mO}}}/{\scaleobj{0.9}{\widehat G_{\mO}}}$ 
is the quotient map, and the coordinates of $\bbz\in\ell^\infty(\mO)$ 
are denoted by $z_i$, $i\in\mO$.) 
 
\begin{lemma}
\label{lem12.8}
$f$ maps each $\pi(r{\scaleobj{0.9}{\overline{\mD}^{\mO}}})$ 
onto $r\overline{\mD},$ $r\in (0,1)$.
\end{lemma}
\begin{proof}
Since $|z_i|\le 1$, $i\in O_1$, 
we have $f(\pi(r{\scaleobj{0.9}{\overline{\mD}^{\mO}}}))\subset r\overline{\mD}$. 
These sets  coincide, 
because if $\bf{r}\in \ell^\infty(\mO)$ has all its coordinates equal to $r$, 
then the image of the disc 
$\{z\cdot {\bf r}\,:\, z\in\overline{\mD}\}\subset r{\scaleobj{0.9}{\overline{\mD}^{\mO}}}$ 
under $\pi^*(f)$ is $r\overline{\mD}$.  
\end{proof}

Fix  a sequence $\{r_i\}_{i\in\mN}\subset (0,1)$, strictly increasing to $1$, 
and define $s_i:=  \frac{r_{i+1}-r_{i}}{3}$, $i\in\mN$. 
Consider open discs $D_i\subset (r_{i+1}\mD)\setminus (r_i\overline{\mD})$ 
centred at points $c_i:=\frac{r_i+r_{i+1}}{2}$ of radii $s_i$, $i\in\mN$. 
By definition, $r_i\overline{\mD}\cap\overline {D}_i=\emptyset$ for every $i$. 
Let $\chi_i$ be the function which is equal to $1$ on $r_i\overline{\mD}$, 
and $0$ on $\overline{D}_i$.

\begin{lemma}
\label{lem12.9}
Given $\varepsilon>0$ there exists a polynomial $p_{i,\varepsilon}\in\mC[z]$ such that
$$
\max_{(r_i\overline{\mD})\sqcup\overline {D}_i}|\chi_i-p_{i,\varepsilon}| < \varepsilon.
$$
\end{lemma}
\begin{proof}
The set $K:=(r_i\overline{\mD})\sqcup\overline {D}_i$ 
is polynomially convex because its complement is a connected subset of $\mC$ 
(see, e.g., \cite[Lm.~1.3]{Gam}). 
Hence, as $\chi_i$ extends holomorphically to an open neighbourhood of $K$, 
by Runge's approximation theorem (see, e.g., \cite[Cor.~1.4]{Gam}), 
$\chi_i$ is uniformly approximated on $K$ by holomorphic polynomials.
\end{proof}

Next, consider the closed ball 
${\scaleobj{0.9}{\overline{\mD}_{s_i}}}({\bf c}_i) 
:= {\bf c}_i+s_i{\scaleobj{0.9}{\overline{\mD}^{\mO}}} 
\subset 
r_{i+1}{\scaleobj{0.9}{\overline{\mD}^{\mO}}}$, 
where ${\bf c}_i\in \ell^\infty(\mO)$ has all coordinates equal to $c_i$. 
Also, let ${\scaleobj{0.9}{\mD_{s_i}}}({\bf c}_i) 
:= {\bf c}_i+s_i {\scaleobj{0.9}{\mD^{\mO}}}$. 
We set $V_i := \pi({\scaleobj{0.9}{\overline{\mD}_{s_i}}}({\bf c}_i))$. 
Clearly, $f(V_i)=\overline{D}_i$, 
and thus by the definition of $D_i$ and Lemma~\ref{lem12.8},  
$V_i\cap \pi(r_i{\scaleobj{0.9}{\overline{\mD}^\mO}}) = \emptyset$. 
Let $A_i$ be the uniform closure of the algebra 
$\mathcal O_w({\scaleobj{0.9}{\mD^{\mO}}}/{\scaleobj{0.9}{\widehat G_{\mO}}})|_{V_i}$.

\begin{lemma}
\label{lemma12.9}
The algebra $A_i$ is isometrically isomorphic 
to the algebra $A({\scaleobj{0.9}{\mD^{\mO}}}/{\scaleobj{0.9}{\widehat G_{\mO}}})$.
\end{lemma}
\begin{proof}
Note that the ball ${\scaleobj{0.9}{\overline{\mD}_{s_i}}}({\bf c}_i)$ 
is invariant with respect to the action of the group 
${\scaleobj{0.9}{\widehat G_{\mO}}}$ 
on ${\scaleobj{0.9}{\overline{\mD}^{\mO}}}$ (by permutations of the coordinates). 
In turn, the map 
$H_i: {\scaleobj{0.9}{\overline{\mD}_{s_i}}}({\bf c}_i) 
\to {\scaleobj{0.9}{\overline{\mD}^{\mO}}}$, 
$H_i(\bbz):=\frac{1}{s_i}(\bbz -{\bf c}_i)$, 
$\bbz\in {\scaleobj{0.9}{\overline{\mD}_{s_i}}}({\bf c}_i)$, 
is equivariant with respect to the action of ${\scaleobj{0.9}{\widehat G_{\mO}}}$, 
and its pullback $H_i^*$ maps 
the algebra $\mC[\bbz]$ of polynomials in the coordinates $z_i$, $i\in\mO$, 
isomorphically onto itself. Hence, $H_i^*$ is an isometric isomorphism 
from the uniform closure  (namely, $A({\scaleobj{0.9}{\mD^\mO}})$) 
of the algebra $\mC[\bbz]|_{{\scaleobj{0.9}{\mD^{\mO}}}}$ 
to the uniform closure (denoted by $A({\scaleobj{0.9}{\mD_{s_i}}}({\bf c}_i))$) 
of the algebra $\mC[\bbz]|_{{\scaleobj{0.9}{\mD_{s_i}}}({\bf c}_i)}$. 
Since the algebra $\mC[\bbz]$ 
is dense in $\mathcal O_w({\scaleobj{0.9}{\mD^\mO}})$ 
and ${\scaleobj{0.9}{\overline{\mD}_{s_i}}}({\bf c}_i) 
\subset r_{i+1}{\scaleobj{0.9}{\overline{\mD}^{\mO}}}$, 
the algebra $A({\scaleobj{0.9}{\mD_{s_i}}}({\bf c}_i))$ 
is also the uniform closure of 
$\mathcal O_w({\scaleobj{0.9}{\mD^\mO}})|_{{\scaleobj{0.9}{\mD_{s_i}}}({\bf c}_i)}$. 
Since $H_i$ is equivariant with respect to 
the action of ${\scaleobj{0.9}{\widehat G_{\mO}}}$, 
$H_i^*$ maps the ${\scaleobj{0.9}{\widehat G_\mO}}$-invariant subalgebra 
$A({\scaleobj{0.9}{\mD^\mO}})_{{\scaleobj{0.9}{\widehat G_\mO}}}$ 
of $A({\scaleobj{0.9}{\mD^\mO}})$ 
isometrically onto the ${\scaleobj{0.9}{\widehat G_\mO}}$-invariant subalgebra 
$A({\scaleobj{0.9}{\mD_{s_i}}}({\bf c}_i))_{{\scaleobj{0.9}{\widehat G_\mO}}}$ 
of $A({\scaleobj{0.9}{\mD_{s_i}}}({\bf c}_i))$. 
Finally, since the algebra 
$\mathcal O_w({\scaleobj{0.9}{\mD^\mO}})|_{{\scaleobj{0.9}{\mD_{s_i}}}({\bf c}_i)}$ 
is dense in $A({\scaleobj{0.9}{\mD_{s_i}}}({\bf c}_i))$, 
applying the projection $P_{{\scaleobj{0.9}{\widehat G_\mO}}}$ given by formula \ref{proj}, 
we get that $A({\scaleobj{0.9}{\mD_{s_i}}}({\bf c}_i))_{{\scaleobj{0.9}{\widehat G_\mO}}}$ 
is the uniform closure of the algebra 
$\mathcal O_w({\scaleobj{0.9}{\mD^\mO}})_{{\scaleobj{0.9}{\widehat{G}_\mO}}}$. 
Thus, as $\mathcal O_w({\scaleobj{0.9}{\mD^\mO}})_{{\scaleobj{0.9}{\widehat{G}_\mO}}}$ 
is isometrically isomorphic to 
$\mathcal O_w({\scaleobj{0.9}{\mD^{\mO}}}/{\scaleobj{0.9}{\widehat G_{\mO}}})$, 
the algebra $A({\scaleobj{0.9}{\mD_{s_i}}}({\bf c}_i))_{{\scaleobj{0.9}{\widehat G_\mO}}}$ 
is isometrically isomorphic to $A_i$. 
This implies that $H_i^*$ determines an isometric isomorphism 
from $A({\scaleobj{0.9}{\mD^{\mO}}}/{\scaleobj{0.9}{\widehat G_{\mO}}})$ 
(which is isometrically isomorphic to 
$A({\scaleobj{0.9}{\mD^\mO}})_{{\scaleobj{0.9}{\widehat G_\mO}}}$) to $A_i$, 
and completes the proof. 
\end{proof}

Since for $|\mO|=\infty$, 
$e_n(A({\scaleobj{0.9}{\mD^{\mO}}}/{\scaleobj{0.9}{\widehat G_{\mO}}})) 
= \infty$ for each $n\ge 2$, 
Lemma~\ref{lemma12.9} implies that there exist $F_i \in 
\text{SL}_n(\mathcal O_w({\scaleobj{0.9}{\mD^{\mO}}}/{\scaleobj{0.9}{\widehat G_{\mO}}}))$ 
and $\ell_i\in\mN$ such that $\lim\limits_{i\to\infty} \ell_i=\infty$, 
and $F_i|_{V_i}$ cannot be represented 
as a product of $\le\ell_i$ exponents of functions from $\text{M}_n(A_i)$. 

We  inductively construct a sequence $G_i \in 
\text{SL}_n( \mathcal O_w({\scaleobj{0.9}{\mD^{\mO}}}/{\scaleobj{0.9}{\widehat G_{\mO}}}) )$, $i\in\mN$, satisfying the following properties:
$$
\begin{array}{cl}
\text{(i)}   & G_{1}:=F_1,\\[0.1cm]
\text{(ii)}  & \|G_{i}-F_{i}\|_{\text{M}_n(C(V_i))}
\le {\scaleobj{1.1}{\frac{1}{2^{i}\|F_{i}^{-1}\|_{\text{M}_n(C(V_i))}}}},\\
\text{(iii)} & \|G_{i}^{\pm 1}-G_{i-1}^{\pm 1}\|_{\text{M}_n(\mathcal B_{r_{i}})}
\le {\scaleobj{1.1}{\frac{1}{2^{i}\max\limits_{1\le j\le i-1}\|F_j^{-1}\|_{\text{M}_n(C(V_j))}}}}.
\end{array}
$$
In (iii) above, recall that for $r\in (0,1)$, 
$\mathcal B_r$ denotes the uniform closure of the algebra 
$\mathcal O_w({\scaleobj{0.9}{\mD^\mO}} / 
{\scaleobj{0.9}{\widehat G_\mO}})|_{\pi(r{\scaleobj{0.9}{\overline{\mD}^\mO}})}$. 
 
If $G_1,\dots, G_{i-1}$, $i\ge 2$, have been constructed, 
then we construct $G_{i}$ as follows.

Since the set $\pi(r_{i+1}{\scaleobj{0.9}{\overline{\mD}^\mO}})$ is contractible, 
$F_i$ can be represented on this set 
as a finite product of unitriangular matrices with entries in $\mathcal B_{r_{i+1}}$: 
\begin{equation}
\label{eq12.49}
F_i = (I_n+F_{i,1})\cdots (I_n+F_{i,k_i}).
\end{equation}
(Here all $F_{i,j}$ are nilpotent, i.e., $F_{i,j}^n=0$.)
 
Similarly, we can represent $G_{i-1}$ on $\pi(r_{i}{\scaleobj{0.9}{\overline{\mD}^\mO}})$ 
as a finite product of  unitriangular matrices with entries in $\mathcal B_{r_i}$:
\begin{equation}
\label{eq12.50}
G_{i-1} = (I_n+G_{i-1,1})\cdots (I_n+G_{i-1,s_{i-1}}).
\end{equation}
Take $\chi_i$ and $p_{i,\varepsilon}$ as in Lemma~\ref{lem12.9} 
for some $\varepsilon>0$ which will be specified later. 
By our definition, see Lemma~\ref{lem12.8}, 
$f^*(\chi_i)$ is well-defined on $\pi(r_i{\scaleobj{0.9}{\overline{\mD}^\mO}})$ 
where it takes the value $1$, and on $V_i$ where it takes the value $0$. 
Moreover, $f^{*}(p_{i,\varepsilon}) \in 
\mathcal O_w({\scaleobj{0.9}{\mD^{\mO}}}/{\scaleobj{0.9}{\widehat G_{\mO}}})$, 
and satisfies
\begin{equation}
\label{eq12.51}
\max_{\pi(r_i\overline{\mD})\sqcup V_i}|f^*(\chi_i)-f^*(p_{i,\varepsilon})| < \varepsilon.
\end{equation}
By the definition of the algebras $\mathcal B_r$, 
the nonzero entries of matrix functions $F_{i,j}$ and $G_{i-1,s}$ 
can be uniformly approximated on $\pi(r_{i+1}{\scaleobj{0.9}{\overline{\mD}^\mO}})$ 
and $\pi(r_{i}{\scaleobj{0.9}{\overline{\mD}^\mO}})$ 
by functions belonging to 
$\mathcal O_w({\scaleobj{0.9}{\mD^{\mO}}}/{\scaleobj{0.9}{\widehat G_{\mO}}})$. 
We denote the approximating nilpotent matrices by 
$F_{i,j,\varepsilon }$ and $G_{i-1,s,\varepsilon }$, 
and assume that they satisfy for all $j$ and $s$ the following inequalities:
$$
\begin{array}{l}
\|F_{i,j,\varepsilon}-F_{i,j}\|_{\text{M}_n(\mathcal B_{r_{i+1}})}
< \varepsilon, \;\text{ and}
\\[0.1cm]
\|G_{i-1,s,\varepsilon}-G_{i-1,s}\|_{\text{M}_n(\mathcal B_{r_{i}})} 
< \varepsilon.
\end{array}
$$
Then we set
\begin{eqnarray}
\label{eq2.52}
F_{i,\varepsilon} 
\!\!\!\!\!&=&\!\!\!\! 
(I_n\!+\!(1\!-\!f^*(p_{i,\varepsilon}))F_{i,1,\varepsilon}) 
\cdots 
(I_n\!+\!(1\!-\!f^*(p_{i,\varepsilon}))F_{i,k_i,\varepsilon}) 
\quad\quad\quad \;\;
\\
\label{eq12.53}
G_{i-1,\varepsilon} 
\!\!\!\!\!&=&\!\!\!\! 
(I_n\!+\!f^*(p_{i,\varepsilon})G_{i-1,1,\varepsilon}) 
\cdots 
(I_n\!+\!f^*(p_{i,\varepsilon})G_{i-1,s_{i-1},\varepsilon}).
\quad \quad\quad  \;\;
\end{eqnarray}
By their definition, the matrix functions 
$F_{i,\varepsilon}\cdot G_{i-1,\varepsilon} \in 
\text{SL}_n(\mathcal O_w({\scaleobj{0.9}{\mD^{\mO}}}/{\scaleobj{0.9}{\widehat G_{\mO}}}))$, 
$\varepsilon>0$, 
approximate $F_i$ on $V_i\, (\subset\pi(r_{i+1}{\scaleobj{0.9}{\overline{\mD}^\mO}}))$ 
and $G_{i-1}$ on $\pi(r_i{\scaleobj{0.9}{\overline{\mD}^\mO}})$, 
and converge uniformly to them on these sets as $\varepsilon$ tends to $0$. 
Thus, we can choose $\varepsilon:=\varepsilon_i$ small enough 
so that for $G_i:=F_{i,\varepsilon_i}\cdot G_{i-1,\varepsilon_i}$, 
$$
\begin{array}{l}
\|G_{i}-F_{i}\|_{\text{M}_n(C(V_i))}
\le
{\scaleobj{1.2}{\frac{1}{2^{i}\|F_{i}^{-1}\|_{\text{M}_n(C(V_i))}}}},
\\[0.33cm]
\|G_{i}^{\pm 1}-G_{i-1}^{\pm 1}\|_{\text{M}_n(\mathcal B_{r_{i}})}
\le
{\scaleobj{1.2}{\frac{1}{2^{i}\max\limits_{1\le j\le i-1}\|F_j^{-1}\|_{\text{M}_n(C(V_j))}}}},
\end{array}
$$
as required.

Now, by the definition, given $k\!\in\!\mN$, for all sufficiently large $i_2\!>\!i_1\!>\!k$, 
we have on $\pi(r_{k}{\scaleobj{0.9}{\overline{\mD}^\mO}})$:
$$
\textstyle 
\|G_{i_2}^{\pm 1}-G_{i_1}^{\pm 1}\|_{M_n(\mathcal B_{r_{k}})}
\le 
\sum\limits_{i=i_{1}}^{i_{2}-1} \|G_{i+1}^{\pm 1}-G_{i}^{\pm 1}\|_{M_n(\mathcal B_{r_{k}})}
\le 
\sum\limits_{i=i_{1}}^{i_{2}-1} {\scaleobj{1.1}{\frac{1}{2^{i+1} }}}
\le 
{\scaleobj{1.1}{\frac{1}{2^{i_1}}}}.
$$
(Note that as $F_i\!\in\! 
\text{SL}_n(\mathcal O_w({\scaleobj{0.9}{\mD^{\mO}}}/{\scaleobj{0.9}{\widehat G_{\mO}}}))$, 
$\|F_{i}^{-1}\|_{\text{M}_n(C(V_i))}\!\ge\! 1$ for all $i$.) 
This shows that the sequences $\{G_i\}_{i\in\mN}$ and $\{G_i^{-1}\}_{i\in\mN}$ 
restricted to $\pi(r_{k}{\scaleobj{0.9}{\overline{\mD}^\mO}})$ 
are Cauchy sequences in the Banach space $\text{M}_n(\mathcal B_{r_{k}})$. 
Hence, the pointwise limits of these sequences 
determine the matrix functions $G$ and $H$ 
on ${\scaleobj{0.9}{\mD^{\mO}}}/{\scaleobj{0.9}{\widehat G_{\mO}}}$. 
Since $G_i\cdot  G_i^{-1}\!=\!I_n$, we get that $H=G^{-1}$. Consequently, 
$G\in \text{SL}_n(C({\scaleobj{0.9}{\mD^{\mO}}}/{\scaleobj{0.9}{\widehat G_{\mO}}},\tau_{hk}))$.

Also, we have
$$
\!\!\! 
\begin{array}{rcl} 
\|G\!-\!F_{i}\|_{\text{M}_n(C(V_i))}
\!\!\!&\le&\!\!\!
{\scaleobj{0.9}{\sum\limits_{m=i}^{\infty}}} 
\|G_{m+1}-G_{m}\|_{\text{M}_n(C(V_i))} 
+ \|G_{i}-F_{i}\|_{\text{M}_n(C(V_i))} 
\\[0.3cm] 
\!\!\!&\le&\!\!\! 
{\scaleobj{0.9}{\sum\limits_{m=i}^\infty}} 
{\scaleobj{1.1}{ \frac{1}{2^{m+1} 
\max\limits_{1\le j\le m} \|F_j^{-1}\|_{\text{M}_n(C(V_j))}}}}
+ {\scaleobj{1.1}{\frac{1}{2^{i}\|F_{i}^{-1}\|_{\text{M}_n(C(V_i))}}}} 
\\[0.3cm] 
\!\!\!&\le&\!\!\! 
\Big(\;\!{\scaleobj{0.9}{\sum\limits_{m=i}^\infty}} 
{\scaleobj{1.1}{\frac{1}{2^{m+1}}}} 
+{\scaleobj{1.1}{\frac{1}{2^{i}}}}\Big) 
{\scaleobj{1.1}{\frac{1}{\|F_{i}^{-1}\|_{\text{M}_n(C(V_i))}}}} 
= {\scaleobj{1.1}{\frac{1}{2^{i-1}\|F_{i}^{-1}\|_{\text{M}_n(C(V_i))}}}}.\;
\end{array}
$$
Thus, we obtain
$$
\textstyle
\|G\cdot F_{i}^{-1}-I_n\|_{\text{M}_n(C(V_i))} 
\le \|G-F_{i}\|_{\text{M}_n(C(V_i))}\cdot\|F_{i}^{-1}\|_{\text{M}_n(C(V_i))}
\le {\scaleobj{1.1}{\frac{1}{2^{i-1}}}}.
$$
This implies that for $i\ge 2$, 
there exists a matrix $H_i\in \text{M}_n(A_i)$ with zero trace 
such that $e^{H_i}=(G\cdot F_{i}^{-1})|_{V_i}$. 
Then by our choice of $F_i$, the matrix function $G|_{V_i}$ 
cannot be represented as a product of $\le\ell_i-1$ exponents 
of functions from $\text{M}_n(A_i)$.

Finally, assume on the contrary that $G\in 
\text{E}_n(\mathcal O_w({\scaleobj{0.9}{\mD^{\mO}}}/{\scaleobj{0.9}{\widehat G_{\mO}}}))$. 
This implies that there is some $\ell\in\mN$ 
and $L_1,\dots, L_\ell\in 
\text{M}_n(\mathcal O_w({\scaleobj{0.9}{\mD^{\mO}}}/{\scaleobj{0.9}{\widehat G_{\mO}}}))$ 
such that $G=e^{L_1}\cdots e^{L_\ell}$, which contradicts the above statement.
 
The proof of Theorem \ref{theorem2.7} is complete.

\end{document}